\documentclass{amsart}
\usepackage[utf8]{inputenc}
\usepackage{amsfonts}
\usepackage{amsmath}
\usepackage{amssymb}
\usepackage{amsthm}
\usepackage{mathtools}
\usepackage{enumerate}
\usepackage{graphicx}
\usepackage{hyperref}
\hypersetup{hidelinks,
 pdftitle={How the p-Bergman metric approaches the Carath\'eodory metric},
 pdfauthor={Gunhee Cho}}
\usepackage[capitalise]{cleveref}
\usepackage[top=3cm,bottom=3cm,left=2.5cm,right=2.5cm,headsep=10pt,letterpaper]{geometry}
\usepackage{mathrsfs}

\usepackage[
    backend=biber,
    style=alphabetic,
    sorting=nyt,
    url=false,
    isbn=false
]{biblatex}
\numberwithin{equation}{section}

\theoremstyle{plain}
\newtheorem{theorem}{Theorem}[section]
\newtheorem{proposition}[theorem]{Proposition}
\newtheorem{lemma}[theorem]{Lemma}
\newtheorem{corollary}[theorem]{Corollary}
\newtheorem{maintheorem}{Theorem}

\crefname{maintheorem}{Theorem}{Theorems}

\theoremstyle{definition}
\newtheorem{definition}[theorem]{Definition}
\newtheorem{example}[theorem]{Example}

\theoremstyle{remark}
\newtheorem{remark}[theorem]{Remark}

\newcommand{\C}{\mathbb{C}}
\newcommand{\R}{\mathbb{R}}
\newcommand{\eps}{\varepsilon}
\newcommand{\Id}{\operatorname{Id}}
\newcommand{\mob}{\mathfrak m}

\usepackage[expansion=false]{microtype}
\usepackage{needspace}
\AtBeginEnvironment{theorem}{\Needspace{5\baselineskip}}
\AtBeginEnvironment{proposition}{\Needspace{4\baselineskip}}
\AtBeginEnvironment{corollary}{\Needspace{4\baselineskip}}
\title{How the $p$-Bergman metric approaches the Carath\'eodory metric}
\author[G. Cho]{Gunhee Cho}
\address{Department of Mathematics, Texas State University, San Marcos, TX 78666, USA}
\email{wvx17@txstate.edu}
\date{5 November 2026}
\begin{document}
\begin{abstract}
We study the rate at which the $p$-Bergman metric approaches the
Carath\'eodory metric on convex domains with measures comparable to
$\delta^\alpha dV$, $\alpha>-1$.
On smooth convex domains of finite type the logarithmic coefficient
is bounded above by the harmonic-measure average of the tangential
multitype dimension plus $1+\alpha$, with equality when almost every
endpoint of the extremal disc is coercive.
For a class of convex domains in $\C^2$ with $C^4$ boundary and
rank-one real tangential contact along an axis disc, a winding number
distinguishes the coefficient $2+\alpha$ from an upper bound
$7/4+\alpha$; quartic contact makes the latter bound exact.
The exact value holds along the entire axis of explicit
real-algebraic strongly pseudoconvex domains, including perturbations
without the weighted circle action, although their pointwise boundary
models agree up to unitary equivalence with those of a domain having
coefficient $2+\alpha$.
A $C^2$ example has negative winding but coefficient $2+\alpha$, with
a $-\frac12\log\log p$ correction, so boundary regularity is essential
to the gap theorem.
\end{abstract}
\maketitle
\section{Introduction}
Let $D\Subset\C^n$ and let $\mu$ be a probability measure locally
comparable to Lebesgue measure. For $p\ge1$, $x\in D$ and $v\in\C^n$,
define, with all suprema over holomorphic functions,
\begin{align}
 E_{p,\mu}(x)&=\sup_{\|f\|_{L^p(\mu)}\le1}|f(x)|,\qquad
 D_{p,\mu}(x;v)=\sup_{\substack{\|f\|_{L^p(\mu)}\le1\\f(x)=0}}|df_x(v)|,
 \label{eq:value-kernel}\\
 B_{p,\mu}(x;v)&=D_{p,\mu}(x;v)/E_{p,\mu}(x),\qquad
 \gamma_D(x;v)=\sup_{\substack{f(x)=0\\\|f\|_\infty\le1}}|df_x(v)|.
 \notag
\end{align}
At $p=2$ the first metric is the weighted Bergman metric, in its scalar
normalization. Chen--Zhang~\cite[Proposition~2.21]{CZ} prove convergence
to the Carath\'eodory metric for Lebesgue measure; the same normal-family
argument applies to the measures used here
(Lemma~\ref{lem:basic-comparison}). Kernel regularity and
metric properties of the $p$-Bergman theory are studied
in~\cite{ChenXiong,ChenXiongRegularity,LiYinji2023}, and
Bao--Guan--Sun~\cite{BaoGuanSun2025} treat $p$-Bergman kernels with
respect to general holomorphic functionals. We ask how fast
$B_{p,\mu}$ converges to $\gamma_D$, and which geometric and measure
data determine the rate.

On the ball with normalized volume,
$(B_p(0;v)/\gamma(0;v))^p=\Gamma(n+1+p/2)/(\Gamma(n+1)\Gamma(1+p/2))$:
by subharmonicity $1$ is the value extremal at $0$, and averaging
$e^{-i\theta}f(e^{i\theta}z)$ and $f\circ U$ over $\theta$ and over
unitary $U$ fixing $v$ turns a jet competitor $f$, without increasing
its norm, into $\langle z,v\rangle/|v|^2$, whose norm is a beta
integral. Hence
$p\log(B_p/\gamma)=n\log p-\log(2^nn!)+o(1)$; the constant term, here
$-\log(2^nn!)$, is not treated in this paper. The natural scale is
therefore $(\log p)/p$, and the object of this paper is the
\emph{coefficient}
\begin{equation}\label{eq:intro-coefficient}
 N=\lim_{p\to\infty}\frac p{\log p}
 \log\frac{B_{p,\mu}(x;v)}{\gamma_D(x;v)},
\end{equation}
together with its lower and upper limits when the limit is not known
to exist. It is governed by the contact geometry of $\partial D$ along
an extremal disc.

Write $\delta_D=\operatorname{dist}(\cdot,\partial D)$. At a smooth
convex finite-type boundary point $q$, Yu's tangential linear multitype
$(\tau_1(q),\ldots,\tau_{n-1}(q))$ records the contact orders along
its complex-line flag; it agrees with the tangential Catlin
multitype~\cite{Yu1992Multitypes}. Put
\[
 T(q)=\sum_{j=1}^{n-1}2/\tau_j(q)\le n-1.
\]
This \emph{tangential dimension} is $n-1$ at strongly convex points and
decreases as complex tangential directions flatten. A complex geodesic
is a Kobayashi-isometric holomorphic disc; by Lempert's
theorem~\cite{Lempert81,Lempert1982Retracts} (see
also~\cite[Chapter~11]{JarnickiPflug2013}
and~\cite[Proposition~4.2]{ArosioFiacchi2025}) one passes through every jet
of a bounded convex domain, and it has a holomorphic left inverse. On a
smooth convex domain of finite type every complex geodesic $\varphi$
extends H\"older-continuously to $\overline\Delta$
(Lemma~\ref{lem:holder}, a case of Mercer's theorem for $m$-convex
domains~\cite[Proposition~2.9]{Mercer1993}; see
also~\cite{Bharali2016,Huang1995}), although in general not as a $C^1$
map~\cite[Example~4.8]{ArosioFiacchi2025}. We
write $\varphi^*$ for its boundary map and $\omega_a$ for harmonic
measure of the unit disc $\Delta$ at $a$.

The following heuristic explains why harmonic measure appears. Let
$\varphi$ be a complex geodesic with $\varphi(0)=x$, let $f$ be a jet
competitor, and put $g=f\circ\varphi$, so $g(0)=0$ and $g'(0)$ is
fixed. Submean inequalities on polydiscs transverse to the disc bound
$\|f\|_p^p$ below by $\int_\Delta|g|^pW\,dA$, where $W(\lambda)$ is,
up to constants, $\delta_D^\alpha$ times the volume of a transverse
polydisc at $\varphi(\lambda)$, divided by $(1-|\lambda|^2)^2$. Near a
boundary point $\varphi^*(\zeta)$ the multitype makes this volume
comparable to $s^{2+T(\varphi^*(\zeta))}$ at depth $s=1-|\lambda|^2$, so
$W(r\zeta)\approx s^{\alpha+T(\varphi^*(\zeta))}$. Because of the factor
$r^p$, the integral is decided near the circle $|\lambda|=1-1/p$. There
Jensen's inequality for $g(\lambda)/\lambda$ and the
arithmetic--geometric mean inequality give
\[
 \int_{\partial\Delta}|g(r\zeta)|^pW(r\zeta)\,\frac{|d\zeta|}{2\pi}
 \ge r^p\exp\int_{\partial\Delta}\log W(r\zeta)\,\frac{|d\zeta|}{2\pi}
 \approx p^{-\alpha-\bar T},\qquad
 \bar T=\int_{\partial\Delta}T(\varphi^*)\,\frac{|d\zeta|}{2\pi},
\]
and the shell $1/p\le s\le2/p$ contributes a further factor $1/p$.
Equality in Jensen's inequality requires $|g|^pW$ to be almost constant
on that circle, that is, $|g(r\zeta)|^p\approx p^{T(\varphi^*(\zeta))}$
up to a constant factor. This is achieved by an outer function: if $F$
is holomorphic on $\Delta$ with $\operatorname{Re}F=T\circ\varphi^*$ on
the circle, then $\lambda\exp((\log p/p)(F(\lambda)-F(0)))$ has the
required size, and in $D$ the corresponding competitor is
$u\exp((\log p/p)(F(u)-F(0)))$ for a left inverse $u$ of $\varphi$.
The correction factor is $1+O(\log p/p)$ in the interior, so it is
invisible at a fixed point, yet it changes the coefficient of $\log p$;
at a base point $\varphi(a)$, $\bar T$ becomes $\int T(\varphi^*)\,d\omega_a$.
Theorem~\ref{thm:general-main}(ii) below makes the first half of this
argument rigorous, and Theorem~\ref{thm:regular-reverse} the second half
at coercive boundary points, where the fibers of a left inverse with
affine fibers can be controlled in Stolz approach regions, defined in
Section~\ref{sec:preliminaries}.

\begin{maintheorem}[The window and the averaged upper bound]
\label{thm:general-main}
Let $D\Subset\C^n$ be a convex domain with $C^2$ boundary and let
$d\mu\asymp\delta_D^\alpha dV$ be a probability measure, $\alpha>-1$.
\begin{enumerate}[(i)]
\item For each compact $K\subset D$ there is $C_K$ such that, for
$x\in K$, $v\ne0$, and $p\ge2$,
\[
 (1+\alpha)\log p-C_K
 \le p\log\frac{B_{p,\mu}(x;v)}{\gamma_D(x;v)}
 \le(n+\alpha)\log p+C_K.
\]
\item If $\partial D$ is smooth and of finite D'Angelo type, then
for every complex geodesic $\varphi$, $a\in\Delta$, and $\eta\ne0$,
\begin{equation}\label{eq:general-averaged-upper}
 \limsup_{p\to\infty}\frac p{\log p}
 \log\frac{B_{p,\mu}(\varphi(a);\eta\varphi'(a))}
               {\gamma_D(\varphi(a);\eta\varphi'(a))}
 \le 1+\alpha+\int_{\partial\Delta}T(\varphi^*(\zeta))\,d\omega_a(\zeta).
\end{equation}
\item The bounds in (i) hold, as lower and upper limits, for moving
jets $(x_p,v_p)$ with $\log(1/\min\{1,\delta_D(x_p)\})=o(\log p)$;
Proposition~\ref{prop:boundary-window} gives an estimate for all $x$
and $p$.
\end{enumerate}
\end{maintheorem}

Part (ii) holds for every complex geodesic $\psi$ through the jet. With
\[
 N_\psi(\lambda_0)=1+\alpha+\int_{\partial\Delta}T(\psi^*(\zeta))\,d\omega_{\lambda_0}(\zeta),
\]
the upper limit at $(x;v)$ is therefore at most the infimum of
$N_\psi(\lambda_0)$ over all complex geodesics $\psi$ and
$\lambda_0\in\Delta$ with $\psi(\lambda_0)=x$ and
$\psi'(\lambda_0)\in\C^*v$.

In every dimension, the averaged upper bound is an equality along every
complex geodesic whose boundary values are coercive almost everywhere,
in the following sense. At a smooth
convex finite-type boundary point $q$, in unitary coordinates centered
at $q$, with $z_0$ the complex normal coordinate and $z'$ adapted to
Yu's flag~\cite[Section~5, (5.1)]{ArosioFiacchi2025}, a local
defining function is
$\operatorname{Re}z_0+H(z')+o(|z_0|+\sum_j|z_j|^{m_j})$, where
$(m_j)$ is a permutation of $(\tau_j(q))$ and $H\ge0$ is a convex
polynomial, weighted homogeneous for the weights $1/m_j$
(see~\eqref{eq:multitype-normal-form}). We call $q$
\emph{coercive} if $H(z')>0$ for $z'\ne0$; strongly convex points are
coercive. In $\C^2$, with a local defining function $\rho$, a point $q$
of tangential type $2m$ is coercive exactly when
\begin{equation}\label{eq:real-tangential-contact}
 \rho(q+we)\ge c_q|w|^{2m}\qquad(w\in\C,\ |w|<r_q)
\end{equation}
for a unit complex tangent vector $e$ at $q$ and some $c_q,r_q>0$, and
then $T(q)=1/m$.

\begin{maintheorem}[Equality at coercive boundary points]
\label{thm:regular-reverse}
Let $D\Subset\C^n$, $n\ge2$, be a smooth convex domain of finite type,
let $d\mu\asymp\delta_D^\alpha dV$, $\alpha>-1$, be a probability
measure, and let $\varphi:\Delta\to D$ be a complex geodesic.
\begin{enumerate}[(i)]
\item Let $\mathsf C\subset\partial\Delta$ be a Borel set such that
$\varphi^*(\zeta)$ is coercive for every $\zeta\in\mathsf C$. Then for
every $a\in\Delta$ and $\eta\ne0$,
\[
 \liminf_{p\to\infty}\frac p{\log p}
 \log\frac{B_{p,\mu}(\varphi(a);\eta\varphi'(a))}
               {\gamma_D(\varphi(a);\eta\varphi'(a))}
 \ge1+\alpha+\int_{\mathsf C}T(\varphi^*(\zeta))\,d\omega_a(\zeta).
\]
\item If $\varphi^*(\zeta)$ is coercive for almost every
$\zeta\in\partial\Delta$, then for every $a\in\Delta$ and $\eta\ne0$
the coefficient exists and
\begin{equation}\label{eq:regular-reverse-law}
 \lim_{p\to\infty}\frac p{\log p}
 \log\frac{B_{p,\mu}(\varphi(a);\eta\varphi'(a))}
               {\gamma_D(\varphi(a);\eta\varphi'(a))}
 =1+\alpha+\int_{\partial\Delta}T(\varphi^*(\zeta))\,d\omega_a(\zeta).
\end{equation}
\item If $\mathsf C$ in (i) has positive measure, then $\varphi$ is the only
complex geodesic $g$ of $D$ with $g(0)=\varphi(0)$ and
$g'(0)=\varphi'(0)$.
\end{enumerate}
\end{maintheorem}

On smooth strongly convex domains $T\equiv n-1$, and (ii) gives the
upper endpoint $n+\alpha$ of the window; every value in
$(1+\alpha,n+\alpha]$ occurs (Example~\ref{ex:continuum}). We do not
know whether the lower endpoint $1+\alpha$ is attained when $n\ge2$.
The coercivity constants ($c_H$ in~\eqref{eq:H-coercive}, or $c_q,r_q$
in~\eqref{eq:real-tangential-contact}) may vary along the circle, and
the multitype need not be locally constant. The proof uses a left
inverse of $\varphi$ with affine fibers. On strongly convex domains
Lempert's left inverses have this property~\cite{Lempert81}, and on
every bounded convex domain one is given by the dual map of Royden and
Wong~\cite[Theorem~2.1(iii)]{KosinskiMcCarthy2019}; we use it in the
refined form of Zwonek~\cite[Theorem~2.1]{Zwonek2025LeftInverses}
(Lemma~\ref{lem:affine-inverse}).

By Lempert's theory every holomorphic disc $g:\Delta\to D$ with
$g(0)=\varphi(0)$ and $g'(0)=\varphi'(0)$ is a complex geodesic (for a
holomorphic left inverse $G$ of $\varphi$, Schwarz's lemma gives
$G\circ g=\Id$), so (iii) is the same as uniqueness among holomorphic
discs with this $1$-jet. Uniqueness of the complex geodesic through a
jet is classical on strongly convex
domains~\cite{Lempert81},~\cite[Corollary~2.6.30]{Abate1989}, on
strictly convex domains~\cite[Section~11.3]{JarnickiPflug2013}, and on
strongly linearly convex
domains~\cite{Lempert1984Varna,KosinskiWarszawski2013}. On strictly
convex domains it follows by averaging: $\frac12(\varphi+g)$ is a
complex geodesic, so $[\varphi^*(\zeta),g^*(\zeta)]\subset\partial D$
for almost every $\zeta$ (see the proof
of~\cite[Corollary~2.6.30]{Abate1989}). The $H^1$ dual map of Royden and
Wong, which on strongly convex domains is unique up to a positive
factor~\cite[Lemma~2.6.31]{Abate1989}, places
$g^*(\zeta)-\varphi^*(\zeta)$ in the complex tangent space, and this
gives Corollary~\ref{cor:uniqueness} below, a uniqueness criterion of a
known type. We call $\overline D\cap(q+T^{1,0}_q\partial D)$ the
\emph{complex face} of $q\in\partial D$. Recall that $q\in\partial D$ is a \emph{complex extreme
point} of $\overline D$ if no $v\ne0$ satisfies $q+\zeta v\in\overline D$
for all $|\zeta|\le1$; for the unit ball of a norm, the geodesic through
$0$ and $z$ is unique up to reparametrization if and only if $z/\|z\|$
is a complex extreme point~\cite[Proposition~2.6.8]{Abate1989}. At a
$C^1$ boundary point the trivial complex face in
Corollary~\ref{cor:uniqueness} is the stronger condition: if
$q+\zeta v\in\overline D$ for $|\zeta|\le1$ and $r$ is a $C^1$ local
defining function, then $\zeta\mapsto r(q+\zeta v)\le0$ attains its
maximum $0$ at $\zeta=0$, so its differential there vanishes,
$v\in T^{1,0}_q\partial D$, and $q+v\ne q$ lies in the complex face. The
converse fails: with $\chi$ as in Section~\ref{sec:flat-convex}, the
domain
$\{(z,w):|z|^2+(\operatorname{Im}w)^2+\chi(|\operatorname{Re}w|-1/2)<1\}$
has D'Angelo type $2$, and its non-coercive boundary points
$(\zeta,x)$, $|\zeta|=1$, $|x|\le1/2$, are complex extreme points whose
complex face is the segment $\{\zeta\}\times[-1/2,1/2]$, while the
points of its boundary segments $\{(0,x\pm i):|x|\le1/2\}$ are
coercive and have trivial complex face (Section~\ref{sec:regular-reverse}).
On the domains $\Omega_\ell$ of
Theorem~\ref{thm:noncoercive-main} every point of the boundary circle of
the axis disc has trivial complex face but is not coercive.

\begin{corollary}\label{cor:uniqueness}
Let $D\Subset\C^n$ be a convex domain and $\varphi:\Delta\to D$ a
complex geodesic. Suppose that there is a measurable set
$\mathsf E\subset\partial\Delta$ of positive measure such that for every
$\zeta\in\mathsf E$ the radial limit $q=\varphi^*(\zeta)$ exists,
$\partial D$ is a $C^1$ hypersurface near $q$, and $q$ has
\emph{trivial complex face}:
\[
 \overline D\cap\bigl(q+T^{1,0}_q\partial D\bigr)=\{q\},\qquad
 T^{1,0}_q\partial D=T_q\partial D\cap iT_q\partial D .
\]
Then every holomorphic $g:\Delta\to D$ with $g(0)=\varphi(0)$ and
$g'(0)=\varphi'(0)$ equals $\varphi$, and, up to reparametrization by a
disc automorphism, $\varphi$ is the unique complex geodesic through each
of its jets. This applies if $D$ is smooth of finite type and
$\varphi^*(\zeta)$ is coercive for all $\zeta$ in a set of positive
measure, and to the axis disc of $\Omega_\ell$.
\end{corollary}
The proof is given at the end of Section~\ref{sec:regular-reverse}. Uniqueness of the geodesic does not entail uniqueness of its left
inverse: on smooth strongly linearly convex domains in $\C^n$, $n>1$,
left inverses are never unique~\cite[Theorem~3.1]{KosinskiZwonek2016}.

For example, in the smooth convex domain $D_\eps\Subset\C^2$ of
Proposition~\ref{prop:three-types} the boundary circle of the axis disc
lies in $\partial D_\eps$ and splits into sets $E_2,E_4,E_6$ of
length $2\pi/3$ each ($E_4$ is two arcs), on which the tangential type
is $2$, $4$, $6$; every point of the circle, including the four
transition points, is coercive. Theorem~\ref{thm:regular-reverse}(ii) gives the coefficient
$N(0)=1+\alpha+\frac13(1+\frac12+\frac13)=\alpha+\frac{29}{18}$ at the
center, a value attained by none of the three contact types, and a
nonconstant harmonic coefficient along the disc; by
Corollary~\ref{cor:uniqueness} the axis disc is, up to
reparametrization, the unique complex geodesic through each of its
jets. The remainder after $N\log p$ tends to $-\infty$
(Corollary~\ref{cor:three-types-remainder}).

At points of type $2$, coercivity is Lempert's strong linear convexity
condition (Section~\ref{sec:regular-reverse}), and it can fail: a convex
boundary given near the origin by
$\operatorname{Re}z_0+(\operatorname{Im}z_1)^2+|z_1|^4=0$ is strongly
pseudoconvex there, but $H(z_1)=(\operatorname{Im}z_1)^2$ vanishes on
$\R$.
Theorem~\ref{thm:noncoercive-main} below shows that then the averaged
bound can be strict. The coefficient is invariant under biholomorphisms
which, with their inverses, are $C^1$ up to the boundary
(Remark~\ref{rem:twisted}), whereas coercivity is a complex-affine but
not a biholomorphic invariant (Section~\ref{sec:regular-reverse}).

\begin{figure}[tbp]
\centering
\includegraphics[width=\textwidth]{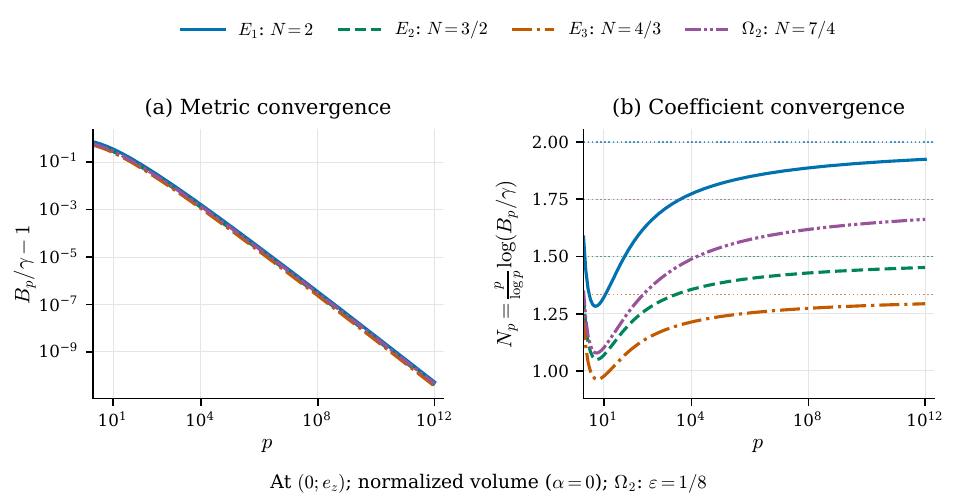}
\caption{Finite-$p$ convergence at the central axis jet with normalized
volume. The curves evaluate
\eqref{eq:illustration-ellipsoid}--\eqref{eq:illustration-twisted},
using $\eps=1/8$ for $\Omega_2$; the dotted lines in (b) mark the
proved limits. All four metric ratios tend to $1$, while their
logarithmic coefficients differ. In particular $E_1$ and $\Omega_2$
both have type $2$ along the extremal circle but have different
coefficients. Here $N_p=N+O(1/\log p)$, and even large finite-$p$ values can differ
substantially from $N$. The quantities $N_p$ need not be monotone or
stay in the limiting window. These curves concern only the specified
central jets and measures. The source package contains
the quadrature, data, and refinement checks; these numerical checks do
not replace the analytic proof.}
\label{fig:finite-p-convergence}
\end{figure}

Without coercivity, pointwise data on the circle that agree up to a
unitary equivalence at each point (multitype, Levi form, model
polynomial) no longer determine the coefficient. Figure~\ref{fig:null-winding}
displays the global distinction, and Figure~\ref{fig:finite-p-convergence}
compares exact finite-$p$ formulas with their limiting coefficients.

\begin{maintheorem}[Non-coercive contact]\label{thm:noncoercive-main}
For an integer $\ell\ge1$ and $0<\eps\le1/(2\ell^2)$ let
\[
 \Omega_\ell=\bigl\{(z,w)\in\Delta\times\C:\ |w|<\sqrt2,\
 |z|^2+\eps\bigl(|w|^2-\operatorname{Re}(\bar z^\ell w^2)\bigr)+|w|^4<1\bigr\}.
\]
Then $\Omega_\ell$ is a bounded convex domain with real-algebraic
strongly pseudoconvex boundary. Every point of the boundary circle $\mathfrak c$ of
the axis disc $\varphi(\lambda)=(\lambda,0)$ is strongly pseudoconvex
and non-coercive, of tangential type $2$; up to reparametrization,
$\varphi$ is the unique complex geodesic through each of its jets, and
the right side of~\eqref{eq:general-averaged-upper} equals $2+\alpha$.
Let $d\mu\asymp\delta_{\Omega_\ell}^\alpha dV$ be a probability measure,
$\alpha>-1$, let $v\ne0$, $\gamma=\gamma_{\Omega_\ell}$, and let
$K\subset\Delta$ be compact.
\begin{enumerate}[(i)]
\item If $\ell=1$, then
$p\log(B_{p,\mu}/\gamma)((x,0);(v,0))=(2+\alpha)\log p+O_K(1)$ for $x\in K$.
\item If $\ell\ge2$, then
$p\log(B_{p,\mu}/\gamma)((x,0);(v,0))=(\frac74+\alpha)\log p+O_K(1)$
for every $x\in K$, with no parity restriction on $\ell$.
\end{enumerate}
\end{maintheorem}

The axis geodesic is unique through each of these jets, so taking an
infimum over extremal discs does not repair the strict inequality in
(ii). At each point of $\mathfrak c$, the domains $\Omega_1$ and
$\Omega_2$ have unitarily equivalent tangential type, Levi form and
model polynomial, yet different coefficients
(Corollary~\ref{cor:noncoercive-answer}). The distinction is global.
For $\Omega_1$, the curved left inverse $z-\eps w^2/2$ removes the real
null direction. For $\Omega_\ell$, $\ell\ge2$, the null line turns in a
way that prevents one holomorphic quadratic correction from doing so
around the circle. This obstruction yields a quantitative estimate
without any rotational symmetry assumption.

\begin{maintheorem}[Winding gap]\label{thm:winding-gap}
Let $D\Subset\C^2$ be convex with $C^2$ boundary and suppose
\[
 \Delta\times\{0\}\subset D\subset\Delta\times\C,\qquad
 \overline D\cap(\partial\Delta\times\C)=\mathfrak c.
\]
Assume that the boundary is $C^4$ near $\mathfrak c$. Write its radial
graph there as
\[
 |z|^2=1-Q_{\arg z}(w),\qquad
 Q_\theta(w)=\mathsf A(e^{i\theta})|w|^2+
       \operatorname{Re}(\mathsf B(e^{i\theta})w^2)+o(|w|^2).
\]
Suppose $\mathsf A,\mathsf B$ are smooth on the circle,
$\mathsf A>0$ and $|\mathsf B|=\mathsf A$. Set
$d=\operatorname{wind}(\mathsf B|_{\partial\Delta})$.
For a probability measure $d\mu\asymp\delta_D^\alpha dV$, $\alpha>-1$,
and every $K\Subset\Delta$, uniformly for $a\in K$,
\begin{enumerate}[(i)]
\item if $d\ge-1$, then
$p\log(B_{p,\mu}/\gamma_D)((a,0);e_z)=(2+\alpha)\log p+O_K(1)$;
\item if $d\le-2$, then
$p\log(B_{p,\mu}/\gamma_D)((a,0);e_z)\le(\frac74+\alpha)\log p+C_K$.
\end{enumerate}
\end{maintheorem}

Here the quadratic form is positive in one real tangential direction
and zero in the other. If contact in the null direction is quartic,
(ii) is an equality with bounded remainder
(Corollary~\ref{cor:quartic-equality}). This gives Theorem~\ref{thm:noncoercive-main}(ii)
and persists under explicit real-algebraic perturbations breaking the
weighted circle action (Proposition~\ref{prop:perturbed-twisting}).
Positive H\"older control of the second derivatives still gives a
strict gap (Corollary~\ref{cor:holder-gap}), but a $C^2$ example with
$d=-2$ has coefficient $2+\alpha$
(Proposition~\ref{prop:C2-obstruction}). The exact value in the negative-winding case is established here
under the quartic contact condition.

The disc problem behind these results is an extremal problem of
Christoffel type: minimize $\int|f|^p\,d\nu$ over holomorphic $f$ with
$f(x)=0$ and $f'(x)=1$. For polynomials of growing degree these are
Christoffel functions, with Szeg\H{o}--Widom
asymptotics~\cite{Widom1969,Totik2000,LevinLubinsky2015}; here the
degree is replaced by $p$, and the weight is sampled at depth $1/p$
from the circle. Unlike the boundary behavior of the invariant
metrics~\cite{Graham1975,Ma1991,FornaessLee2009}, the coefficient is
attached to a fixed interior jet. Harmonic measure along complex
geodesics also appears in the pluricomplex Poisson kernel of a convex
domain of finite type, which on a complex geodesic ending at its pole
is a negative multiple of the Poisson kernel of $\Delta$ with pole
$1$~\cite[Corollary~4.4]{ArosioBracciFiacchi2026}; here, dually, the
Poisson kernel at a fixed interior point weights all boundary points of
the extremal disc.

The proofs compare two sides of the minimum $L^p$ norm for a fixed
jet. Restricting every competitor to an extremal disc and its
multitype polydiscs gives the lower norm bound in
Theorem~\ref{thm:general-main}. For Theorem~\ref{thm:regular-reverse},
coercivity bounds the full scaled affine fibers of a left inverse;
the resulting pushforward density allows a $p$-dependent outer
function to give the opposite norm bound. Theorem~\ref{thm:winding-gap}
uses a different restriction family: on a double cover of the axis
disc, holomorphic curves follow the real null direction around the
entire boundary. The cover accommodates odd winding, and symmetric
averaging preserves the prescribed jet. This yields a bound for every
holomorphic competitor, rather than only for functions of a fixed
left inverse.

Sections~\ref{sec:geodesic-restriction} and~\ref{sec:regular-reverse}
prove Theorems~\ref{thm:general-main} and~\ref{thm:regular-reverse};
Section~\ref{sec:noncoercive} proves
Theorems~\ref{thm:noncoercive-main} and~\ref{thm:winding-gap}.
The scaling input is due to Gaussier and Arosio--Fiacchi
\cite{Gaussier1997,ArosioFiacchi2025}; the polydisc and weighted
geometric-mean inputs are from
\cite{McNeal1994Kernels,Yu1992Multitypes,Garnett,Nakazi2011}.
Smooth means $C^\infty$, $dA,dV$ denote Lebesgue measure, and comparison
constants are independent of $p$.
\section{Preliminaries}
\label{sec:preliminaries}

\paragraph{Notation}
Throughout, $\Delta$ is the unit disc, $s=s(z)=1-|z|^2$ is the depth of
$z\in\Delta$, and $\omega_x$ is harmonic measure of $\Delta$ at $x$:
$d\omega_x(e^{i\theta})=\mathcal P_x(\theta)\,d\theta/2\pi$ with the
Poisson kernel $\mathcal P_x(\theta)=(1-|x|^2)/|e^{i\theta}-x|^2$. We
write $\mathcal P[\phi]$ for the Poisson extension of a function $\phi$
on $\partial\Delta$ and $\mob_a(\zeta)=(\zeta-a)/(1-\bar a\zeta)$ for
the disc automorphisms, and $\mathbb B_b$ is the open Euclidean ball of
radius $b$ about $0$. For a domain $D$, $\delta_D$ is the Euclidean
distance to $\partial D$ (in a manifold, the distance in a fixed smooth
Hermitian metric). In fibered domains $(z,w)\in\Delta\times\C^k$ the
axis disc is $\Delta\times\{0\}$, its boundary circle is
$\mathfrak c=\partial\Delta\times\{0\}$, and $e_z=(1,0)$. The letter $N$
always denotes a coefficient of $\log p$ or, as in $N_\psi$, its
averaged upper bound.

\paragraph{Stolz approach regions}
For $\zeta\in\partial\Delta$ and an aperture parameter $M>1$, define
\begin{equation}\label{eq:stolz-region}
 \mathcal A_M(\zeta)
 =\{\lambda\in\Delta:|\lambda-\zeta|<M(1-|\lambda|)\}.
\end{equation}
We call this a \emph{Stolz approach region}, also called a Stolz angle.
Here $M$ is a ratio bound, not an angle in radians. A limit at $\zeta$
is \emph{nontangential} if it holds as $\lambda\to\zeta$ within every
fixed $\mathcal A_M(\zeta)$. An approach sequence is nontangential exactly when
$|\lambda-\zeta|/(1-|\lambda|)$ remains bounded.
The truncated regions used below are
$\mathcal A_M(\zeta)\cap\{s(\lambda)<s_0\}$, where
$1-|\lambda|\le s(\lambda)\le2(1-|\lambda|)$. Thus replacing the
boundary distance by $s$ gives the same family of approaches after
changing $M$. It suffices to use integer $M\ge2$, since these regions
contain every region with a fixed aperture. Constants may depend on
$M$ and $\zeta$; no uniformity as $M\to\infty$ is assumed.
These regions lie in the parameter disc, not in the ambient domain $D$.

\paragraph{Results of~\cite{ArosioFiacchi2025}}
Statement numbers in~\cite{ArosioFiacchi2025},
\cite{ArosioBracciFiacchi2026} and~\cite{Zimmer2016} refer to
arXiv:2407.09199v2, arXiv:2509.26230v3 and arXiv:1405.2858v3. We use:
Definition~2.1 ($\C$-proper convex domains, containing no complex
affine line); Definition~2.3 (boundary points of finite line type);
Proposition~4.2 (on such domains extremal discs are complex geodesics,
any two points lie on one, and each has a left inverse);
Definitions~4.4--4.5 (reference rays in $\Delta$ and in the left
half-plane, and regular contact points); Remark~4.7 (a complex geodesic
with an endpoint of locally finite type has a regular contact point
there); Example~4.8 (complex geodesics of egg domains which are not
$C^1$ up to the boundary); Proposition~4.9 (at a regular contact point
the normal component of the derivative has a positive nontangential
limit); Section~5 and~(5.1) (Yu's flag of complex subspaces and the
multitype normal form with a nonnegative, weighted homogeneous,
non-degenerate convex polynomial $H$); Theorem~6.2 (the scaled domains
converge to the model $\{\operatorname{Re}z_0+H(z')<0\}$);
Proposition~6.6 (scaled holomorphic maps from the left half-plane with a
regular contact point converge to $\sigma\mapsto(\kappa\sigma,0)$).

\subsection{Extremal norms}
At a point and direction with $\gamma_D(x;v)>0$ put
\[
 I_{p,\mu}(x;v)=\inf\{\|f\|_{L^p(\mu)}:
 f\in\mathcal O(D),\ f(x)=0,\ df_x(v)=1\}.
\]
Normalizing a competitor in~\eqref{eq:value-kernel} gives
$D_{p,\mu}=I_{p,\mu}^{-1}$.

\begin{lemma}[Value norms and comparison]\label{lem:basic-comparison}
Let $\mu$ be a probability measure on a connected complex manifold $D$,
locally comparable to smooth positive volume. For each compact
$K\subset D$ there is $C_K$ with $1\le E_{p,\mu}(x)^p\le C_K$ for
$x\in K$ and $p\ge1$. Consequently
\begin{equation}\label{eq:B-from-I}
 B_{p,\mu}=\frac1{I_{p,\mu}E_{p,\mu}},\qquad
 \log\frac{B_{p,\mu}}{\gamma_D}
        =-\log(\gamma_D I_{p,\mu})+O_K(p^{-1}).
\end{equation}
One always has $B_{p,\mu}\ge\gamma_D$, and
$B_{p,\mu}(x;v)\to\gamma_D(x;v)$ at every fixed jet as $p\to\infty$.
If $C^{-1}\mu_0\le\mu\le C\mu_0$
for probability measures, then
$C^{-2/p}B_{p,\mu_0}\le B_{p,\mu}\le C^{2/p}B_{p,\mu_0}$.
\end{lemma}
\begin{proof}
The constant function gives $E_p\ge1$. On a fixed coordinate ball
about $x$, the submean property of the plurisubharmonic function
$|f|^p$ gives $|f(x)|^p\le C\|f\|_p^p$ with $C$ independent of $p$; a
finite cover makes this uniform on $K$. If $u(x)=0$ and
$\|u\|_\infty\le1$, then $uf$ is an admissible competitor whenever $f$
is, and $d(uf)_x(v)=f(x)du_x(v)$; taking suprema gives
$D_p\ge E_p\gamma_D$. For convergence, take $p_j\to\infty$ realizing the upper limit and
$f_j(x)=0$, $\|f_j\|_{p_j}\le1$, rotated so that their derivatives approximate
$D_{p_j,\mu}(x;v)$ with error tending to zero. The local submean bounds
make $(f_j)$ normal, and every locally uniform limit $f$ satisfies
$|f(y)|\le\lim C_y^{1/p_j}=1$ for every $y$. Derivatives converge, so
$\limsup D_{p_j,\mu}(x;v)\le\gamma_D(x;v)$; use $E_{p,\mu}(x)\to1$
and the opposite inequality. This is the argument
of~\cite[Proposition~2.21]{CZ} with the local volume comparison included.
Finally $D_p$ and $E_p$ change by at most the factor $C^{1/p}$ under
the comparison of measures.
\end{proof}

In particular $p\log(B_{p,\mu}/\gamma_D)=N\log p+O(1)$ is equivalent
to $(\gamma_DI_{p,\mu})^p\asymp p^{-N}$. For a $C^2$ defining function
$r<0$ of a relatively compact domain with $dr\ne0$ on the boundary,
$-r\asymp\delta$, by a normal collar expansion.

\subsection{Line types and multitype}\label{sec:multitype}
Let $D$ have smooth boundary near $q\in\partial D$ and be convex near
$q$, with local defining function $r$. For a unit vector $v$ in the
complex tangent space $T^{1,0}_q$, let $t(v)\in\{2,3,\ldots\}\cup\{\infty\}$
be the vanishing order at $0$ of $\zeta\mapsto r(q+\zeta v)$, and put
$t(0)=\infty$. The \emph{linear multitype} is obtained by successively maximizing
$t$ over unit tangent vectors orthogonal to those already chosen.
For convex finite-type boundaries, Yu~\cite{Yu1992Multitypes} identifies
these entries, in increasing order and preceded by $1$, with Catlin's
multitype; the sets $\{t\ge d\}$ are complex subspaces
\cite[Section~5]{ArosioFiacchi2025}. D'Angelo type equals maximal
line type~\cite{McNeal1992Convex}; see~\cite{DAngelo1982,Catlin1984}
for the general definitions.

\begin{lemma}\label{lem:linear-multitype}
Suppose $\max_vt(v)<\infty$ and every set $S_d=\{v\in T^{1,0}_q:t(v)\ge d\}$,
$d\ge2$, is a complex linear subspace. Then the linear multitype
$(\tau_1,\ldots,\tau_{n-1})$ satisfies
$\#\{j:\tau_j\ge d\}=\dim_{\C}S_d$ for every $d$, and it is the
tangential part of the Catlin multitype.
\end{lemma}
\begin{proof}
The maxima over the decreasing sets $\{v_1,\ldots,v_{j-1}\}^\perp$ do
not increase, so the vectors with $t\ge d$ are chosen first, and they
lie in $S_d$. They are orthonormal, so there are at most $\dim S_d$
of them. If fewer have been chosen, the orthogonal complement of the
chosen vectors meets $S_d$ in a nonzero vector, so the next choice
again has $t\ge d$. This proves the count; the last assertion is Yu's
theorem.
\end{proof}

\subsection{The collar estimate for a fixed extremal}
\label{sec:collar}
For a bounded domain $D\subset\C^n$, $x\in D$ and $v\ne0$, let
$\mathcal E(x,v)$ be the set of holomorphic $u:D\to\Delta$ with
$u(x)=0$, $|du_x(v)|=\gamma_D(x;v)$, continuous on $\overline D$ and
$C^2$ up to the boundary near $\mathcal K_u=\{q\in\partial D:|u(q)|=1\}$.
Put $G_u=-\log|u|^2\ge0$ on $\partial D$, with $G_u=+\infty$ at zeros of
$u$, and let $\beta(u)$ be the supremum of the exponents $\beta\ge0$
with $\operatorname{vol}_{\partial D}\{q\in\partial D:G_u(q)<t\}=O(t^\beta)$ as $t\to0$, where
$\operatorname{vol}_{\partial D}$ is surface measure.

\begin{lemma}[Collar estimate]\label{lem:collar-competitor}
Let $D\subset\C^n$ be a bounded domain with $C^2$ boundary, let
$d\mu\asymp\delta_D^\alpha dV$ be a probability measure, $\alpha>-1$,
and let $u\in\mathcal E(x,v)$. For every $\beta\ge0$ with
$\operatorname{vol}_{\partial D}\{G_u<t\}=O(t^\beta)$,
\[
 \int_D|u|^p\,d\mu\le C_\beta\,p^{-(1+\alpha+\beta)},\qquad
 \log\frac{B_{p,\mu}(x;v)}{\gamma_D(x;v)}
 \ge\frac{(1+\alpha+\beta)\log p-C_\beta}p .
\]
In particular
$\liminf_{p\to\infty}(p/\log p)\log(B_{p,\mu}/\gamma_D)\ge1+\alpha+\beta(u)$.
For every compact $K\subset\Delta$ the moment bound is uniform for the
functions $\mob_a\circ u$, $a\in K$.
\end{lemma}
\begin{proof}
Extend $G_u=-\log|u|^2$ to $\overline D\cap\{u\ne0\}$. Since $u$ is
continuous with $|u|=1$ on $\mathcal K_u$, the set
$\mathcal V=\{z\in\overline D:|u(z)|>1/2\}$ is a relatively open
neighborhood of $\mathcal K_u$, and, shrinking it, $G_u$ is $C^2$ on
$\mathcal V$ up to the boundary, pluriharmonic, hence harmonic, in
$\mathcal V\cap D$, positive there because $|u|<1$ on $D$ by the maximum
principle, and zero on $\mathcal K_u$. Let $\nu$ be the outer unit
normal. At $q\in\mathcal K_u$ the interior ball condition holds since
$\partial D$ is $C^2$, so Hopf's lemma gives $-\partial_\nu G_u(q)>0$. As
$\nabla G_u$ is continuous up to the boundary and $\mathcal K_u$ is
compact, there are $c,t_0>0$ and a compact neighborhood $\mathcal N$ of
$\mathcal K_u$ in $\partial D$ such that $(q,t)\mapsto q-t\nu(q)$ is a
$C^1$ diffeomorphism of $\mathcal N\times[0,t_0]$ onto a collar
$\mathcal U\subset\mathcal V$ and $\frac d{dt}G_u(q-t\nu(q))\ge2c$
there. Integrating in $t$ gives $G_u(q-t\nu(q))\ge G_u(q)+2ct$; as
$\delta_D\asymp t$ on $\mathcal U$, in particular $G_u\ge c'\delta_D$
there. The compact set $\overline D\setminus\operatorname{int}\mathcal U$
does not meet $\mathcal K_u$, so $|u|\le e^{-c_0}<1$ on it. With
$|u|^p=e^{-pG_u/2}$, normal coordinates on $\mathcal U$ and the sublevel
hypothesis yield
\[
 \int_D|u|^p\,d\mu\le
 Cp^{-1-\alpha}\int_0^\infty\tfrac p2e^{-p\tau/2}
            \operatorname{vol}_{\partial D}\{G_u<\tau\}\,d\tau+e^{-pc_0}
 \le C_\beta p^{-1-\alpha-\beta}.
\]
The jet competitor $u/du_x(v)$ and~\eqref{eq:B-from-I} give the result.
For $a\in K$ put $r_K=\max_{a\in K}|a|$ and
$c_K=(1-r_K)/(1+r_K)>0$. The identity
$1-|\mob_a(u)|^2=(1-|a|^2)(1-|u|^2)/|1-\bar a u|^2$
gives $1-|\mob_a(u)|^2\ge c_K(1-|u|^2)$ throughout $D$.
On the collar shrink $\mathcal U$ about $\mathcal K_u$ so that
$0\le G_u\le1$. There $1-|u|^2=1-e^{-G_u}\ge(1-e^{-1})G_u$,
hence $|\mob_a(u)|^p\le e^{-c_K(1-e^{-1})pG_u/2}$.
On the complement $1-|u|^2$ has a positive lower bound.
The collar integral above, with $p$ multiplied by the fixed factor
$c_K(1-e^{-1})$, and the exponential bound on the complement
give the uniform estimate.
\end{proof}

\subsection{Flat convex functions}\label{sec:flat-convex}
For $\nu>0$ define the smooth convex function
\begin{equation}\label{eq:chi-nu}
 \chi_\nu(t)=0\quad(t\le0),\qquad
 \chi_\nu(t)=\int_0^t(t-u)e^{-u^{-\nu}}\,du\quad(t>0),
\end{equation}
and write $\chi=\chi_{1/2}$. All derivatives vanish at $0$, and
$0\le\chi_\nu''\le1$ with $\chi_\nu'(0)=0$ gives
$(\chi_\nu')^2\le2\chi_\nu$, since $(2\chi_\nu-(\chi_\nu')^2)'
=2\chi_\nu'(1-\chi_\nu'')\ge0$. In particular $\sqrt{\chi_\nu}$ is
Lipschitz. For small $t>0$,
\begin{equation}\label{eq:chi-bounds}
 \log\chi_\nu(t)=-t^{-\nu}+O(1+|\log t|).
\end{equation}
Indeed, convexity of $u\mapsto u^{-\nu}$ gives
$u^{-\nu}\ge t^{-\nu}+\nu t^{-\nu-1}(t-u)$, hence
$\chi_\nu(t)\le\nu^{-2}t^{2\nu+2}e^{-t^{-\nu}}$; on $0\le t-u\le
t^{\nu+1}\le t/2$ the mean value theorem gives
$u^{-\nu}\le t^{-\nu}+\nu2^{\nu+1}$, hence
$\chi_\nu(t)\ge c_\nu t^{2\nu+2}e^{-t^{-\nu}}$.

\section{The window and the averaged upper bound}
\label{sec:geodesic-restriction}
We use the Kobayashi distance $k_D$ with
$k_\Delta(0,r)=\operatorname{arctanh}r$, and write
$\mathcal B_D(y,b)=\{z:k_D(y,z)<b\}$ for the Kobayashi balls. Only
bounded convex domains occur in this section. By Montel's theorem a
Kobayashi extremal disc $\varphi$ exists for every $(x,v)$, and so does
a Carath\'eodory extremal $G:D\to\Delta$ with $G(x)=0$ and
$|dG_x(v)|=\gamma_D(x;v)$. By Lempert's
theorem~\cite{Lempert81,Lempert1982Retracts} (see
also~\cite[Chapter~11]{JarnickiPflug2013}
and~\cite[Proposition~4.2]{ArosioFiacchi2025}) the Carath\'eodory and
Kobayashi metrics agree, so $G\circ\varphi$ is a self-map of $\Delta$
fixing $0$ with derivative of modulus one. Hence $\varphi$ is a complex
geodesic and, after a rotation, $G$ is a holomorphic left inverse of
$\varphi$.

\begin{lemma}[Depth along a geodesic]\label{lem:geodesic-depth}
Let $D\Subset\C^n$ be convex with $C^2$ boundary. Then
\[
 |\log\delta_D(y)-\log\delta_D(z)|\le2k_D(y,z).
\]
For each compact $K\subset D$ there are $c_K,C_K>0$ such that every
complex geodesic $\varphi$ with $\varphi(0)\in K$ satisfies
\begin{equation}\label{eq:geodesic-depth}
 c_K(1-|\zeta|^2)\le\delta_D(\varphi(\zeta))
                  \le C_K(1-|\zeta|^2).
\end{equation}
\end{lemma}
\begin{proof}
Suppose $\delta_D(z)\ge\delta_D(y)$ and take a closest boundary point
$q$ to $y$. The complex-linear projection to the inward normal half-plane
sends $y$ to $\delta_D(y)$ and sends $z$ to a point whose real part is at
least $\delta_D(z)$: the ball of radius $\delta_D(z)$ about $z$ lies on
the inner side of the supporting plane. The half-plane distance is at
least half the absolute logarithm of the ratio of the real parts.
Holomorphic contraction proves the first assertion; exchange $y,z$ for
the other order. This uses only a supporting hyperplane, so the first
inequality holds on every convex domain; the $C^2$ hypothesis enters
only in the upper bound of~\eqref{eq:geodesic-depth}. Applied to
$y=\varphi(0)$, it proves the lower bound in~\eqref{eq:geodesic-depth}.

There is a fixed radius $a>0$ such that the ball centered at
$c_q=q-a\nu(q)$ of radius $a$ lies in $D$ for every $q\in\partial D$.
The centers form a compact subset of $D$. Thus, uniformly for $x\in K$
and $y=q-d\nu(q)$ with $d=\delta_D(y)$ small,
\[
 k_D(x,y)\le C_K+k_{c_q+\mathbb B_a}(c_q,y)
       \le C_K+\tfrac12\log(2a/d).
\]
Use $k_D(\varphi(0),\varphi(\zeta))=
\frac12\log((1+|\zeta|)/(1-|\zeta|))$ and enlarge the constant on the
remaining compact set. This gives the upper bound.
\end{proof}

For $z\in D$ and $v\ne0$ let $\delta_D(z;v)$ be the distance from $z$
to $\partial D$ along the complex line $z+\C v$.

\begin{lemma}[Boundary continuity of geodesics]\label{lem:holder}
Let $D\Subset\C^n$ be a convex domain with smooth boundary of finite
D'Angelo type $L$. Every complex geodesic $\varphi:\Delta\to D$ extends
to a map $\varphi^*:\overline\Delta\to\overline D$ which is H\"older
continuous of exponent $1/L$.
\end{lemma}
\begin{proof}
H\"older continuity is a case of Mercer's theorem~\cite[Proposition~2.9]{Mercer1993} for
$m$-convex domains; we include a short argument giving the exponent $1/L$. Let $z\in D$,
$v\ne0$, and let $y=z+\delta_D(z;v)e^{i\theta}v/|v|\in\partial D$. A
real supporting hyperplane of $D$ at $y$ has the form
$\{\operatorname{Re}\Lambda(\cdot-y)=0\}$ with $\Lambda$ complex linear and
$\operatorname{Re}\Lambda(\cdot-y)<0$ on $D$; then $\Lambda(v)\ne0$,
since otherwise $\operatorname{Re}\Lambda(\cdot-y)$ would be constant on
$(z+\C v)\cap D$ and tend to $0$ at $y$. The map
$\Lambda(\cdot-y)$ sends $D$ into a half-plane whose boundary is at
distance at most $|\Lambda(v)|\delta_D(z;v)/|v|$ from $\Lambda(z-y)$, so
contraction of the Kobayashi metric gives
$k_D(z;v)\ge|v|/(2\delta_D(z;v))$. For a geodesic,
$k_D(\varphi(\zeta);\varphi'(\zeta))=(1-|\zeta|^2)^{-1}$, hence
$|\varphi'(\zeta)|\le2\delta_D(\varphi(\zeta);\varphi'(\zeta))/(1-|\zeta|^2)$.
By~\cite[Proposition~9.1]{Zimmer2016} (a local estimate near each
boundary point of finite line type), McNeal's equality of line type
and D'Angelo type~\cite{McNeal1992Convex}, and compactness of
$\partial D$, $\delta_D(z;v)\le C\delta_D(z)^{1/L}$ on $D$, and by
Lemma~\ref{lem:geodesic-depth} $\delta_D(\varphi(\zeta))\le
C_\varphi(1-|\zeta|)$. Thus $|\varphi'(\zeta)|\le
C(1-|\zeta|)^{1/L-1}$. Let $\zeta_1,\zeta_2\in\Delta$ and
$d=|\zeta_1-\zeta_2|\le1/2$. A point of modulus greater than $1-d$ is
moved radially to the circle of radius $1-d$, at cost at most
$\int_{1-d}^1C(1-r)^{1/L-1}dr=CLd^{1/L}$; the two resulting points are
at distance at most $3d$ and are joined by a segment in
$\{|\zeta|\le1-d\}$, on which $|\varphi'|\le Cd^{1/L-1}$. Hence
$|\varphi(\zeta_1)-\varphi(\zeta_2)|\le C'd^{1/L}$, and $\varphi$
extends uniquely as stated.
\end{proof}

For $y\in D$ define the inscribed complex-polydisc volume
\begin{equation}\label{eq:inscribed-volume}
 \mathcal V_D(y)=\sup\{\operatorname{vol}_{2n}(P):
 P=y+A(\Delta^n),\ A\in\operatorname{GL}(n,\C),\
 \overline{2P}\subset D\},
 \qquad 2P=y+2A(\Delta^n).
\end{equation}
It is positive, finite and lower semicontinuous, hence Borel: each
admissible matrix remains admissible when its center is moved slightly.
For a geodesic $\varphi$ put
\begin{equation}\label{eq:restriction-weight}
 W_\varphi(\zeta)=
 \frac{\delta_D(\varphi(\zeta))^\alpha
                    \mathcal V_D(\varphi(\zeta))}{(1-|\zeta|^2)^2}.
\end{equation}

\begin{lemma}[Restriction by bounded overlap]\label{lem:geodesic-restriction}
Let $D\Subset\C^n$ be convex with $C^2$ boundary and
$d\mu\asymp\delta_D^\alpha dV$, $\alpha>-1$. There is a constant $C$,
independent of $p\ge1$, the geodesic $\varphi$, and $f\in\mathcal O(D)$,
such that
\begin{equation}\label{eq:geodesic-restriction}
 \int_\Delta |f\circ\varphi|^p W_\varphi\,dA
                           \le C\int_D|f|^p\,d\mu.
\end{equation}
Moreover $\mathcal V_D(y)\ge c\delta_D(y)^{n+1}$, and for
$\varphi(0)\in K\Subset D$ one may replace $W_\varphi$ on the left
by $(1-|\zeta|^2)^{n-1+\alpha}$, with a constant depending on $K$.
\end{lemma}
\begin{proof}
Put $b=\operatorname{arctanh}(1/2)$. Every admissible $P$ lies in
$\mathcal B_D(y,b)$ by holomorphic contraction from $2P$ and the product formula
for the polydisc distance. Lemma~\ref{lem:geodesic-depth} makes
$\delta_D/\delta_D(y)$ uniformly bounded above and below on $P$.
Submean in the affine complex polydisc gives, for either sign of $\alpha$,
\[
 |f(y)|^p\le
 \frac{C}{\delta_D(y)^\alpha\operatorname{vol}(P)}
           \int_{\mathcal B_D(y,b)}|f|^p\,d\mu.
\]
Taking admissible volumes tending to their supremum proves the same
bound with $\operatorname{vol}(P)$ replaced by $\mathcal V_D(y)$.
Multiply by $W_\varphi(\zeta)$ at $y=\varphi(\zeta)$ and use Tonelli.
For every fixed $z\in D$, the set
$S_z=\{\zeta:k_D(z,\varphi(\zeta))<b\}$, if nonempty, is contained
in a disc Kobayashi ball of radius $2b$: choose $\zeta_0\in S_z$ and
use the triangle inequality and the isometry property of $\varphi$.
The invariant area of that disc ball is finite and independent of its
center. Hence
\[
 \int_{S_z}\frac{dA(\zeta)}{(1-|\zeta|^2)^2}\le C_b,
\]
which proves~\eqref{eq:geodesic-restriction}.

Taylor's theorem for signed distance supplies an admissible polydisc at
a point of depth $d$, with one complex-normal radius $cd$ and $n-1$
complex-tangential radii $c\sqrt d$, for a fixed sufficiently small $c$.
Indeed the linear term is at most $O(cd)$ and the quadratic term at most
$O(c^2d)$ even on the doubled polydisc. Its volume is $c'd^{n+1}$.
Interior balls treat the remaining compact set. The last assertion
follows from~\eqref{eq:geodesic-depth}.
\end{proof}

The following proof of Theorem~\ref{thm:general-main}(i) also gives
its lower bound without convexity.
\begin{proof}
For the lower bound choose a Carath\'eodory extremal $u(x)=0$.
An interior Schwarz estimate gives $1-|u|^2\ge3/4$ on a fixed-radius
ball about $x$, uniformly for $x\in K$. The positive superharmonic
function $1-|u|^2$ dominates $3/4$ times the harmonic function equal
to one on that inner sphere and zero on $\partial D$. Hopf's lemma,
compactness and a finite cover in $x$ give $1-|u|^2\ge c_K\delta_D$.
Consequently
$\int|u|^p\,d\mu\le C_K\int e^{-c_Kp\delta_D/2}\delta_D^\alpha dV
\le C_Kp^{-1-\alpha}$. The normalized first jet and
Lemma~\ref{lem:basic-comparison} prove the assertion.

For the upper bound write $\gamma=\gamma_D(x;v)$ and choose a complex
geodesic with $\varphi(0)=x$, $\varphi'(0)=v/\gamma$. For every jet
competitor $f$, the function $g=f\circ\varphi$ has $g(0)=0$ and
$g'(0)=\gamma^{-1}$. Lemma~\ref{lem:geodesic-restriction} and Cauchy's
coefficient inequality on circles give
\[
 \|f\|_p^p\ge c_K\gamma^{-p}
     \int_0^1r^{p+1}(1-r^2)^{n-1+\alpha}\,dr
                 \ge c_K'\gamma^{-p}p^{-(n+\alpha)}.
\]
Take the infimum and use $E_p\ge1$. This proves
the upper bound in Theorem~\ref{thm:general-main}(i). Smooth strongly
convex domains attain the upper endpoint
(Theorem~\ref{thm:regular-reverse}(ii)); in dimension one the two
endpoints coincide.
\end{proof}

\begin{proposition}[Dependence on the distance to the boundary]
\label{prop:boundary-window}
Under the hypotheses of Theorem~\ref{thm:general-main}(i), set
$d(x)=\min\{1,\delta_D(x)\}$ and $\mathfrak l(x)=\log(1/d(x))$. There is a
constant $C$, depending only on $D$, $\alpha$ and the two density
comparison constants, such that for every $x\in D$, $v\ne0$ and $p\ge2$,
\begin{align}
 (1+\alpha)\log p-(n+2+2\alpha)\mathfrak l(x)-C
 &\le p\log\frac{B_{p,\mu}(x;v)}{\gamma_D(x;v)}\notag\\
 &\le(n+\alpha)\log p+(n+1+\alpha)\mathfrak l(x)+C.
 \label{eq:boundary-window}
\end{align}
Consequently the coefficient window $[1+\alpha,n+\alpha]$ remains
valid for arbitrary moving jets $(x_p,v_p)$ with $v_p\ne0$ and
$\mathfrak l(x_p)=o(\log p)$, in the sense of lower and upper limits.
\end{proposition}
\begin{proof}
A uniform interior tangent ball at each boundary point, joined to a
fixed compact core of $D$, gives
$k_D(x,y)\le C+\tfrac12 \mathfrak l(x)+\tfrac12 \mathfrak l(y)$.
Indeed, along the radius of a tangent ball of radius $r_0$, the
distance from its center to a point of depth $t$ is
$\frac12\log((2r_0-t)/t)$; the centers form a compact subset of $D$.
For a disc-valued extremal $u$ with $u(x)=0$, Schwarz--Pick therefore
gives $1-|u(y)|^2\ge c\,d(x)d(y)$. Boundary normal coordinates, and
$\alpha>-1$, imply
\begin{equation}\label{eq:global-moment-depth}
 \int_D|u|^p\,d\mu\le C\bigl(p\,d(x)\bigr)^{-1-\alpha}.
\end{equation}
For $p\,d(x)<1$ enlarge $C$ and use $\mu(D)=1$; otherwise integrate
$e^{-cpd(x)\delta_D}\delta_D^\alpha$ in the collar, with exponential
decay on its complement.

The inscribed-polydisc construction in the proof of
Lemma~\ref{lem:geodesic-restriction}, with normal radius $c\delta_D(x)$
and tangential radii $c\sqrt{\delta_D(x)}$, has volume at least
$c\delta_D(x)^{n+1}$ and depth comparable to $\delta_D(x)$ throughout.
For $x$ in the compact core use a fixed-size polydisc. The submean
inequality gives
$E_{p,\mu}(x)^p\le C d(x)^{-n-1-\alpha}$.
Combining this with~\eqref{eq:global-moment-depth} and the normalized
jet $u/du_x(v)$ proves the lower bound in~\eqref{eq:boundary-window}.

For the upper bound take an extremal geodesic $\varphi(0)=x$ in the
direction $v$. The depth inequality in Lemma~\ref{lem:geodesic-depth}
gives $\delta_D(\varphi(\zeta))\ge\delta_D(x)(1-|\zeta|^2)/4$.
The restriction weight thus satisfies, since $n+1+\alpha>0$,
\[
 W_\varphi(\zeta)\ge c d(x)^{n+1+\alpha}
                         (1-|\zeta|^2)^{n-1+\alpha}.
\]
The circular Cauchy inequality in the proof of
Theorem~\ref{thm:general-main}(i) gives
$(\gamma_D I_{p,\mu})^p\ge
c d(x)^{n+1+\alpha}p^{-n-\alpha}$. Use $E_{p,\mu}\ge1$.
\end{proof}

We first record a finite-$p$ estimate which uses more information than
the limiting multitype. With $W_\varphi$ as in~\eqref{eq:restriction-weight},
put
\[
 \mathcal L_\varphi(r)=\int_{\partial\Delta}\log W_\varphi(r\zeta)
                                      \,\frac{|d\zeta|}{2\pi}.
\]
For every $p\ge2$ and every geodesic in a bounded $C^2$ convex domain,
\begin{equation}\label{eq:finite-p-geodesic-entropy}
 I_{p,\mu}(\varphi(0);\varphi'(0))^p
       \ge c\int_{1/2}^1r^{p+1}\exp \mathcal L_\varphi(r)\,dr.
\end{equation}
Indeed restriction and the arithmetic--geometric mean inequality on
each circle, applied to $g(\zeta)/\zeta$ with $g'(0)=1$, give this
formula directly. All logarithms are integrable on each fixed circle:
$\mathcal V_D$ has the positive lower bound from the restriction lemma
and is at most $2^{-2n}\operatorname{vol}(D)$.

\begin{lemma}[Endpoint scaling along a complex geodesic]
\label{lem:endpoint-scaling}
Let $D\Subset\C^n$ be a convex domain with smooth boundary of finite
D'Angelo type, let $\varphi:\Delta\to D$ be a complex geodesic, and let
$\zeta\in\partial\Delta$ and $q=\varphi^*(\zeta)$. Then
$q\in\partial D$, and there are $c_\zeta>0$ and $r_\zeta<1$ such that
\begin{equation}\label{eq:volume-along-geodesic}
 \mathcal V_D(\varphi(r\zeta))\ge c_\zeta
          (1-r^2)^{2+T(q)}\qquad(r_\zeta<r<1).
\end{equation}
\end{lemma}
\begin{proof}
We use~\cite{ArosioFiacchi2025}, whose standing assumption is that the
domain is convex and contains no complex affine line ($\C$-proper
in the sense of~\cite[Definition~2.1]{ArosioFiacchi2025}); a bounded
convex domain has this property. A holomorphic left inverse $u:D\to\Delta$
satisfies $u(\varphi(r\zeta))=r\zeta$, so $q$ cannot be interior, since
$|u(q)|$ would equal one. Rotating the parameter, we take $\zeta=1$.

\emph{Local finite type.} Let $\rho$ be a local defining function. A
complex line is a holomorphic curve of order one, so every line type at
$q$ is at most the D'Angelo type $L$ of $\partial D$ at $q$ (for
directions transverse to $T^{1,0}_q$ it is one). Thus, for every unit
$v\in\C^n$, some Taylor coefficient of order at most $L$ of
$\zeta\mapsto\rho(p+\zeta v)$ at $\zeta=0$ is nonzero when $p=q$. These
coefficients are continuous in $(p,v)$, so by compactness of the unit
sphere the same holds for all $p\in\partial D$ near $q$. Hence $q$ is a
point of locally finite type in the sense
of~\cite[Definition~2.3]{ArosioFiacchi2025}.

\emph{Regular contact point.} In~\cite{ArosioFiacchi2025} the
Kobayashi distance is normalized by $k_\Delta(0,r)=2\operatorname{arctanh}r$,
twice ours, and regular contact points are defined
in~\cite[Definition~4.5]{ArosioFiacchi2025} through the rays
$\tanh(t/2)$ in $\Delta$ and $-e^{-t}$ in the left half-plane
$\mathbb L$~\cite[Definition~4.4]{ArosioFiacchi2025}. A constant factor
in the distance does not change the class of complex geodesics, and
only qualitative conclusions are used. The curve
$t\mapsto\varphi(\tanh(t/2))$ has endpoint $q$, and a complex geodesic
whose endpoint is a point of locally finite type has a regular contact
point there~\cite[Remark~4.7]{ArosioFiacchi2025}; so $1$ is a regular
contact point of $\varphi$ with endpoint $q$.

\emph{Coordinates.} By~\cite[Section~5, (5.1)]{ArosioFiacchi2025}, after a
unitary change of coordinates and a translation, $q=0$, the outer unit
normal is $e_0$, and $e_1,\ldots,e_{n-1}$ is an orthonormal basis of
$T^{1,0}_q$ adapted to Yu's flag $S_d=\{v\in T^{1,0}_q:t(v)\ge d\}$, $d\ge2$, in the sense
that each $S_d$ is spanned by the basis vectors it contains. With
$m_j=t(e_j)$ a defining function has the form
\[
 \operatorname{Re}z_0+H(z_1,\ldots,z_{n-1})
   +o\Bigl(|z_0|+\sum_{j=1}^{n-1}|z_j|^{m_j}\Bigr),
\]
where $H$ is a nonnegative convex polynomial, weighted homogeneous for
the weights $1/m_j$, and non-degenerate in the sense
of~\cite{ArosioFiacchi2025}: $\{H=0\}$ contains no complex line (for
instance $H=(\operatorname{Im}z_1)^2$ is allowed). The scaling model is
$\mathcal M_H=\{\operatorname{Re}z_0+H(z')<0\}$. These are the linear multitype
coordinates of Section~\ref{sec:multitype}: the sets $S_d$ are the
complex subspaces of Lemma~\ref{lem:linear-multitype}, an adapted
basis has exactly $\dim_\C S_d$ vectors with $m_j\ge d$, and by that
lemma the same count holds for $(\tau_j)$. Thus $(m_j)$ is a permutation
of $(\tau_j(q))$, and $\sum_j2/m_j=T(q)$.

\emph{Scaling.} Put
\[
 \mathcal S_t(z_0,z')=(t^{-1}z_0,t^{-1/m_1}z_1,\ldots,
                                    t^{-1/m_{n-1}}z_{n-1}),\qquad t>0.
\]
For every sequence $t_k\downarrow0$, the scaling
theorem~\cite[Theorem~6.2]{ArosioFiacchi2025} (Gaussier's
scaling~\cite{Gaussier1997}), applied with $\lambda_k=1/t_k$, for which
its linear maps $A_k$ are our $\mathcal S_{t_k}$, shows that the
domains $\mathcal S_{t_k}D$ converge to $\mathcal M_H$ in the local Hausdorff
topology, and in particular
each compact subset of $\mathcal M_H$ is contained in $\mathcal S_{t_k}D$ for
large $k$. Since the sequence is arbitrary, each compact subset of
$\mathcal M_H$ lies in $\mathcal S_tD$ for all sufficiently small $t$.

\emph{The scaled geodesic.} Let $\mathscr C(\zeta)=(\zeta-1)/(\zeta+1)$,
which maps $\Delta$ onto the left half-plane $\mathbb L$ and $1$ to
$0$, and put $\psi=\varphi\circ\mathscr C^{-1}$. Since
$\mathscr C^{-1}(-e^{-t})=\tanh(t/2)$, the rays
of~\cite[Definition~4.4]{ArosioFiacchi2025} for $\varphi$ at $1$ and
for $\psi$ at $0$ coincide, so $0$ is a regular contact point of $\psi$
with endpoint $q$. By~\cite[Proposition~4.9]{ArosioFiacchi2025} the
function $\langle\psi',e_0\rangle$ has a positive nontangential limit
$\kappa$ at $0$, and
by~\cite[Proposition~6.6]{ArosioFiacchi2025}, which concerns
holomorphic maps from $\mathbb L$ with a regular contact point at $0$
and is applied with $\lambda_k=1/t_k$, so that
$A_k\psi(\lambda_k^{-1}\,\cdot\,)=\mathcal S_{t_k}\psi(t_k\,\cdot\,)$,
we have $\mathcal S_{t_k}\psi(t_k\,\cdot)\to(\kappa\,\cdot\,,0)$
uniformly on compact subsets of $\mathbb L$ for every sequence
$t_k\downarrow0$. As $\varphi(1-t)=\psi(tw_t)$ with
$w_t=-1/(2-t)\to-1/2$, we obtain
\begin{equation}\label{eq:scaled-geodesic-normal}
 \mathcal S_t\varphi(1-t)\longrightarrow y_0=(-\kappa/2,0)\quad(t\downarrow0),\qquad
 \kappa>0.
\end{equation}

\emph{Volume.} Choose $r_1>0$ with $y_0+4r_1\overline\Delta^n
\subset \mathcal M_H$. For small $t$ this compact polydisc lies in $\mathcal S_tD$, and
$\mathcal S_t\varphi(1-t)$ lies within $r_1$ of $y_0$ in each coordinate.
Then $P_t=\mathcal S_t\varphi(1-t)+r_1\Delta^n$ satisfies
$\overline{2P_t}\subset y_0+3r_1\overline\Delta^n\subset \mathcal S_tD$.
Thus $\mathcal S_t^{-1}P_t$ is admissible in~\eqref{eq:inscribed-volume} at
$\varphi(1-t)$; unitary changes preserve this class and volume, and
$\det_\R \mathcal S_t^{-1}=t^{2+T(q)}$. With $t=1-r\ge(1-r^2)/2$ this
gives~\eqref{eq:volume-along-geodesic}.
\end{proof}

\begin{proof}[Proof of Theorem~\ref{thm:general-main}(ii)]
First take $a=0$. By Lemmas~\ref{lem:holder}
and~\ref{lem:endpoint-scaling}, $\varphi^*$ is continuous, and for each
$\zeta\in\partial\Delta$ one has $\varphi^*(\zeta)\in\partial D$
and~\eqref{eq:volume-along-geodesic} for $r>r_\zeta$.

Write $s=1-r^2$. The depth lemma and
\eqref{eq:volume-along-geodesic} give, for every $\zeta\in\partial\Delta$,
\[
 \liminf_{r\to1}\frac{\log W_\varphi(r\zeta)}{\log(1/s)}
                     \ge-\alpha-T(\varphi^*(\zeta)).
\]
The normalized logarithms are bounded above and below independently
of $r$ close to one and of $\zeta$. This follows from
$c\delta_D^{n+1}\le\mathcal V_D\le2^{-2n}\operatorname{vol}(D)$ and
$\delta_D\circ\varphi\asymp s$. The function $T$ is Borel on
$\partial D$, so $T\circ\varphi^*$ is Borel, since $\varphi^*$ is
continuous. Indeed, let $t_q$ be the line type at $q$. By
\cite[Section~5]{ArosioFiacchi2025} (Yu's theorem), at every
$q\in\partial D$ each $S_m(q)=\{v\in T^{1,0}_q:t_q(v)\ge m\}$ is a
complex subspace, and $v\in S_m(q)$ if and only if the Taylor
coefficients of orders below $m$ of $\zeta\mapsto r(q+\zeta v)$ vanish;
these coefficients are continuous in $(q,v)$. Each set
$\mathcal C_{m,j}=\{q:\dim_\C S_m(q)\ge j\}$ is closed: if $q_k\to q$ and
$P_k\subset S_m(q_k)$ are $j$-dimensional, then, in a local smooth
frame of $T^{1,0}\partial D$, a subsequence of $(P_k)$ converges in the
compact Grassmannian to a $j$-plane $P\subset T^{1,0}_q$, every $v\in P$
is a limit of $v_k\in P_k$, and continuity gives $P\subset S_m(q)$.
So $q\mapsto\dim S_m(q)=\sum_j\mathbf 1_{\mathcal C_{m,j}}(q)$ is upper
semicontinuous, hence Borel. By Lemma~\ref{lem:linear-multitype},
$\#\{j:\tau_j(q)\ge m\}=\dim S_m(q)$, so, with $L$ the D'Angelo type of $\partial D$,
$T(q)=\sum_{m=2}^{L}\frac2m(\dim S_m(q)-\dim S_{m+1}(q))$ is Borel.
Also $0<T\le n-1$ when $n>1$, and $T=0$ when $n=1$. Since
$r\mapsto\log W_\varphi(r\zeta)$ is only lower semicontinuous, choose a
sequence $r_k\to1$ realizing $\liminf_{r\to1}\mathcal L_\varphi(r)/\log(1/s)$ and apply
Fatou's lemma to the sequence $\log W_\varphi(r_k\zeta)/\log(1/s_k)$, after
adding a constant lower bound; this yields
\[
 \liminf_{r\to1}\frac{\mathcal L_\varphi(r)}{\log(1/s)}
       \ge-\alpha-\int_{\partial\Delta}T(\varphi^*)\,\frac{|d\zeta|}{2\pi}.
\]
Put $N_\varphi(0)=1+\alpha+\int T(\varphi^*)|d\zeta|/(2\pi)$, the right
side of~\eqref{eq:general-averaged-upper} at $a=0$.
For every $\theta>0$, at all sufficiently small depths one has
$\exp \mathcal L_\varphi(r)\ge s^{N_\varphi(0)-1+\theta}$.
Integrate~\eqref{eq:finite-p-geodesic-entropy} only over
$1/p\le s\le2/p$. It gives
$I_{p,\mu}(\varphi(0);\varphi'(0))^p\ge c_\theta
p^{-N_\varphi(0)-\theta}$ for large $p$.
Since $\gamma_D(\varphi(0);\varphi'(0))=1$ and $E_p\ge1$, this is
\eqref{eq:general-averaged-upper} at zero after $\theta\downarrow0$.
Reparametrize by $\mob_{-a}$ to prove the assertion at $a$; its boundary
map carries normalized arclength to $\omega_a$. Homogeneity removes $\eta$.
\end{proof}

\section{Equality at coercive boundary points}
\label{sec:regular-reverse}

This section proves Theorem~\ref{thm:regular-reverse} in every
dimension $n\ge2$. Lemma~\ref{lem:stolz-transfer} reduces the lower
bound, for any holomorphic $u:D\to\Delta$ with $u(x)=0$,
$du_x(v)=\gamma_D(x;v)$, to bounds for $u_*\mu$ in Stolz regions;
Lemmas~\ref{lem:model-differential} and~\ref{lem:nontangential-fibers}
prove them at coercive points when $u$ is a left inverse with affine
fibers, which exists by Lemma~\ref{lem:affine-inverse}.

Throughout, $\mathbb L=\{\operatorname{Re}\sigma<0\}$ and
$\mathscr C(\lambda)=(\lambda-1)/(\lambda+1)$.
Stolz regions always mean~\eqref{eq:stolz-region}, with depth
$s(\lambda)=1-|\lambda|^2$. Under the Cayley map, a fixed Stolz
region at $1$ lies in a fixed sector at $0$ of $\mathbb L$, since
\[
 \frac{|\mathscr C(\lambda)|}{-\operatorname{Re}\mathscr C(\lambda)}
 =\frac{|1-\lambda|\,|1+\lambda|}{1-|\lambda|^2}<M
 \qquad(\lambda\in\mathcal A_M(1)).
\]
Conversely, the inverse Cayley map sends a truncated sector
$|\sigma|\le M(-\operatorname{Re}\sigma)$ near $0$ into a Stolz region
at $1$, with a possibly larger aperture. This follows from the same
identity and $|1+\lambda|/(1+|\lambda|)\to1$ as $\lambda\to1$.
For a measure $\nu$ on $\Delta$, a Borel set $W\subset\Delta$ and a
Borel function $w\ge0$, we write $\nu\le w\,dA$ on $W$ if
$\nu(S)\le\int_Sw\,dA$ for all Borel $S\subset W$. If this holds on
each set $W_k$ of a countable family of Borel sets, it holds on their
union: apply it to the disjoint sets
$W_k\setminus\bigcup_{i<k}W_i$ and add. Since $\Delta$ is second countable, every
union of open sets is a countable subunion, so the same holds for
every family of open sets.

\emph{Multitype coordinates.} Let $q\in\partial D$ be a point of a
smooth convex finite-type boundary. By the normal
form~\cite[(5.1)]{ArosioFiacchi2025} (see the proof of
Lemma~\ref{lem:endpoint-scaling}) there are a unitary frame
$e_0,\ldots,e_{n-1}$, with $e_0$ the outer unit normal at $q$, and
integers $m_1,\ldots,m_{n-1}\ge2$, a permutation of
$(\tau_1(q),\ldots,\tau_{n-1}(q))$, such that in the coordinates
$z=q+\sum_{j=0}^{n-1}z_je_j$ centered at $q$ a local defining function
is
\begin{equation}\label{eq:multitype-normal-form}
 r(z)=\operatorname{Re}z_0+H(z')+R(z),\qquad
 R(z)=o\Bigl(|z_0|+\sum_{j=1}^{n-1}|z_j|^{m_j}\Bigr),
\end{equation}
where $H\ge0$ is a convex polynomial with
$H(t^{1/m_1}z_1,\ldots,t^{1/m_{n-1}}z_{n-1})=tH(z')$ for $t>0$. Thus
$T(q)=\sum_j2/m_j$. Multiplying $r$ by a positive smooth function
equal to one at $q$ changes $R$ only, so $H$ does not depend on the
defining function, and in given adapted coordinates $H$ is determined
by $D$: if $H_1,H_2$ both occur, then $H_1-H_2$ is weighted homogeneous
and $o(\sum_j|z_j|^{m_j})$ on $\{z_0=0\}$, hence zero.

\emph{Affine invariance.} Coercivity is a complex-affine invariant of
$(D,q)$. After translations it suffices to consider $D'=AD$, $q'=Aq$,
with $A\in\operatorname{GL}(n,\C)$. As $r\circ A^{-1}$ defines $D'$ and
$r(q+\zeta v)=(r\circ A^{-1})(q'+\zeta Av)$, $A$ maps $T^{1,0}_q$ onto
$T^{1,0}_{q'}$ and preserves line types, hence maps each $S_d(q)$ onto
$S_d(q')$; so the multitypes agree, and we order adapted coordinates
$(z_0,z')$ at $q$ and $(w_0,w')$ at $q'$ so that the weights $m_j$
agree. Write $z=A^{-1}w$. Since $A^{-1}$ maps $\{w_0=0\}$ onto
$\{z_0=0\}$, the normal coordinate is $z_0=w_0/a$, $a\in\C^*$, with no
$w'$-terms, and comparing the linear part $\operatorname{Re}(w_0/a)$ of
$r\circ A^{-1}$ with the outer normal $e_0$ at $q'$ gives $a>0$. Since
$A^{-1}$ maps the $k$th adapted basis vector at $q'$, $k\ge1$, into
$S_{m_k}(q)$, we have $z'=A_0w'+Ow'+bw_0$, where
$A_0$ is block diagonal for $T^{1,0}=\bigoplus_d(S_d\ominus S_{d+1})$
and $O$ sends $w_k$ only into coordinates $z_j$ with $m_j>m_k$; thus
$A_0+O$ is block triangular and invertible, and so is $A_0$. Give $w_0$ the weight $1$ and $w_k$ the weight $1/m_k$. A
monomial of $H$ has weight one, and replacing one of its factors
$z_j$ (or $\bar z_j$) by a term of $Ow'$ or $bw_0$ replaces the weight
$1/m_j$ by $1/m_k>1/m_j$ or by $1>1/m_j$. Terms of weight above one are
$o(|w_0|+\sum_k|w_k|^{m_k})$, so
$H(z')=H(A_0w')+o(|w_0|+\sum_k|w_k|^{m_k})$; the same bound holds for
$R(A^{-1}w)$, because $|z_0|\lesssim|w_0|$ and
$|z_j|^{m_j}\lesssim|w_0|^{m_j}+\sum_{m_k\le m_j}|w_k|^{m_j}$. Thus
$a\,r\circ A^{-1}=\operatorname{Re}w_0+aH(A_0w')+o(|w_0|+\sum_k|w_k|^{m_k})$,
the model at $q'$ is $aH\circ A_0$, and it vanishes only at $0$ if and
only if $H$ does. For two adapted frames at $q$, $A$ is unitary and
$H$ is replaced by $H\circ A_0$. A biholomorphism, in contrast, can add a
term of weight one to the normal coordinate, so coercivity is not a
biholomorphic invariant: $w_0=z_0-z_1^2/2$ turns
$\operatorname{Re}z_0+(\operatorname{Im}z_1)^2$ into
$\operatorname{Re}w_0+|z_1|^2/2$.

In $\C^2$, $H$ is the degree-$2m$ Taylor term of
$w\mapsto r(q+we_1)$, where $2m=m_1$ is the tangential type, which
gives the equivalence of coercivity
with~\eqref{eq:real-tangential-contact}. If $q$ is coercive,
compactness of $\{\sum_j|z_j|^{m_j}=1\}$ and weighted homogeneity give
$c_H>0$ with
\begin{equation}\label{eq:H-coercive}
 H(z')\ge c_H\sum_{j=1}^{n-1}|z_j|^{m_j}\qquad(z'\in\C^{n-1}).
\end{equation}

The non-degeneracy of $H$ in~\cite{ArosioFiacchi2025}, that $\{H=0\}$
contains no complex line, holds at every point of finite type, coercive
or not; coercivity is the stronger condition $\{H=0\}=\{0\}$.

\begin{figure}[tbp]
\centering
\includegraphics[width=\textwidth]{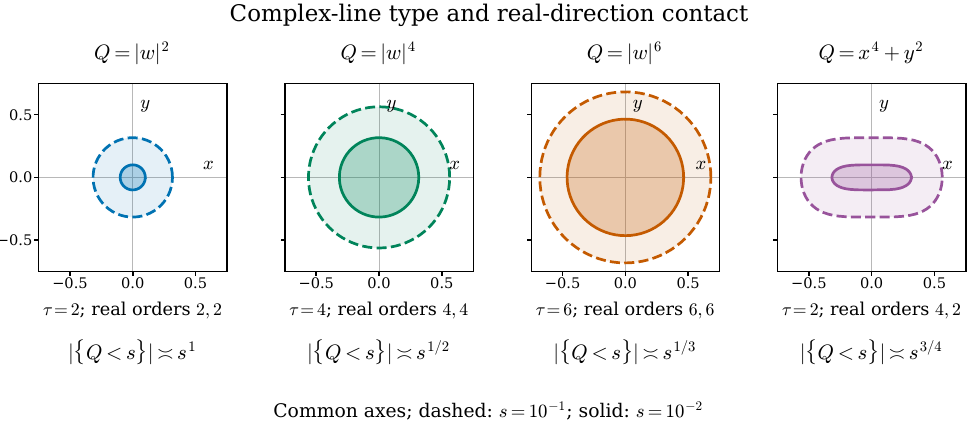}
\caption{Tangential sublevel sets in one complex variable $w=x+iy$.
Here $\tau$ is the vanishing order of the full real Taylor expansion
in the complex variable, not the order along a chosen real ray.
The first three profiles are coercive; the last has a non-coercive
quadratic model and quartic contact in its null direction. Real orders
refer to the displayed $x$- and $y$-axes. The areas follow by the
changes of variables $w=s^{1/\tau}W$ and
$(x,y)=(s^{1/4}X,s^{1/2}Y)$, respectively. All panels use the same
coordinate scale.}
\label{fig:contact-orders}
\end{figure}

Figure~\ref{fig:contact-orders} separates complex-line type from contact
along individual real directions. In the last profile the quadratic
model is $H=y^2$, although the full tangential profile $Q=x^4+y^2$
has a bounded sublevel set. Its area exponent is $3/4$, not $2/\tau=1$.
Such a local area exponent alone does not determine $N$: the global
correction in Section~\ref{sec:noncoercive} must also be taken into account.

\emph{Type-two points and strong linear convexity.} Take $r$ smooth,
and for $z'\in T^{1,0}_q=\{z_0=0\}$, the complex tangent space at $q$,
let $\mathcal L_q$ and $\mathcal Q_q$ be the Levi form and the complex
Hessian of $r$ at $q$, displayed below. The quadratic part of
$\lambda\mapsto r(q+\lambda z')$ is
$|\lambda|^2\mathcal L_q(z')+\operatorname{Re}(\lambda^2\mathcal Q_q(z'))$;
by convexity it is nonnegative, so $|\mathcal Q_q|\le\mathcal L_q$, and
it vanishes identically exactly when $\mathcal L_q(z')=0$. Hence
$m_j=2$ for all $j$ if and only if $q$ is strongly pseudoconvex. In this
case $H$ is a quadratic form, and since $R(0,z')=o(|z'|^2)$ and the
linear term of $r$ vanishes on $T^{1,0}_q$, $H$ is the second-order
Taylor polynomial of $z'\mapsto r(0,z')$:
\[
 H(z')=\mathcal L_q(z')+\operatorname{Re}\mathcal Q_q(z'),\qquad
 \mathcal L_q(z')=\sum_{j,k\ge1}\frac{\partial^2r(q)}{\partial z_j\partial\bar z_k}z_j\bar z_k,\qquad
 \mathcal Q_q(z')=\sum_{j,k\ge1}\frac{\partial^2r(q)}{\partial z_j\partial z_k}z_jz_k .
\]
As $\min_\theta H(e^{i\theta}z')=\mathcal L_q(z')-|\mathcal Q_q(z')|$,
the point $q$ is coercive if and only if
$|\mathcal Q_q(z')|<\mathcal L_q(z')$ for $z'\ne0$, that is, if and
only if the condition in Lempert's definition of strongly linearly
convex domains~\cite{Lempert1984Varna,KosinskiWarszawski2013} holds at
$q$: the Levi form dominates the complex Hessian on $T^{1,0}_q$. The
condition does not depend on $r$, since multiplying $r$ by a positive
function multiplies both forms on $T^{1,0}_q$ by its value at $q$. The
non-coercive type-two points are those where
$|\mathcal Q_q(z')|=\mathcal L_q(z')$ for some $z'\ne0$, as for
$(\operatorname{Im}z_1)^2=\frac12|z_1|^2-\frac12\operatorname{Re}z_1^2$.
Thus a smooth convex domain with only strongly pseudoconvex coercive
boundary points is strongly linearly convex, whereas
Theorem~\ref{thm:regular-reverse}(iii) asks for coercivity only on a set
of positive measure of one boundary circle.

\emph{The domain with boundary segments in the introduction.} Write
$w=x+iy$, $c(x)=\chi(|x|-1/2)$ and
$r=|z|^2+y^2+c(x)-1$. By Section~\ref{sec:flat-convex}, $c$ is smooth,
convex, zero exactly on $[-1/2,1/2]$, with $c'\ne0$ for $|x|>1/2$,
$0\le c''\le1$, and $c(x)\to\infty$ as $|x|\to\infty$; so $D=\{r<0\}$ is bounded and
convex, and $r(0)=-1$ gives $dr(Z)[Z]\ge1$ on $\{r=0\}$, so $\partial D$
is smooth. As $r_z=\bar z$ and $r_w=c'(x)/2-iy$, the complex tangent
space at a boundary point is $\{(a,b):\bar za+(c'/2-iy)b=0\}$, and
\[
 \mathcal L(a,b)=|a|^2+\tfrac14(2+c'')|b|^2,\qquad
 \mathcal Q(a,b)=\tfrac14(c''-2)b^2,\qquad
 \mathcal L-|\mathcal Q|=|a|^2+\tfrac12c''|b|^2 .
\]
Thus $\mathcal L\ge|a|^2+|b|^2/2$: every boundary point is strongly
pseudoconvex, of D'Angelo type $2$. A point is non-coercive if and only
if $\mathcal L=|\mathcal Q|$ at some tangent $(a,b)\ne0$, that is, $a=0$,
$c''(x)=0$ and $(c'/2-iy)b=0$ with $b\ne0$; this means $y=0$,
$|x|\le1/2$, hence $|z|=1$. At these points $(\zeta,x)$ the tangent space
is $\{0\}\times\C$, and the model of $r/2$ is
$\frac12(\mathcal L+\operatorname{Re}\mathcal Q)(0,b)
=\frac14(|b|^2-\operatorname{Re}b^2)=\frac12(\operatorname{Im}b)^2$. On
the segments $\{(0,x\pm i):|x|\le1/2\}\subset\partial D$ one has
$r_w=\mp i$, tangent space $\C\times\{0\}$, $\mathcal Q=0$, and model
$|a|^2/2$. Since $r(\zeta,w)=y^2+c(x)$ for $|\zeta|=1$, the complex face
$\overline D\cap(\{\zeta\}\times\C)$ at $(\zeta,x)$ is
$\{\zeta\}\times[-1/2,1/2]$; since $r(z,x\pm i)=|z|^2$ for
$|x|\le1/2$, the complex face at $(0,x\pm i)$ is trivial. In particular
$D$ is neither strictly nor strongly linearly convex.

\begin{lemma}[Transfer from Stolz regions]\label{lem:stolz-transfer}
Let $D\Subset\C^n$ be convex with $C^2$ boundary,
$d\mu\asymp\delta_D^\alpha dV$ a probability measure, $\alpha>-1$,
$x\in D$, $v\ne0$, and let $u:D\to\Delta$ be holomorphic with $u(x)=0$
and $du_x(v)=\gamma_D(x;v)$. Put $\nu=u_*\mu$. Let
$\mathcal T\subset[0,\infty)$ be finite, $\mathsf C\subset\partial\Delta$ a
Borel set and $\beta:\mathsf C\to\mathcal T$ measurable. Suppose that for every
$\zeta\in\mathsf C$ and every $M\ge2$ there are $c,s_0>0$ with
\begin{equation}\label{eq:stolz-hypothesis}
 \nu\le c\,s^{\alpha+\beta(\zeta)}\,dA\quad\text{on }
 \mathcal A_M(\zeta)\cap\{s<s_0\}.
\end{equation}
Then
\[
 \liminf_{p\to\infty}\frac p{\log p}
 \log\frac{B_{p,\mu}(x;v)}{\gamma_D(x;v)}
 \ge1+\alpha+\int_{\mathsf C}\beta\,d\omega_0 .
\]
\end{lemma}

\begin{proof}
Write $dm=d\omega_0$, $\Theta=\max\mathcal T$, and normalize
$\gamma_D(x;v)=1$. Fix integers $M,j\ge2$ and put
$W(\zeta)=\mathcal A_M(\zeta)\cap\{s<1/j\}$. For $\theta\in\mathcal T$ let
$E_\theta$ be the set of $\zeta\in\partial\Delta$ with
$\nu\le js^{\alpha+\theta}dA$ on $W(\zeta)$.

\emph{$E_\theta$ is closed.} Let $\zeta_k\to\zeta$ with
$\zeta_k\in E_\theta$. Since $W(\zeta_k)=(\zeta_k/\zeta)W(\zeta)$, every
compact subset of $W(\zeta)$ lies in $W(\zeta_k)$ for large $k$, so the
bound holds on every compact subset of $W(\zeta)$, hence on $W(\zeta)$
by inner regularity.

\emph{The weight.} Let $\hat\beta(\xi)=\max\{\theta:\xi\in E_\theta\}$, and
$\hat\beta(\xi)=0$ if $\xi$ lies in no $E_\theta$. This is a Borel function with
values in $[0,\Theta]$. Let $F$ be its Herglotz integral, so $F$ is
holomorphic on $\Delta$, $\operatorname{Re}F=\mathcal P[\hat\beta]$ and
$F(0)=\langle\hat\beta\rangle:=\int\hat\beta\,dm$. The function
\[
 f_p=u\exp\Bigl\{\frac{\log p}p\bigl(F(u)-\langle\hat\beta\rangle\bigr)\Bigr\}
\]
is holomorphic on $D$, $f_p(x)=0$ and $d(f_p)_x(v)=1$; moreover
$|f_p|^p=|u|^pp^{\mathcal P[\hat\beta](u)-\langle\hat\beta\rangle}$. Hence
\begin{equation}\label{eq:transfer-norm}
 I_{p,\mu}(x;v)^p\le\|f_p\|_p^p
 =p^{-\langle\hat\beta\rangle}\int_\Delta|\lambda|^pp^{\mathcal P[\hat\beta](\lambda)}\,d\nu(\lambda).
\end{equation}

\emph{Harmonic measure outside cones.} Let $|\lambda|\ge1/2$,
$d=1-|\lambda|$ and $\xi=e^{it}\lambda/|\lambda|$, $|t|\le\pi$. Then
$|\xi-\lambda|^2=d^2+|\lambda||\xi-\lambda/|\lambda||^2\ge2t^2/\pi^2$, since
$|e^{it}-1|=2|\sin(t/2)|\ge2|t|/\pi$. On
$\{\xi:|\xi-\lambda|\ge Md\}$ therefore
$|\xi-\lambda|^2\ge\frac12(M^2d^2+2t^2/\pi^2)$, and with
$\mathcal P_\lambda=(1-|\lambda|^2)/|\xi-\lambda|^2\le2d/|\xi-\lambda|^2$,
\[
 \omega_\lambda\{\xi:|\xi-\lambda|\ge Md\}\le\frac1{2\pi}\int_{\R}
 \frac{4d\,dt}{M^2d^2+2t^2/\pi^2}
 =\frac{4d}{2\pi}\cdot\frac{\pi^2}{\sqrt2\,Md}=\frac{C_0}M,
 \qquad C_0=\sqrt2\,\pi .
\]

\emph{Splitting.} Let $S_\theta=\bigcup_{\zeta\in E_\theta}W(\zeta)$;
then $\nu\le js^{\alpha+\theta}dA$ on $S_\theta$. Points with $s<1/j$
satisfy $|\lambda|\ge1/2$. Let $\lambda$ satisfy $s(\lambda)<1/j$.
If $\lambda\notin S_\theta$, then $\lambda\notin\mathcal A_M(\zeta)$, that is
$|\lambda-\zeta|\ge M(1-|\lambda|)$, for every $\zeta\in E_\theta$.
Suppose $\lambda\in S:=\bigcup_\theta S_\theta$ and let $\theta_*(\lambda)$
be the largest $\theta$ with $\lambda\in S_\theta$. Then $\hat\beta\le\theta_*$
off $E_>=\bigcup_{\theta>\theta_*}E_\theta$, $\hat\beta\le\Theta$ on $E_>$, and $E_>$
lies in the set $\{\xi:|\xi-\lambda|\ge M(1-|\lambda|)\}$. Hence
$\mathcal P[\hat\beta](\lambda)\le\theta_*(\lambda)+\Theta C_0/M$. If
$\lambda\notin S$, the same argument gives
$\mathcal P[\hat\beta](\lambda)\le\Theta C_0/M$.

Split the integral in~\eqref{eq:transfer-norm} accordingly. On
$\{\lambda\in S:\theta_*(\lambda)=\theta\}\subset S_\theta$, since
$\int_\Delta|\lambda|^ps^{\alpha+\theta}dA\le Cp^{-1-\alpha-\theta}$,
the contribution is at most
\[
 jp^{\theta+\Theta C_0/M}\int_\Delta|\lambda|^ps^{\alpha+\theta}\,dA
 \le C_{j}\,p^{-1-\alpha+\Theta C_0/M}.
\]
On $\{s<1/j\}\setminus S$, \eqref{eq:global-moment-depth} gives the
bound $p^{\Theta C_0/M}\int_D|u|^p\,d\mu\le Cp^{-1-\alpha+\Theta C_0/M}$.
On $\{s\ge1/j\}$ the integrand is at most
$(1-1/j)^{p/2}p^{\Theta}$. Summing over the finitely many $\theta$,
\[
 I_{p,\mu}(x;v)^p\le C_{M,j}\,p^{-1-\alpha-\langle\hat\beta\rangle+\Theta C_0/M}.
\]
By Lemma~\ref{lem:basic-comparison},
$\liminf(p/\log p)\log(B_{p,\mu}/\gamma_D)\ge1+\alpha+\langle\hat\beta\rangle-\Theta C_0/M$.

\emph{Limits.} Write $\hat\beta=\hat\beta_{M,j}$. For fixed $M$, $W(\zeta)$ decreases
as $j$ increases, so each $E_\theta$ and hence $\hat\beta_{M,j}$ increases with
$j$. If $\zeta\in\mathsf C$, then~\eqref{eq:stolz-hypothesis} with this $M$
gives $\zeta\in E_{\beta(\zeta)}$ for all $j\ge\max\{c,1/s_0\}$, so
$\hat\beta_{M,j}(\zeta)\ge\beta(\zeta)$ for large $j$. Monotone convergence gives
$\liminf_j\langle\hat\beta_{M,j}\rangle\ge\int_{\mathsf C}\beta\,dm$. Let $j\to\infty$ and then
$M\to\infty$.
\end{proof}

\begin{lemma}[Differentials of left inverses on the model]
\label{lem:model-differential}
Let $H:\C^{n-1}\to[0,\infty)$ be continuous with
$H(t^{1/m_1}z_1,\ldots,t^{1/m_{n-1}}z_{n-1})=tH(z')$ for $t>0$, where
all $m_j>1$, and let $\mathcal M_H=\{\operatorname{Re}z_0+H(z')<0\}$. If
$\kappa>0$ and $U:\mathcal M_H\to\mathbb L$ is holomorphic with $U(\kappa\sigma,0)=\sigma$
for $\sigma\in\mathbb L$, then $dU_{(\kappa\sigma,0)}=\kappa^{-1}dz_0$ for every
$\sigma\in\mathbb L$.
\end{lemma}

\begin{proof}
Differentiating $U(\kappa\sigma,0)=\sigma$ gives
$\partial_0U(\kappa\sigma,0)=1/\kappa$. Fix $\sigma_0$, put
$x_0=-\operatorname{Re}\sigma_0>0$, fix $j\ge1$, and let
$c=\partial_jU(\kappa\sigma_0,0)$. For $w\in\C$ and $0<\eta<1$ define,
for $\xi\in\mathbb L$,
\[
 \varpi(\xi)=(1-\eta)\xi-\eta x_0+i\eta\operatorname{Im}\sigma_0,\qquad
 g(\xi)=\frac{\xi-\sigma_0}{\xi+\bar\sigma_0},\qquad
 f(\xi)=\kappa\varpi(\xi)e_0+wg(\xi)e_j .
\]
Then $\varpi(\sigma_0)=\sigma_0$, $\varpi'=1-\eta$,
$\operatorname{Re}\varpi<-\eta x_0$, $g(\sigma_0)=0$,
$g'(\sigma_0)=-1/(2x_0)$ and $|g|<1$ on $\mathbb L$, since
$-\bar\sigma_0$ is the reflection of $\sigma_0$ in $i\R$. By weighted
homogeneity $H(wg\,e_j)\le C_j|w|^{m_j}$ with
$C_j=\max_\vartheta H(e^{i\vartheta}e_j)$. Choose
$\eta=2C_j|w|^{m_j}/(\kappa x_0)+|w|^{m_j}$; for small $|w|$ this lies in
$(0,1)$, and $f(\mathbb L)\subset \mathcal M_H$. The holomorphic self-map
$U\circ f$ of $\mathbb L$ fixes $\sigma_0$, so by Schwarz--Pick
$|(U\circ f)'(\sigma_0)|\le1$, that is,
\[
 \Bigl|1-\eta-\frac{c\,w}{2x_0}\Bigr|\le1 .
\]
If $c\ne0$, take $w=-\tau\bar c/|c|$ with $\tau>0$ small: this gives
$|c|\tau/(2x_0)\le\eta=O(\tau^{m_j})$, impossible as $\tau\to0$
because $m_j>1$. Hence $c=0$.
\end{proof}

Every complex geodesic of a bounded convex domain has a holomorphic left
inverse (Section~\ref{sec:geodesic-restriction}), and also one with
affine fibers. On strongly convex domains the left inverses constructed
by Lempert have affine fibers~\cite{Lempert81}
(see~\cite[Section~1]{Zwonek2025LeftInverses}). On every bounded convex
domain, the theorem of Royden and Wong in the form
of~\cite[Theorem~2.1]{KosinskiMcCarthy2019} gives a nonzero
$\mathfrak d\in H^1(\Delta,\C^n)$ (the dual map) and a holomorphic left
inverse $u_0$ of $\varphi$ with
$(z-\varphi(u_0(z)))\cdot\mathfrak d(u_0(z))=0$ for $z\in D$, the dot
being complex bilinear; thus the fiber of $u_0$ over $\lambda$ lies in
the affine hyperplane $\varphi(\lambda)+\ker\mathfrak d(\lambda)$
whenever $\mathfrak d(\lambda)\ne0$. We use the following refinement,
due to Zwonek~\cite[Theorem~2.1]{Zwonek2025LeftInverses}: from any
holomorphic left inverse $G$ one obtains a left inverse with affine
fibers whose differential along $\varphi$ is $G'\circ\varphi$. We state
it for an arbitrary holomorphic $\varphi$ with $G\circ\varphi=\Id$ and
include a proof for completeness.

\begin{lemma}[Left inverses with affine fibers]\label{lem:affine-inverse}
Let $D\subset\C^n$ be a bounded convex domain, let $\varphi:\Delta\to D$
be holomorphic, and let $G:D\to\Delta$ be holomorphic with
$G\circ\varphi=\Id_\Delta$. Write $G'=(\partial_1G,\ldots,\partial_nG)$,
so that $dG_z(v)=G'(z)\cdot v$ with the complex bilinear pairing, and
put
\[
 g_z(\lambda)=(z-\varphi(\lambda))\cdot G'(\varphi(\lambda))
 \qquad(z\in D,\ \lambda\in\Delta).
\]
For every $z\in D$ the function $g_z$ has exactly one zero
$u(z)\in\Delta$, and this zero is simple. The map $u:D\to\Delta$ is
holomorphic, $u\circ\varphi=\Id_\Delta$, and for every $\lambda\in\Delta$
\begin{equation}\label{eq:affine-fibers}
 u^{-1}(\lambda)=D\cap\{z:(z-\varphi(\lambda))\cdot h(\lambda)=0\},
 \qquad h(\lambda)=du_{\varphi(\lambda)}=G'(\varphi(\lambda));
\end{equation}
$h:\Delta\to(\C^n)^*$ is holomorphic and $h\cdot\varphi'=1$.
\end{lemma}

\begin{proof}
Put $x_0=\varphi(0)$ and $h=G'\circ\varphi$. Differentiating
$G\circ\varphi=\Id$ gives $h\cdot\varphi'=1$, so $h$ has no zeros, and
$(z,\lambda)\mapsto g_z(\lambda)=h(\lambda)\cdot(z-\varphi(\lambda))$ is
holomorphic on $D\times\Delta$.

\emph{Step 1: the sets $\mathcal B_\varrho$.} For $0<\varrho\le1$ let
$\mathcal B_\varrho$ be the set of points $\psi(t)$, where $\psi:\Delta\to D$
is holomorphic, $\psi(0)=x_0$ and $t\in[0,\varrho)$; rotating the parameter,
every $t$ with $|t|<\varrho$ may be used. (a) $\mathcal B_\varrho$ is convex: if
$z_i=\psi_i(t_i)$, $i=1,2$, and $\varrho_1\in(\max\{t_1,t_2\},\varrho)$, $\varrho_1>0$,
the maps $\psi_i(t_i\lambda/\varrho_1)$ take the value $x_0$ at $0$ and $z_i$
at $\varrho_1$, and every convex combination of them maps $\Delta$ into $D$.
(b) $\mathcal B_\varrho$ is open: if $z=\psi(t)$, $\varrho_1\in(t,\varrho)$, $\varrho_1>0$, and
$d$ is the distance from $\psi(\{|\xi|\le t/\varrho_1\})$ to $\C^n\setminus D$,
then $\psi(t\lambda/\varrho_1)+(\lambda/\varrho_1)(z'-z)$ shows $z'\in\mathcal B_\varrho$
for $|z'-z|<\varrho_1d$. (c) $\mathcal B_\varrho$ increases with $\varrho$ and
$\bigcup_{\varrho<1}\mathcal B_\varrho=D$: use a Riemann map of the bounded convex
planar domain $D\cap(x_0+\C(z-x_0))$ centered at $x_0$. (d) $|G|<\varrho$ on
$\mathcal B_\varrho$, by Schwarz's lemma for $G\circ\psi$. (e) $\varphi(\lambda)
\in\mathcal B_\varrho$ for $|\lambda|<\varrho$, and for $|\lambda|=\varrho<1$,
$\varphi(\lambda)=\lim_{r\uparrow1}\varphi(r\lambda)$ lies in
$\overline{\mathcal B_\varrho}$ but not in $\mathcal B_\varrho$, by (d).

\emph{Step 2: a supporting inequality.} Let $0<\varrho<1$, $|\lambda|=\varrho$ and
$y=\varphi(\lambda)$. For $y'\in\mathcal B_\varrho$ the points
$y+\tau(y'-y)$, $0<\tau\le1$, lie in $\mathcal B_\varrho$ by (a), (b) and (e),
where $|G|^2<\varrho^2=|G(y)|^2$ by (d). Differentiating at $\tau=0$ shows that
the real affine function
$\Lambda(y')=\operatorname{Re}(\bar\lambda\,h(\lambda)\cdot(y'-y))$ is
$\le0$ on $\mathcal B_\varrho$. It is nonconstant, so it has no maximum on the
open set $\mathcal B_\varrho$, and therefore, for $z\in\mathcal B_\varrho$,
\begin{equation}\label{eq:affine-support}
 \operatorname{Re}\bigl(\bar\lambda\,g_z(\lambda)\bigr)<0 .
\end{equation}

\emph{Step 3: exactly one zero.} Let $z\in\mathcal B_\varrho$. Since
$\bar\lambda=\varrho^2/\lambda$ on $|\lambda|=\varrho$,
\eqref{eq:affine-support} says that $g_z(\lambda)/\lambda$ has negative
real part on this circle, so its winding number about $0$ vanishes. By
the argument principle $g_z=\lambda\cdot(g_z/\lambda)$ has exactly one
zero in $\{|\lambda|<\varrho\}$, and it is simple. Given $z\in D$ and zeros
$\lambda_1,\lambda_2$ of $g_z$, choose by (c) a radius $\varrho<1$ with
$z\in\mathcal B_\varrho$ and $|\lambda_1|,|\lambda_2|<\varrho$; then
$\lambda_1=\lambda_2$, and by (c) a zero exists. This defines $u(z)$.

\emph{Step 4: holomorphy.} For $z\in\mathcal B_\varrho$,
$u(z)=\frac1{2\pi i}\oint_{|\lambda|=\varrho}\lambda\,\partial_\lambda
g_z(\lambda)\,g_z(\lambda)^{-1}\,d\lambda$, and by Step~3 the integrand
is continuous on $\mathcal B_\varrho\times\{|\lambda|=\varrho\}$ and holomorphic in
$z$. Hence $u$ is holomorphic on each $\mathcal B_\varrho$, and so on $D$.

\emph{Step 5: left inverse and fibers.} Since $g_{\varphi(\lambda)}(\lambda)=0$,
$u\circ\varphi=\Id$. By definition $u(z)=\lambda$ if and only if
$g_z(\lambda)=0$, which is~\eqref{eq:affine-fibers}. Thus $u$ is constant
on the relatively open subset $u^{-1}(\lambda)$ of the hyperplane
$\varphi(\lambda)+\ker h(\lambda)$, so $du_{\varphi(\lambda)}$ vanishes
on $\ker h(\lambda)$. Both $du_{\varphi(\lambda)}$ and $h(\lambda)$ take
the value $1$ on $\varphi'(\lambda)$, hence they are equal.
\end{proof}

Lemma~\ref{lem:nontangential-fibers} uses only that $u$ is a
holomorphic left inverse of $\varphi$ satisfying~\eqref{eq:affine-fibers}
with $h$ holomorphic and $h\cdot\varphi'=1$.

The distinction between coercive and non-coercive models is already
visible in one tangential variable. A slice of
$\{\operatorname{Re}z_0+|w|^2<0\}$ at $z_0=-1$ is bounded, whereas the
corresponding slice for $(\operatorname{Im}w)^2$ is an unbounded strip.
The argument below first controls the direction of the scaled fibers,
then uses coercivity and convexity to control their entire extent.
Only after that step is the Jacobian used to estimate their measure.

\begin{lemma}[Nontangential fibers]\label{lem:nontangential-fibers}
Let $D\Subset\C^n$ be a smooth convex domain of finite type,
$d\mu\asymp\delta_D^\alpha dV$, $\alpha>-1$, $\varphi:\Delta\to D$ a
complex geodesic, and $u:D\to\Delta$ a left inverse with affine fibers,
\[
 u^{-1}(\lambda)=D\cap\{z:(z-\varphi(\lambda))\cdot h(\lambda)=0\},
 \qquad h(\lambda)=du_{\varphi(\lambda)},
\]
the dot being complex bilinear, so that $h\cdot\varphi'=1$. Put
$\nu=u_*\mu$. Let $\zeta_0\in\partial\Delta$ be a point at which
$q=\varphi^*(\zeta_0)$ is coercive, and use the
multitype coordinates~\eqref{eq:multitype-normal-form} centered at $q$.
For every $M\ge2$ there are $c,s_0>0$ such that, for
$\lambda\in\mathcal A_M(\zeta_0)$ with $s(\lambda)<s_0$,
\begin{enumerate}[(a)]
\item $u^{-1}(\lambda)\subset\{|z_0|\le c\,s(\lambda),\
|z_j|\le c\,s(\lambda)^{1/m_j}\ (1\le j\le n-1)\}$, and
\item $\nu\le c\,s^{\alpha+T(q)}dA$ on $\mathcal A_M(\zeta_0)\cap\{s<s_0\}$.
\end{enumerate}
\end{lemma}

\begin{proof}
Rotating the parameter we take $\zeta_0=1$. Put
$\psi=\varphi\circ\mathscr C^{-1}:\mathbb L\to D$ and
$\tilde u=\mathscr C\circ u:D\to\mathbb L$, so that
$\tilde u\circ\psi=\Id$ and the fibers of $\tilde u$ are those of $u$.
As in the proof of Lemma~\ref{lem:endpoint-scaling}, $0$ is a regular
contact point of $\psi$ with endpoint $q$, and $\langle\psi',e_0\rangle$ has a nontangential limit $\kappa>0$ at $0$
\cite[Remark~4.7, Proposition~4.9]{ArosioFiacchi2025}. With
$\mathcal S_t$ as in the proof of Lemma~\ref{lem:endpoint-scaling} put
\[
 \Phi_t(\sigma)=\mathcal S_t\psi(t\sigma),\qquad
 U_t(y)=t^{-1}\tilde u(\mathcal S_t^{-1}y)\quad(y\in \mathcal S_tD),
\]
so that $U_t:\mathcal S_tD\to\mathbb L$ and $U_t\circ\Phi_t=\Id_{\mathbb L}$.
By~\cite[Proposition~6.6]{ArosioFiacchi2025}, applied along every
sequence $t_k\downarrow0$, $\Phi_t(\sigma)\to(\kappa\sigma,0)$ uniformly on
compact subsets of $\mathbb L$ as $t\downarrow0$. Fix $M$ and the
compact sets
\begin{align*}
 Q_M&=\{\sigma:\operatorname{Re}\sigma<0,\ \tfrac12\le|\sigma|\le1,\
       |\sigma|\le M|\operatorname{Re}\sigma|\},\\
 Q'_M&=\{\sigma:\operatorname{Re}\sigma<0,\ \tfrac14\le|\sigma|\le2,\
       |\sigma|\le2M|\operatorname{Re}\sigma|\};
\end{align*}
on them $\operatorname{Re}\sigma\le-1/(8M)$, so
$Q_M\subset\operatorname{int}Q'_M$ and $Q'_M\subset\mathbb L$.

\emph{Step 1: the scaled domains.} Let $\rho_t(y)=r(\mathcal S_t^{-1}y)/t$, with
$r$ as in~\eqref{eq:multitype-normal-form}. Weighted homogeneity gives
$\rho_t(y)=\operatorname{Re}y_0+H(y')+R(\mathcal S_t^{-1}y)/t$, and the
remainder tends to zero uniformly on bounded sets. Since
$\mathcal S_t^{-1}$ shrinks bounded sets to $q$, for each $b_0$ we have, for
small $t$, $\mathcal S_tD\cap \mathbb B_{b_0}=\{\rho_t<0\}\cap \mathbb B_{b_0}$. Consequently
each compact subset of $\mathcal M_H$ lies in $\mathcal S_tD$ for small $t$
(cf.~\cite[Theorem~6.2]{ArosioFiacchi2025}), and
\begin{equation}\label{eq:scaled-domain-upper}
 \mathcal S_tD\cap \mathbb B_{b_0}\subset\{\operatorname{Re}y_0+H(y')<e_t(b_0)\},
 \qquad e_t(b_0)\to0 .
\end{equation}

\emph{Step 2: the scaled differentials.} Put
$h_t(\sigma)=d(U_t)_{\Phi_t(\sigma)}$, a holomorphic function of
$\sigma$ with values in $(\C^n)^*$, so that $h_t\cdot\Phi_t'=1$. We
claim $h_t\to \kappa^{-1}dz_0$ uniformly on $Q'_M$. Otherwise there are
$c_1>0$, $t_k\downarrow0$ and $\sigma_k\to\sigma$ in $Q'_M$ with
$|h_{t_k}(\sigma_k)-\kappa^{-1}dz_0|\ge c_1$. The maps
$\mathscr C^{-1}\circ U_{t_k}$ are bounded by one and, by Step~1,
defined on each compact subset of $\mathcal M_H$ for large $k$. A subsequence
converges locally uniformly on $\mathcal M_H$ to a holomorphic $V$ with
$|V|\le1$. Since $\Phi_{t_k}\to(\kappa\,\cdot\,,0)$ locally uniformly,
$V(\kappa\sigma',0)=\mathscr C^{-1}(\sigma')\in\Delta$ for all
$\sigma'\in\mathbb L$, so $|V|<1$ by the maximum principle, and
$U=\mathscr C\circ V:\mathcal M_H\to\mathbb L$ satisfies $U(\kappa\sigma',0)=\sigma'$.
Derivatives converge locally uniformly, so
$h_{t_k}(\sigma_k)\to dU_{(\kappa\sigma,0)}=\kappa^{-1}dz_0$ by
Lemma~\ref{lem:model-differential}, a contradiction.

Write $h_t=(h_{t,0},\ldots,h_{t,n-1})$. For small $t$, $h_{t,0}\ne0$
on $Q'_M$, and $\theta_j=h_{t,j}/h_{t,0}\to0$ uniformly on $Q'_M$, hence,
by Cauchy's estimates, $\theta_j'\to0$ uniformly on $Q_M$. The vectors
$\mathbf k_j=e_j-\theta_je_0$, $1\le j\le n-1$, span $\ker h_t$.

\emph{Step 3: bounded fibers.} By the affine-fiber hypothesis, for
$\sigma\in Q'_M$,
\[
 U_t^{-1}(\sigma)=\mathcal S_tD\cap\Pi_t(\sigma),\qquad
 \Pi_t(\sigma)=\Bigl\{\Phi_t(\sigma)+\sum_{j\ge1}w_j\mathbf k_j(\sigma):
 w\in\C^{n-1}\Bigr\}.
\]
Let $r_M=\max_{Q_M}|\sigma|$ and let $b_2$ exceed
$\kappa r_M+2+\sum_j((\kappa r_M+2)/c_H)^{1/m_j}$, with $c_H$
from~\eqref{eq:H-coercive}; put $b_0=2b_2$. Let $\sigma\in Q_M$ and
$y\in\Pi_t(\sigma)\cap \mathcal S_tD\cap \mathbb B_{b_0}$. On $\Pi_t(\sigma)$,
$y_0=\Phi_{t,0}-\sum_j\theta_j(y_j-\Phi_{t,j})$, so
$\operatorname{Re}y_0\ge\operatorname{Re}(\kappa\sigma)-o(1)$
uniformly. By~\eqref{eq:scaled-domain-upper}
and~\eqref{eq:H-coercive}, for small $t$,
$c_H\sum_j|y_j|^{m_j}\le H(y')<\kappa r_M+1$, and
$|y_0|\le \kappa r_M+1$. Thus $|y|<b_2$. The set $\mathcal S_tD\cap\Pi_t(\sigma)$ is
convex and contains $\Phi_t(\sigma)\in \mathbb B_{b_2}$; if it had a point outside
$\mathbb B_{b_0}$, the segment to that point would meet the sphere of radius
$(b_0+b_2)/2$ inside $\mathcal S_tD\cap\Pi_t(\sigma)\cap \mathbb B_{b_0}$, which is
impossible. Hence $U_t^{-1}(\sigma)\subset \mathbb B_{b_2}$ for $\sigma\in Q_M$
and $t<t_0$.

\emph{Step 4: the density.} The map
$\Psi_t(\sigma,w)=\Phi_t(\sigma)+\sum_jw_j\mathbf k_j(\sigma)$ is
holomorphic and injective on
$\{(\sigma,w)\in Q'_M\times\C^{n-1}:\Psi_t\in \mathcal S_tD\}$, since
$U_t\circ\Psi_t=\sigma$; its image is $U_t^{-1}(Q'_M)$. Adding
$\theta_j$ times row $j$ to row $0$ and using $h_t\cdot\Phi_t'=1$,
\[
 \det{}_\C d\Psi_t
 =\Phi'_{t,0}+\sum_j\theta_j\Phi'_{t,j}-\sum_jw_j\theta_j'
 =\frac1{h_{t,0}}-\sum_jw_j\theta_j' .
\]
For $\sigma\in Q_M$ the section
$\mathcal F_t(\sigma)=\{w:\Psi_t(\sigma,w)\in \mathcal S_tD\}$ is convex. For
$w\in \mathcal F_t(\sigma)$, both $\Psi_t(\sigma,w)$ and $\Phi_t(\sigma)$ lie in
$\mathbb B_{b_2}$ by Step~3, and the $e_j$-component ($j\ge1$) of
$\Psi_t(\sigma,w)-\Phi_t(\sigma)=\sum_jw_j\mathbf k_j$ is exactly $w_j$,
so $|w|\le|\sum_jw_j\mathbf k_j|<2b_2$. Thus
$\mathcal F_t(\sigma)\subset\{|w|<2b_2\}$, and there
$|\det d\Psi_t|\le2\kappa+1$ for small $t$. In fact
$\det d\Psi_t\to\kappa$ uniformly on these full sections: the
$w_j$ are uniformly bounded and $\theta_j'\to0$ on $Q_M$. No
boundary derivative of the original dual data is assumed here.

Let $\sigma\in Q_M$, $0<t<t_0$, and $\lambda=\mathscr C^{-1}(t\sigma)$.
Since $s(\lambda)=-4t\operatorname{Re}\sigma/|1-t\sigma|^2$, we have
$s(\lambda)\asymp t$ with constants depending on $Q_M$. The function
$d(w)=\delta_D(\mathcal S_t^{-1}\Psi_t(\sigma,w))$ is concave on $\mathcal F_t(\sigma)$,
since boundary distance in a convex domain is the infimum of affine
supporting distances. By the depth lemma,
$d(0)=\delta_D(\varphi(\lambda))\ge cs(\lambda)$, and by
Schwarz--Pick, as in~\eqref{eq:global-moment-depth},
$\delta_D\le Cs(\lambda)$ on $u^{-1}(\lambda)$. Thus
$d_*=\max d\asymp t$. Contracting $\mathcal F_t(\sigma)$ about a maximum point by
the factor $1-\tau/d_*$ gives $|\{d<\tau\}|\le(2n-2)\tau|\mathcal F_t(\sigma)|/d_*$,
and for $-1<\alpha<0$ integration of this distribution bound gives
$\int_{\mathcal F_t(\sigma)}d^\alpha\le(1+(2n-2)(-\alpha)/(1+\alpha))d_*^\alpha
|\mathcal F_t(\sigma)|$; for $\alpha\ge0$ use $d\le d_*$. In both cases
$\int_{\mathcal F_t(\sigma)}d^\alpha\,dV(w)\le Ct^\alpha$.

The real Jacobian of $\mathcal S_t^{-1}$ is $t^{2+T(q)}$. For Borel $S\subset Q_M$,
change of variables by $\mathcal S_t^{-1}$ and $\Psi_t$ gives
\[
 \tilde u_*\mu(tS)\le C\,t^{2+T(q)}\int_S\int_{\mathcal F_t(\sigma)}
 d^\alpha|\det d\Psi_t|^2\,dV(w)\,dA(\sigma)
 \le C\,t^{2+\alpha+T(q)}|S| .
\]
Since $|tS|=t^2|S|$ and $t/2\le|\sigma'|\le t$ on $tQ_M$,
$\tilde u_*\mu\le C|\sigma'|^{\alpha+T(q)}dA(\sigma')$ on $tQ_M$ for each
$t<t_0$, with $C$ independent of $t$. Let
$\Sigma=\{\sigma':\operatorname{Re}\sigma'<0,\ 0<|\sigma'|<t_0/2,\
|\sigma'|\le M|\operatorname{Re}\sigma'|\}$. For $\sigma'\in\Sigma$ and every $t\in[|\sigma'|,2|\sigma'|]$ we have
$t<t_0$ and $\sigma'/t\in Q_M$; choosing $t$ rational, the countably
many compact sets $tQ_M$, $t\in\mathbb Q\cap(0,t_0)$, cover $\Sigma$. By the
remark on countable unions,
$\tilde u_*\mu\le C|\sigma'|^{\alpha+T(q)}dA(\sigma')$ on $\Sigma$.
Now $\lambda\in\mathcal A_M(1)$ implies
$|\mathscr C(\lambda)|\le M|\operatorname{Re}\mathscr C(\lambda)|$,
because $\operatorname{Re}\mathscr C(\lambda)=-s(\lambda)/|1+\lambda|^2$
and $|\mathscr C(\lambda)|=|\lambda-1|/|\lambda+1|$. Hence
$\mathscr C(\mathcal A_M(1)\cap\{s<s_0\})\subset \Sigma$ for small $s_0$. On
this set $|\mathscr C'|\le2$ and $|\mathscr C(\lambda)|\asymp_M s(\lambda)$.
Pulling back by $\mathscr C$ proves (b).

For (a), let $\lambda\in\mathcal A_M(1)$ with $s(\lambda)<s_0$, and write
$\mathscr C(\lambda)=t\sigma$ with $t=|\mathscr C(\lambda)|$, so
$\sigma\in Q_M$. By Step~3,
$\mathcal S_tu^{-1}(\lambda)=U_t^{-1}(\sigma)\subset \mathbb B_{b_2}$, that is,
$|z_0|\le b_2t$ and $|z_j|\le b_2t^{1/m_j}$ on $u^{-1}(\lambda)$, and
$t\asymp_Ms(\lambda)$.
\end{proof}

Steps~1 and~2 of this proof do not use coercivity: at every
$\zeta_0\in\partial\Delta$, $h_t\to\kappa^{-1}dz_0$, so after scaling
the fibers of $u$ become parallel to the complex tangent space at
$\varphi^*(\zeta_0)$. Coercivity enters only in Step~3. At a
non-coercive point the slices of the model domain by these hyperplanes
are unbounded; this is relevant to the conjecture in
Section~\ref{sec:questions}.

\begin{proof}[Proof of Theorem~\ref{thm:regular-reverse}]
By Section~\ref{sec:geodesic-restriction}, $\varphi$ has a
holomorphic left inverse $G:D\to\Delta$. By
Lemma~\ref{lem:affine-inverse}, for every $z\in D$ there is exactly one
$\lambda\in\Delta$ with $(z-\varphi(\lambda))\cdot G'(\varphi(\lambda))=0$,
and the map $u:z\mapsto\lambda$ is a holomorphic left inverse of
$\varphi$ with affine fibers
$u^{-1}(\lambda)=D\cap\{z:(z-\varphi(\lambda))\cdot h(\lambda)=0\}$,
where $h(\lambda)=du_{\varphi(\lambda)}=G'(\varphi(\lambda))$. Then $u(\varphi(0))=0$,
$du(\varphi'(0))=1$, and $\gamma_D(\varphi(0);\varphi'(0))=1$, since
the Carath\'eodory metric is at least $|du(\varphi'(0))|$ and at most
the Kobayashi metric.

(i) The multitype takes values in $\{2,\ldots,L\}^{n-1}$, $L$ the
maximal D'Angelo type, so $T$ takes values in a finite set $\mathcal T$, and $T\circ\varphi^*$ is measurable (proof of
Theorem~\ref{thm:general-main}(ii)). Put $u_a=\mob_a\circ u$, so that
$u_a(\varphi(a))=0$ and
$du_a(\varphi'(a))=(1-|a|^2)^{-1}=\gamma_D(\varphi(a);\varphi'(a))$.
Let $\nu=u_*\mu$ and $\nu_a=u_{a*}\mu=(\mob_a)_*\nu$, and let
$\zeta\in\mathsf C$. Since $\mob_{-a}$ is a diffeomorphism of a
neighborhood of $\overline\Delta$ preserving $\partial\Delta$, for every
$M'\ge2$ there is $M=M(M',a)\ge2$ with
$\mob_{-a}(\mathcal A_{M'}(\mob_a\zeta))\subset\mathcal A_M(\zeta)$ near
$\zeta$. With $c,s_0$ from Lemma~\ref{lem:nontangential-fibers}(b) for
this $M$, so that $\nu\le c\,s^{\alpha+T(\varphi^*(\zeta))}dA$ on
$\mathcal A_M(\zeta)\cap\{s<s_0\}$, choose $s_1>0$ with
$\mob_{-a}(\mathcal A_{M'}(\mob_a\zeta)\cap\{s<s_1\})\subset
\mathcal A_M(\zeta)\cap\{s<s_0\}$; moreover $s\circ\mob_{-a}\asymp_as$ and
$|\mob_{-a}'|\asymp_a1$. Therefore
$\nu_a\le c'\,s^{\alpha+T(\varphi^*(\zeta))}dA$ on
$\mathcal A_{M'}(\mob_a\zeta)\cap\{s<s_1\}$. Apply
Lemma~\ref{lem:stolz-transfer} with $x=\varphi(a)$, $v=\varphi'(a)$, the
function $u_a$, the Borel set $\mob_a(\mathsf C)$ and
$\beta=T\circ\varphi^*\circ\mob_{-a}$. The boundary map of $\mob_{-a}$
carries $\omega_0$ to $\omega_a$, so
$\int_{\mob_a(\mathsf C)}\beta\,d\omega_0=\int_{\mathsf C}T(\varphi^*)\,d\omega_a$.
Homogeneity removes $\eta$.

(ii) By hypothesis, the set of $\zeta$ at which $\varphi^*(\zeta)$ is
not coercive is contained in a Borel null set $Z$. Take
$\mathsf C=\partial\Delta\setminus Z$ in (i), and combine it with
Theorem~\ref{thm:general-main}(ii).

(iii) By Lemma~\ref{lem:holder} the radial limits $\varphi^*(\zeta)$
exist, $\partial D$ is smooth, and by Lemma~\ref{lem:coercive-face}
every coercive point has trivial complex face. Corollary~\ref{cor:uniqueness}
applies.
\end{proof}

\begin{proposition}[Three contact orders along one circle]
\label{prop:three-types}
With $\chi_2$ as in~\eqref{eq:chi-nu}, put
\[
 a_+(x)=\chi_2(x-1/2),\qquad a_-(x)=\chi_2(x+1/2),\qquad
 Q(w)=|w|^4+\tfrac1{12}\operatorname{Re}w^4.
\]
For $0<\eps\le1/100$, put
\begin{equation}\label{eq:three-types-domain}
 \rho_\eps(z,w)=|z|^2+\eps a_+(\operatorname{Re}z)|w|^2
       +\eps a_-(\operatorname{Re}z)Q(w)+|w|^6-1,
 \qquad D_\eps=\{\rho_\eps<0\}.
\end{equation}
Then $D_\eps\Subset\C^2$ is smooth, convex and of finite type.
The axis disc $\varphi(\zeta)=(\zeta,0)$ is a complex geodesic, and
$\varphi^*(\zeta)=(\zeta,0)$ is coercive for every
$\zeta\in\partial\Delta$, including the four transition points
$\operatorname{Re}\zeta=\pm1/2$. For any probability measure
$d\mu\asymp\delta_{D_\eps}^\alpha dV$, $\alpha>-1$, its axis
coefficient is
\begin{equation}\label{eq:three-types-average}
 N(a)=1+\alpha+\omega_a(E_2)+\tfrac12\omega_a(E_4)
                      +\tfrac13\omega_a(E_6),\qquad a\in\Delta,
\end{equation}
where $E_2=\{\operatorname{Re}\zeta>1/2\}$,
$E_4=\{-1/2<\operatorname{Re}\zeta<1/2\}$ and
$E_6=\{\operatorname{Re}\zeta<-1/2\}$ are subsets of $\partial\Delta$.
In particular $N(0)=\alpha+29/18$, and $N$ is nonconstant.
\end{proposition}
\begin{proof}
By Section~\ref{sec:flat-convex}, $\chi_2$ is smooth, flat at zero,
positive on $(0,\infty)$, and satisfies $0\le\chi_2''\le1$ and
$(\chi_2')^2\le2\chi_2$. Write $r=|w|$. Real derivatives give
\[
 \tfrac{11}{12}r^4\le Q(w)\le\tfrac{13}{12}r^4,
 \qquad \nabla_w^2Q\ge3r^2 I,\qquad |\nabla Q|\le5r^3.
\]
Here $\nabla^2|w|^4\ge4r^2I$ and
$\|\nabla^2\operatorname{Re}w^4\|=12r^2$.
Since $\rho_\eps\ge|z|^2+|w|^6-1$, its nonpositive sublevel
lies in the convex cylinder $|z|\le1$, $r\le1$. On this cylinder
the real $z$-block of its Hessian is at least $2I$, and the $w$-block
is at least $\Lambda I$, where
$\Lambda=2\eps a_++3\eps a_- r^2+6r^4$.
The mixed block has only one nonzero row,
$\mathbf m=\eps(a_+'\nabla r^2+a_-'\nabla Q)$, and
\[
 |\mathbf m|^2\le16\eps^2 a_+ r^2+100\eps^2 a_- r^6
          \le\tfrac{100}{3}\eps\Lambda\le\tfrac13\Lambda.
\]
Consequently, for real vectors $\xi,\upsilon\in\R^2$,
\begin{equation}\label{eq:three-types-hessian}
 \nabla^2\rho_\eps[(\xi,\upsilon)]
 \ge |\xi|^2+\tfrac23\Lambda|\upsilon|^2.
\end{equation}
The Hessian bound proves convexity of the sublevel. Since
$\rho_\eps(0)=-1$, convexity gives $d\rho_\eps(Z)[Z]\ge1$
on its zero level, proving smoothness. Strong convexity holds where
$w\ne0$ or $a_+>0$.

Let $\mathfrak c=\{(\zeta,0):|\zeta|=1\}$. At $q_\zeta=(\zeta,0)$
the complex tangent line is vertical, and for every unit complex
number $e$ the exact restriction is
\begin{equation}\label{eq:three-types-contact}
 \rho_\eps(q_\zeta+(0,we))
  =\eps a_+(\operatorname{Re}\zeta)|w|^2
     +\eps a_-(\operatorname{Re}\zeta)Q(we)+|w|^6.
\end{equation}
Its complex-line order is $2$ on $E_2$, $4$ on $E_4$ and at
$\operatorname{Re}\zeta=1/2$, where $a_+=0<a_-$, and $6$ on $E_6$ and at
$\operatorname{Re}\zeta=-1/2$, where $a_+=a_-=0$. The convex line-type
criterion~\cite{McNeal1992Convex} gives these D'Angelo types and
finite type everywhere.

\emph{Normal form.} In the unitary coordinates $z=\zeta(1+z_0)$, $w$
centered at $q_\zeta$ one has $|z|^2-1=2\operatorname{Re}z_0+|z_0|^2$
and $|\operatorname{Re}z-\operatorname{Re}\zeta|\le|z_0|$, so
\[
 \tfrac12\rho_\eps=\operatorname{Re}z_0+H_\zeta(w)+R_\zeta(z_0,w),
\]
with $m_1=2$ and $H_\zeta=\frac\eps2a_+(\operatorname{Re}\zeta)|w|^2$ on
$E_2$; $m_1=4$ and $H_\zeta=\frac\eps2a_-(\operatorname{Re}\zeta)Q(w)\ge
\frac{11\eps}{24}a_-(\operatorname{Re}\zeta)|w|^4$ on
$E_4\cup\{\operatorname{Re}\zeta=1/2\}$; and $m_1=6$ and
$H_\zeta=\frac12|w|^6$ on $E_6\cup\{\operatorname{Re}\zeta=-1/2\}$. The
remainder $R_\zeta$ consists of $\frac12|z_0|^2$, of the terms
of~\eqref{eq:three-types-contact} of order above $m_1$ in $w$, and of
$\frac\eps2(a_\pm(\operatorname{Re}z)-a_\pm(\operatorname{Re}\zeta))$
times $|w|^2$, respectively $Q(w)$. A difference with
$a_\pm(\operatorname{Re}\zeta)>0$ is $O(|z_0|)$, and its product with a
power of $|w|$ is $o(|z_0|)$. If $a_\pm(\operatorname{Re}\zeta)=0$, then
$a_\pm$ vanishes to infinite order at $\operatorname{Re}\zeta$ (as a
function of $x=\operatorname{Re}z$, $a_\pm$ vanishes for $x<\pm1/2$
and is flat at $x=\pm1/2$), so
$a_\pm(\operatorname{Re}z)=O(|z_0|^2)$ and the term is $o(|z_0|)$.
Hence $R_\zeta=o(|z_0|+|w|^{m_1})$, which
is~\eqref{eq:multitype-normal-form} with $H=H_\zeta$. Each $H_\zeta$ is
positive for $w\ne0$, so every point of $\mathfrak c$ is coercive.
Since $D_\eps\subset\Delta\times\C$ and contains $\Delta\times\{0\}$,
the coordinate $u(z,w)=z$ satisfies $u\circ\varphi=\Id$; contraction of
the Kobayashi metric proves that $\varphi$ is a geodesic.
Theorem~\ref{thm:regular-reverse}(ii) proves
\eqref{eq:three-types-average}, and Corollary~\ref{cor:uniqueness}
shows that the axis disc is, up to disc automorphisms, the unique
complex geodesic through each of its jets. The three sets have arclength $2\pi/3$,
so $N(0)=\alpha+29/18$. As real $a\to1$ or $a\to-1$, the respective
limits are $2+\alpha$ and $4/3+\alpha$, proving nonconstancy.
\end{proof}

\begin{corollary}[An unbounded remainder in the three-type example]
\label{cor:three-types-remainder}
For each $a\in\Delta$ and $v\ne0$ in
Proposition~\ref{prop:three-types},
\[
 p\log\frac{B_{p,\mu}((a,0);(v,0))}{\gamma_D((a,0);(v,0))}
                    -N(a)\log p\longrightarrow-\infty.
\]
\end{corollary}
\begin{proof}
Write $\hat a_\pm(z)=\eps a_\pm(\operatorname{Re}z)$, $s=1-|z|^2$, and
$d_s(z)=\hat a_+(z)+\sqrt{\hat a_-(z)s}+s^{2/3}$. The fiber contains
$|w|<c\sqrt{s/d_s}$; there each of the three terms in
\eqref{eq:three-types-domain} is at most a fixed small multiple of
$s$, so $\delta_D(z,w)\asymp s$. Submean on this disc gives, for
$g(z)=f(z,0)$ and either sign of $\alpha$,
\begin{equation}\label{eq:three-type-restriction}
 \|f\|_p^p\ge c\int_{|z|>r_0}|g(z)|^p
                         \frac{s^{\alpha+1}}{d_s(z)}\,dA(z).
\end{equation}
On $s\in[1/p,2/p]$ one has, uniformly in $\zeta\in\partial\Delta$,
\[
 d_s(r\zeta)\asymp d_p(\zeta):=
 \hat a_+(\zeta)+\sqrt{\hat a_-(\zeta)}\,p^{-1/2}+p^{-2/3}.
\]
Indeed $a_+$ and $\sqrt{a_-}$ are Lipschitz, the latter by
$(a_-')^2\le2a_-$; their radial errors are absorbed by the floor
$p^{-2/3}$. For $g(0)=0$, $g'(0)=1$, weighted Jensen on each circle gives
$\int|g(r\zeta)|^p/d_s(r\zeta)\,dm
\ge r^p\exp[-\int\log d_s(r\zeta)\,dm]$, where $dm=|d\zeta|/(2\pi)$.
Integrating~\eqref{eq:three-type-restriction} over these depths, and
using $E_p\ge1$, gives
\begin{equation}\label{eq:three-type-remainder-bound}
 p\log(B_p/\gamma_D)\le(2+\alpha)\log p
                       +\int\log d_p\,d\omega_a+O_a(1).
\end{equation}
For general $a$, precompose $g$ by a disc automorphism taking $0$ to
$a$; its depth and area factors are comparable by constants depending
on $a$, and its angular measure becomes $\omega_a$. The same radial
comparison proves~\eqref{eq:three-type-remainder-bound}.

Set $\vartheta=0,1/2,2/3$ on $E_2,E_4,E_6$, respectively.
Then $N(a)=2+\alpha-\int\vartheta\,d\omega_a$. The functions
$\log(p^\vartheta d_p)$ are bounded above uniformly in $p$ and tend to
$\log\hat a_+$ on $E_2$, $\frac12\log\hat a_-$ on $E_4$, and
$0$ on $E_6$. By~\eqref{eq:chi-bounds},
$\log\chi_2(t)=-t^{-2}+O(1+|\log t|)$ for small $t>0$, so the
positive-side logarithms are nonintegrable at the four transverse
crossings $\operatorname{Re}\zeta=\pm1/2$. Reverse Fatou therefore
gives $\int\log(p^\vartheta d_p)\,d\omega_a\to-\infty$;
\eqref{eq:three-type-remainder-bound} proves the assertion.
\end{proof}

\begin{example}[Every coefficient in an interval]
\label{ex:continuum}
Let $k=n-1\ge1$, $m\ge2$, $-1<c<1$, $a_c(z)=\chi(\operatorname{Re}z-c)$
with $\chi=\chi_{1/2}$, and
\[
 D_c=\{(z,w)\in\C\times\C^k:\rho_c(z,w)<0\},\qquad
 \rho_c=|z|^2+\tfrac18a_c(z)|w|^2+|w|^{2m}-1.
\]
Then $D_c$ is a smooth bounded convex domain of finite type, the axis
disc $\varphi(\lambda)=(\lambda,0)$ is a complex geodesic, and every
point of $\mathfrak c$ is coercive, with tangential multitype
$(2m,\ldots,2m)$ on $Z\times\{0\}$, $Z=\{\zeta\in\partial\Delta:
\operatorname{Re}\zeta\le c\}$, and $(2,\ldots,2)$ elsewhere. By
Theorem~\ref{thm:regular-reverse}(ii) the coefficient along the axis is
\[
 N(x)=1+\alpha+k\Bigl(1-d_c(x)+\frac{d_c(x)}m\Bigr),\qquad
 d_c(x)=\omega_x(Z),\qquad d_c(0)=1-\frac{\arccos c}\pi .
\]
As $c$ ranges over $(-1,1)$ and $m\ge2$, $N(0)$ takes every value in
$(1+\alpha,n+\alpha)$, and $n+\alpha$ is attained by the ball. For
fixed $c$, $N$ is a nonconstant harmonic function on the axis disc,
tending to $n+\alpha$ and $1+\alpha+k/m$ at $x=1$ and $x=-1$. For $k=1$,
$m=2$, $c=0$ and $\alpha=0$, $N(0)=7/4$ lies strictly between the
values $3/2$ and $2$ given by the two multitypes, and on the real
diameter $N(x)=7/4+\arctan(x)/\pi$.
\end{example}

\begin{proof}
The function $a=a_c$ is smooth, convex and nonnegative, with
$|\nabla a|^2=\chi'(\operatorname{Re}z-c)^2\le2a$
(Section~\ref{sec:flat-convex}). On the convex cylinder
$\{|z|\le1,\ |w|\le1\}$ and for real vectors $(\xi,\upsilon)$,
$\nabla^2(a|w|^2)[(\xi,\upsilon)]\ge2a|\upsilon|^2-4|\nabla a||\xi||w||\upsilon|$
and $\frac12|\nabla a||\xi||w||\upsilon|\le|\xi|^2+\frac1{16}|\nabla a|^2
|\upsilon|^2\le|\xi|^2+\frac18a|\upsilon|^2$, so
\[
 \nabla^2\rho_c[(\xi,\upsilon)]\ge|\xi|^2+\tfrac18a|\upsilon|^2
 +2m|w|^{2m-2}|\upsilon|^2 .
\]
Since $\rho_c\ge|z|^2+|w|^{2m}-1$, the set $D_c$ lies in this cylinder, so it
is convex and bounded; $\rho_c(0)=-1$ and convexity give $d\rho_c(Z)[Z]\ge1$
where $\rho_c(Z)=0$, so $\partial D_c$ is smooth. Where $w\ne0$ the Hessian
is positive definite and $\partial D_c$ is strongly convex. The
boundary points with $w=0$ are those of $\mathfrak c$; there the complex
tangent space is $\{0\}\times\C^k$, and
$\rho_c(\zeta,t\upsilon)=\frac18a(\zeta)|t\upsilon|^2+|t\upsilon|^{2m}$,
so every complex tangent line has type $2$ if $a(\zeta)>0$ and $2m$ if
$a(\zeta)=0$, that is, if $\zeta\in Z$. By~\cite{McNeal1992Convex}
the domain has finite type, and Lemma~\ref{lem:linear-multitype} gives
the multitype. In the coordinates $z=\zeta(1+z_0)$, $w$,
\[
 \tfrac12\rho_c=\operatorname{Re}z_0+\tfrac1{16}a(\zeta)|w|^2
 +\tfrac12|w|^{2m}+\tfrac12|z_0|^2
 +\tfrac1{16}\bigl(a(\zeta(1+z_0))-a(\zeta)\bigr)|w|^2 .
\]
If $a(\zeta)>0$ the last two terms and $\frac12|w|^{2m}$ are
$o(|z_0|+|w|^2)$; if $a(\zeta)=0$, then $\nabla a(\zeta)=0$, so
$a(\zeta(1+z_0))=O(|z_0|^2)$ and the last two terms are $o(|z_0|)$.
Hence~\eqref{eq:multitype-normal-form} holds with
$H=\frac1{16}a(\zeta)|w|^2$, respectively $H=\frac12|w|^{2m}$; in either
case every orthonormal frame of $\{0\}\times\C^k$ is adapted, $H$ is
the weight-one part of the defining function on the complex tangent
space, and $H>0$ for $w\ne0$. Thus every point of $\mathfrak c$ is coercive,
with $T=k$ off $Z$ and $T=k/m$ on $Z$. Since
$\Delta\times\{0\}\subset D_c\subset\Delta\times\C^k$, the projection
$z$ is a left inverse of $\varphi$, which is therefore a complex
geodesic. Theorem~\ref{thm:regular-reverse}(ii) gives the formula for
$N$; $d_c(0)$ is the normalized length of $Z$. Since $1\notin Z$ and
$-1$ is an interior point of $Z$, $\omega_x(Z)\to0$ and $\to1$ as real
$x\to1$ and $x\to-1$. For fixed $m$, $d_c(0)$ increases continuously
from $0$ to $1$ as $c$ increases, so $N(0)$ fills
$(1+\alpha+k/m,1+\alpha+k)$; let $m\to\infty$. The ball has $T\equiv
n-1$.
\end{proof}

\begin{lemma}[Coercive points have trivial complex face]
\label{lem:coercive-face}
Let $D\subset\C^n$ be a convex domain and let $q\in\partial D$ be a
point near which $\partial D$ is smooth and of finite type. If $q$ is
coercive, then $\overline D\cap(q+T^{1,0}_q\partial D)=\{q\}$.
\end{lemma}
\begin{proof}
In the coordinates~\eqref{eq:multitype-normal-form} centered at $q$,
$T^{1,0}_q\partial D=\{z_0=0\}$ and $\overline D=\{r\le0\}$ near $q$. Let
$v=\sum_{j\ge1}v_je_j\ne0$ and $v'=(v_1,\ldots,v_{n-1})$.
By~\eqref{eq:H-coercive}, $H(tv')\ge c_HS(t)$ with
$S(t)=\sum_j|tv_j|^{m_j}>0$ for $t>0$, and $R(0,tv')=o(S(t))$ as
$t\downarrow0$, so $r(q+tv)\ge(c_H-o(1))S(t)>0$ for small $t>0$. If
$q+v\in\overline D$, convexity would give $q+tv\in\overline D$ for
$0<t\le1$, which is impossible.
\end{proof}

\begin{proof}[Proof of Corollary~\ref{cor:uniqueness}]
A complex geodesic is a Kobayashi extremal for the pair
$(\varphi(0),\varphi(1/2))$, since the Kobayashi distance does not
exceed the Lempert function. By the theorem of Royden and Wong in the
form of~\cite[Theorem~2.1(ii)]{KosinskiMcCarthy2019} (proofs
in~\cite[Theorem~2.6.17]{Abate1989}
and~\cite[Lemmata~8.2.2, 8.2.4, Remark~8.2.3]{JarnickiPflug1993}; see
also~\cite[Theorem~4]{JarnickiPflugZeinstra1993}) there is
$\mathfrak d\in H^1(\Delta,\C^n)$, $\mathfrak d\not\equiv0$, such that
for every $z\in D$,
$\operatorname{Re}[(z-\varphi^*(\zeta))\cdot\bar\zeta\mathfrak d^*(\zeta)]<0$
for almost every $\zeta\in\partial\Delta$. Apply this to a countable
dense subset of $D$. As $\mathfrak d^*\ne0$ almost
everywhere~\cite[Theorem~2.6.12]{Abate1989}, for $\zeta$ in a set
$\mathsf E_1$ of full measure $\varphi^*(\zeta)$ and
$\mathfrak d^*(\zeta)\ne0$ exist and, by continuity in $z$, the complex
linear form $\Lambda_\zeta=\bar\zeta\mathfrak d^*(\zeta)\cdot(\,\cdot\,)$
satisfies $\operatorname{Re}\Lambda_\zeta(z-\varphi^*(\zeta))\le0$ for
$z\in\overline D$.

The product $(g-\varphi)\cdot\mathfrak d$ of a bounded and an $H^1$ map
lies in $H^1$ and vanishes at $0$ with derivative
$(g'(0)-\varphi'(0))\cdot\mathfrak d(0)=0$, so
$\mathfrak g=\lambda^{-1}(g-\varphi)\cdot\mathfrak d\in H^1$ and
$\mathfrak g(0)=0$. For almost every $\zeta\in\mathsf E_1$,
$g^*(\zeta)\in\overline D$ exists and
$\operatorname{Re}\mathfrak g^*(\zeta)
=\operatorname{Re}\Lambda_\zeta(g^*(\zeta)-\varphi^*(\zeta))\le0$. As
$\mathfrak g$ is the Poisson integral of
$\mathfrak g^*$~\cite[Theorem~2.6.12]{Abate1989},
$\int_{\partial\Delta}\operatorname{Re}\mathfrak g^*=2\pi\operatorname{Re}\mathfrak g(0)=0$,
so $\operatorname{Re}\mathfrak g^*=0$ almost everywhere,
$\mathfrak g\equiv0$, and
$g^*(\zeta)-\varphi^*(\zeta)\in\ker\Lambda_\zeta$ for almost every
$\zeta$. Let $\zeta\in\mathsf E\cap\mathsf E_1$, $q=\varphi^*(\zeta)$,
and $r$ a $C^1$ local defining function at $q$. If $dr_q(v)<0$, then
$q+tv\in D$ for small $t>0$, so $\operatorname{Re}\Lambda_\zeta(v)\le0$.
Hence $\operatorname{Re}\Lambda_\zeta=c\,dr_q$ with $c>0$, as
$\Lambda_\zeta\ne0$, and since $\Lambda_\zeta$ is complex linear,
$\ker\Lambda_\zeta=\{v:dr_q(v)=dr_q(iv)=0\}=T^{1,0}_q\partial D$. So
$g^*(\zeta)\in\overline D\cap(q+T^{1,0}_q\partial D)=\{q\}$ for almost
every $\zeta\in\mathsf E$. As the logarithm of the modulus of the
boundary values of a nonzero bounded holomorphic function is
integrable~\cite{Garnett}, the coordinates of $g-\varphi$ vanish, and
$g=\varphi$.

If $g$ is a complex geodesic with $g(0)=\varphi(a)$ and $g'(0)$
parallel to $\varphi'(a)$, then $g'(0)=e^{i\theta}(1-|a|^2)\varphi'(a)$
by the isometry property, and with the disc automorphism
$\mathfrak a(\lambda)=\mob_{-a}(e^{i\theta}\lambda)$ the geodesic
$\varphi\circ\mathfrak a$ has the $1$-jet of $g$ at $0$. By
Lindel\"of's theorem~\cite{Garnett},
$(\varphi\circ\mathfrak a)^*=\varphi^*\circ\mathfrak a$ on
$\mathfrak a^{-1}(\mathsf E)$, so $g=\varphi\circ\mathfrak a$.

If $D$ is smooth of finite type, $\varphi^*$ exists everywhere
(Lemma~\ref{lem:holder}), and coercive points have trivial complex face
(Lemma~\ref{lem:coercive-face}). On $\Omega_\ell$, the complex tangent
line at $(\zeta,0)\in\mathfrak c$ is $\{0\}\times\C$ and
$\overline{\Omega_\ell}\cap\{|z|=1\}=\mathfrak c$ (proof of
Proposition~\ref{prop:twisted-geometry}), so the complex face at
$(\zeta,0)$ is trivial.
\end{proof}

\section{Non-coercive contact}\label{sec:noncoercive}

This section proves Theorem~\ref{thm:noncoercive-main} and gives a
global condition under which the averaged upper bound is an equality
without coercivity. We write $\C^n=\C_z\times\C^k_w$, $k=n-1$,
$\varphi(\lambda)=(\lambda,0)$ for the axis disc, and $N_\psi$ for the
right side of~\eqref{eq:general-averaged-upper} along a complex
geodesic $\psi$, as in the introduction.

\subsection{A correction criterion}
Fix integers $m_1,\ldots,m_k\ge2$ and put
$\mathcal D_t(w)=(t^{1/m_1}w_1,\ldots,t^{1/m_k}w_k)$,
$|w|_{\mathbf m}=\sum_j|w_j|^{m_j}$ and $T_{\mathbf m}=\sum_j2/m_j$ for $\mathbf m=(m_1,\ldots,m_k)$, so that
$\det_\R\mathcal D_t=t^{T_{\mathbf m}}$. Let $\mathcal W$ be the
finite-dimensional space of holomorphic polynomials $Y$ on $\C^k$ with
$Y\circ\mathcal D_t=tY$ for $t>0$, normed by $\|Y\|=\max_{|w|_{\mathbf m}=1}|Y(w)|$, so
that $|Y(w)|\le\|Y\||w|_{\mathbf m}$.

\begin{definition}\label{def:class-F}
A domain $D\subset\C^{1+k}$ is of \emph{class (F)} if:
\begin{enumerate}[(F1)]
\item $D$ is bounded and convex with $C^2$ boundary,
$\Delta\times\{0\}\subset D\subset\Delta\times\C^k$, and
$\overline D\cap(\partial\Delta\times\C^k)=\mathfrak c$;
\item there are $\eta\in(0,1)$ and, on
$U=\{1-\eta<|z|<1+\eta,\ |w|<\eta\}$, a continuous function
$\rho=|z|^2-1+\mathcal H(z,w)+R(z,w)$ with $D\cap U=\{\rho<0\}\cap U$,
where $\mathcal H$ is continuous on
$\{1-\eta\le|z|\le1+\eta\}\times\C^k$, each $\mathcal H(z,\cdot)$ is a
real polynomial with $\mathcal H(z,\mathcal D_tw)=t\mathcal H(z,w)$, and
$e_R(r)=\sup\{|R(z,w)|/|w|_{\mathbf m}:(z,w)\in U,\ 0<|w|\le r\}\to0$ as $r\to0$.
\end{enumerate}
Such a domain is \emph{correctable along $\mathfrak c$} if there is a
continuous map $\Upsilon:\overline\Delta\to\mathcal W$, holomorphic on
$\Delta$, with
\begin{equation}\label{eq:correctable}
 \mathcal H(\zeta,w)+2\operatorname{Re}\bigl(\bar\zeta\,\Upsilon(\zeta)(w)\bigr)>0
 \qquad(\zeta\in\partial\Delta,\ w\in\C^k\setminus\{0\}).
\end{equation}
\end{definition}

In the type-two case $m_j\equiv2$ one has
$\mathcal H(z,w)=w^*\mathsf A(z)w+\operatorname{Re}(w^{\mathsf T}\mathsf B(z)w)$ with
$\mathsf A$ Hermitian and $\mathsf B$ complex symmetric, and $\mathcal W$ consists of the
holomorphic quadratic forms on $\C^k$. The pointwise choice
$\Upsilon(\zeta)(w)=-\frac12\zeta\,w^{\mathsf T}\mathsf B(\zeta)w$ turns the left
side of~\eqref{eq:correctable} into $w^*\mathsf A(\zeta)w$, so a correction
always exists at a single point with $\mathsf A>0$; correctability asks that
the corrections be, up to a strict inequality, the boundary values of
one holomorphic map on the disc. By Lemma~\ref{lem:type-two} below, if
$\mathsf A(\zeta)>0$ the point $(\zeta,0)$ has multitype $(2,\ldots,2)$ and
model polynomial $\frac12\mathcal H(\zeta,\cdot)$
in~\eqref{eq:multitype-normal-form}. Hence~\eqref{eq:correctable} holds
with $\Upsilon\equiv0$ exactly when every point of $\mathfrak c$ is coercive
with multitype $(2,\ldots,2)$. After replacing $\Upsilon$ by a polynomial approximation and then
multiplying it by a sufficiently small positive constant, the
holomorphic map $z-\Upsilon(z)(w)$ is a left inverse of $\varphi$ with curved fibers,
and~\eqref{eq:correctable} says that these fibers are coercive in the
scaling limit at every point of $\mathfrak c$.

\begin{theorem}[Correction criterion]\label{thm:correction}
Let $D$ be of class (F) and correctable along $\mathfrak c$, and let
$d\mu\asymp\delta_D^\alpha dV$ be a probability measure, $\alpha>-1$.
For every compact $K\subset\Delta$, $x\in K$ and $v\ne0$,
\[
 p\log\frac{B_{p,\mu}((x,0);(v,0))}{\gamma_D((x,0);(v,0))}
 =(1+\alpha+T_{\mathbf m})\log p+O_K(1).
\]
If $\partial D$ is smooth and of finite type, then
$1+\alpha+T_{\mathbf m}\le N_\varphi(x)$, with equality when the multitype
at almost every point of $\mathfrak c$ is a permutation of $(m_j)$. If
moreover $m_j\equiv2$, every point of $\mathfrak c$ has multitype
$(2,\ldots,2)$, and the coefficient is $n+\alpha=N_\varphi(x)$.
\end{theorem}

The theorem allows every point of $\mathfrak c$ to be non-coercive. If
$m_j\equiv2$, $\partial D$ is smooth and $\mathsf A(\zeta)>0$, the model
polynomial at $(\zeta,0)$ is $\frac12\mathcal H(\zeta,\cdot)$
(Lemma~\ref{lem:type-two} below), so by the discussion of type-two
points in Section~\ref{sec:regular-reverse}, $(\zeta,0)$ is coercive if
and only if $|v^{\mathsf T}\mathsf B(\zeta)v|<v^*\mathsf A(\zeta)v$ for $v\ne0$, which
is Lempert's strong linear convexity condition at $(\zeta,0)$. In
$\C^2$ it reads $|\mathsf B(\zeta)|<\mathsf A(\zeta)$. The hypothesis $|\mathsf B|=\mathsf A$ below is
the borderline case, strongly pseudoconvex without strong linear
convexity, and correctability is then a Riemann--Hilbert-type condition
decided by a winding number.

\begin{corollary}[Winding criterion]\label{cor:winding}
Let $D\subset\C^2$ be of class (F) with $m_1=2$, so that
$\mathcal H(z,w)=\mathsf A(z)|w|^2+\operatorname{Re}(\mathsf B(z)w^2)$, and suppose that
$\mathsf A|_{\partial\Delta}$ and $\mathsf B|_{\partial\Delta}$ are smooth, with $\mathsf A>0$
and $|\mathsf B|=\mathsf A$ on $\partial\Delta$. Then every point of $\mathfrak c$ is
non-coercive of type $2$ if $\partial D$ is smooth, and $D$ is
correctable along $\mathfrak c$ if and only if
$\operatorname{wind}(\mathsf B|_{\partial\Delta})\ge-1$. In that case the
coefficient along the axis is $2+\alpha$, with bounded remainder.
\end{corollary}

\begin{lemma}\label{lem:class-F}
Let $D$ be of class (F). Then:
\begin{enumerate}[(a)]
\item $\mathcal H(\zeta,w)\ge0$ for $\zeta\in\partial\Delta$, $w\in\C^k$;
\item there is $C_{\mathcal H}$ with $|\mathcal H(z,w)|\le C_{\mathcal H}|w|_{\mathbf m}$
for $1-\eta\le|z|\le1+\eta$;
\item on $\partial D\cap U$, $\rho=0$ and
$1-|z|^2=\mathcal H(z,w)+R(z,w)\le(C_{\mathcal H}+1)|w|_{\mathbf m}$ wherever
$e_R(|w|)\le1$;
\item for every $r_1>0$ there is $r_0>0$ with
$\partial D\cap\{|w|\le r_0\}\subset\{q:\operatorname{dist}(q,\mathfrak c)<r_1\}$;
\item there are a neighborhood $\mathcal G$ of $\mathfrak c$ in $\partial D$
and $C$ with $\operatorname{vol}_{\partial D}\{q\in\mathcal G:|w(q)|_{\mathbf m}<\tau\}\le C\tau^{T_{\mathbf m}}$ for
$\tau>0$; moreover $\mathcal G$ is the graph
$\{|z|=1-g(\arg z,w)\}$ of a $C^2$ function $g\ge0$ with $g(\cdot,0)=0$,
smooth if $\partial D$ is smooth.
\end{enumerate}
\end{lemma}

\begin{proof}
(a) If $\mathcal H(\zeta,w_0)<0$, then for small $t>0$ the point
$(\zeta,\mathcal D_tw_0)$ lies in $U$ and
$\rho(\zeta,\mathcal D_tw_0)\le t(\mathcal H(\zeta,w_0)+e_R(|\mathcal D_tw_0|)|w_0|_{\mathbf m})<0$,
so it lies in $D\subset\Delta\times\C^k$, which is absurd since
$|\zeta|=1$. (b) follows from continuity, compactness of
$\{1-\eta\le|z|\le1+\eta\}\times\{w:|w|_{\mathbf m}=1\}$ and weighted homogeneity.
(c) A point of $\partial D\cap U$ is not in $D$ and is a limit of points
of $D\cap U$, so $\rho=0$ there. (d) Since
$\Delta\times\{0\}\subset D\subset\Delta\times\C^k$, the points of
$\partial D$ with $w=0$ are those of $\mathfrak c$; use compactness.
(e) The $C^2$ hypersurfaces $\partial D$ and $\partial\Delta\times\C^k$
touch along $\mathfrak c$, and $D$ lies inside the cylinder, so their
tangent hyperplanes coincide there. Hence the radial projection
$(z,w)\mapsto(z/|z|,w)$ from $\partial D$ to $\partial\Delta\times\C^k$
has invertible differential along $\mathfrak c$; it is injective on $\mathfrak c$,
so it is a $C^1$ diffeomorphism from a neighborhood $\mathcal G$ of
$\mathfrak c$ in $\partial D$ onto its image, with bounded Jacobian, and
$\mathcal G$ is a graph as stated, with $g\ge0$ because
$D\subset\Delta\times\C^k$. Thus
$\operatorname{vol}_{\partial D}\{q\in\mathcal G:|w(q)|_{\mathbf m}<\tau\}\le
C\cdot2\pi\operatorname{vol}_{2k}\{w:|w|_{\mathbf m}<\tau\}=C'\tau^{T_{\mathbf m}}$.
\end{proof}

\begin{proof}[Proof of Theorem~\ref{thm:correction}]
\emph{Step 1: a polynomial correction.} The Fej\'er means of the
boundary values of $\Upsilon$ are polynomials in $\zeta$ with
coefficients in $\mathcal W$, converging to $\Upsilon$ uniformly on
$\partial\Delta$, hence on $\overline\Delta$. The left side
of~\eqref{eq:correctable} is continuous and positive on the compact set
$\partial\Delta\times\{w:|w|_{\mathbf m}=1\}$, and replacing $\Upsilon$ by a Fej\'er mean
changes it there by at most twice the uniform distance. Hence we may
assume that $\Upsilon(z)(w)$ is a holomorphic polynomial on $\C^{1+k}$
and, by weighted homogeneity, that for some $c_*>0$
\begin{equation}\label{eq:correctable-strict}
 \mathcal H(\zeta,w)+2\operatorname{Re}(\bar\zeta\,\Upsilon(\zeta)(w))\ge2c_*|w|_{\mathbf m}
 \qquad(\zeta\in\partial\Delta,\ w\in\C^k).
\end{equation}
Put $\|\Upsilon\|_\infty=\max_{z\in\overline\Delta}\|\Upsilon(z)\|$.

\emph{Step 2: the fully corrected map near $\mathfrak c$.} Put
$u_1=z-\Upsilon(z)(w)$. On $\partial D\cap U$, by
Lemma~\ref{lem:class-F}(c),
\[
 1-|u_1|^2=\mathcal H(z,w)+R(z,w)+2\operatorname{Re}(\bar z\,\Upsilon(z)(w))
 -|\Upsilon(z)(w)|^2 .
\]
Let $e_{\mathcal H}(r)$ be the maximum, over
$1-r\le|z|\le1$ and $|\omega|_{\mathbf m}=1$, of
$|\mathcal H(z,\omega)-\mathcal H(z/|z|,\omega)|
+2|\bar z\Upsilon(z)(\omega)-\overline{(z/|z|)}\Upsilon(z/|z|)(\omega)|$;
then $e_{\mathcal H}(r)\to0$ as $r\to0$. On
$\partial D\cap U$, where $e_R(|w|)\le1$, we have
$1-|z|\le(C_{\mathcal H}+1)|w|_{\mathbf m}$, and weighted homogeneity
and~\eqref{eq:correctable-strict} give
\[
 1-|u_1|^2\ge\bigl(2c_*-e_{\mathcal H}((C_{\mathcal H}+1)|w|_{\mathbf m})-e_R(|w|)
 -\|\Upsilon\|_\infty^2|w|_{\mathbf m}\bigr)|w|_{\mathbf m}.
\]
Choose $r_0>0$ with $\partial D\cap\{|w|\le r_0\}\subset U\cap\mathcal G$
(Lemma~\ref{lem:class-F}(d),(e)) and such that the bracket is at least
$c_*$ for $|w|\le r_0$, and put $V_0=\partial D\cap\{|w|<r_0\}$. Then
$1-|u_1|^2\ge c_*|w|_{\mathbf m}$ on $V_0$.

\emph{Step 3: a scaled correction is a Carath\'eodory extremal.} For
$0<t\le1$ put $u_t=z-t\Upsilon(z)(w)=(1-t)z+tu_1$. Convexity of
$|\cdot|^2$ gives
\begin{equation}\label{eq:convex-combination}
 1-|u_t|^2\ge(1-t)(1-|z|^2)+t(1-|u_1|^2).
\end{equation}
On $V_0$ both terms are nonnegative, so $1-|u_t|^2\ge tc_*|w|_{\mathbf m}$. The
compact set $\partial D\setminus V_0$ does not meet $\{|z|=1\}$, by (F1),
so $1-|z|^2\ge c_2>0$ there, while $|u_1|\le C_3$ on $\overline D$.
By~\eqref{eq:convex-combination}, $1-|u_t|^2\ge c_2/2$ on
$\partial D\setminus V_0$ for small $t$. Fix such a $t$ and put $u=u_t$.
By the maximum principle $|u|<1$ on $D$, and
$\mathcal K_u=\{q\in\partial D:|u(q)|=1\}=\mathfrak c$. Each monomial of
$\Upsilon(z)(w)$ has degree at least $2$ in $w$, since its weighted
degree is one and every $m_j\ge2$. Hence $u(z,0)=z$ and
$du_{(x,0)}(e_z)=1$.

\emph{Step 4: the lower bound.} The axis disc and the projection $z$
show that the Kobayashi and Carath\'eodory metrics at $(x,0)$ in the
direction $e_z$ both equal $1/(1-|x|^2)$. Thus $u\in\mathcal E((0,0),e_z)$,
and Lemma~\ref{lem:collar-competitor} applies to the functions
$\mob_x\circ u$, $x\in K$. On $\partial D$,
$G_u=-\log|u|^2\ge1-|u|^2$, which is at least $tc_*|w|_{\mathbf m}$ on $V_0$ and
at least $c_2/2$ elsewhere, so Lemma~\ref{lem:class-F}(e) gives
$\operatorname{vol}_{\partial D}\{G_u<\tau\}\le C\tau^{T_{\mathbf m}}$. Lemma~\ref{lem:collar-competitor}
yields $p\log(B_{p,\mu}/\gamma_D)\ge(1+\alpha+T_{\mathbf m})\log p-C_K$ for
$x\in K$; homogeneity removes $v$.

\emph{Step 5: the upper bound.} Let $\lambda\in\Delta$ with
$s=1-|\lambda|^2<s_0$, and let $P_\lambda$ be the polydisc
$\{|z-\lambda|<cs,\ |w_j|<cs^{1/m_j}\}$. On $2P_\lambda$,
$|z|^2-|\lambda|^2\le5cs$ and $|w|_{\mathbf m}\le k(2c)^2s$, so for small $c$ and
$s_0$, $2P_\lambda\subset U$ and, by Lemma~\ref{lem:class-F}(b),
$\rho\le-s+5cs+(C_{\mathcal H}+1)k(2c)^2s<-s/2$. Hence
$\overline{2P_\lambda}\subset D$,
$\mathcal V_D(\lambda,0)\ge c's^{2+T_{\mathbf m}}$ and
$cs\le\delta_D(\lambda,0)\le s$, so $W_\varphi\ge c''s^{\alpha+T_{\mathbf m}}$ for
$s<s_0$, for either sign of $\alpha$.
By~\eqref{eq:finite-p-geodesic-entropy},
$I_{p,\mu}((0,0);e_z)^p\ge c\int_{r^2>1-s_0}r^{p+1}(1-r^2)^{\alpha+T_{\mathbf m}}dr
\ge c'p^{-1-\alpha-T_{\mathbf m}}$. For the base point $(x,0)$ apply the same
argument to $\varphi\circ\mob_{-x}$. With $s(\lambda)=1-|\lambda|^2$
and $b=\mob_{-x}(\lambda)$, one has $s(b)\asymp_K s(\lambda)$ and
$W_{\varphi\circ\mob_{-x}}(\lambda)
=W_\varphi(b)s(b)^2/s(\lambda)^2
\ge c_K s(\lambda)^{\alpha+T_{\mathbf m}}$ near the circle.
The derivative of $\mob_{-x}$ at zero is $1-|x|^2$, the reciprocal
of the axis metric in direction $e_z$. With $E_{p,\mu}\ge1$ this gives
$p\log(B_{p,\mu}/\gamma_D)\le(1+\alpha+T_{\mathbf m})\log p+C_K$.

\emph{Step 6: the multitype.} Theorem~\ref{thm:general-main}(ii) and
Steps~4--5 give $1+\alpha+T_{\mathbf m}\le N_\varphi(x)$, and the equality
case is the definition of $T$. For $m_j\equiv2$ use
Lemma~\ref{lem:type-two}.
\end{proof}

\begin{lemma}\label{lem:type-two}
Let $D$ be of class (F) with $m_j\equiv2$ and smooth boundary, and let
$\zeta\in\partial\Delta$. The point $(\zeta,0)$ has tangential multitype
$(2,\ldots,2)$ if and only if $\mathsf A(\zeta)>0$, and then, in the unitary
coordinates $z=\zeta(1+z_0)$, $w$ centered at $(\zeta,0)$, its model
polynomial in~\eqref{eq:multitype-normal-form} is
$\frac12\mathcal H(\zeta,\cdot)$. Condition~\eqref{eq:correctable}
implies $\mathsf A(\zeta)>0$.
\end{lemma}

\begin{proof}
The outer normal at $(\zeta,0)$ is $(\zeta,0)$, so $z_0$ is the complex
normal coordinate. Since $\partial D$ is smooth, the function $g$ of
Lemma~\ref{lem:class-F}(e) is smooth; as $g\ge0$ and $g(\cdot,0)=0$,
also $d_wg(\cdot,0)=0$. By Lemma~\ref{lem:class-F}(c),
$2g-g^2=\mathcal H(z,w)+R(z,w)$ with $z=(1-g)e^{i\arg z}$, and
$0\le g\le C|w|^2$; continuity of $\mathcal H$ in $z$ gives
$g(\arg\zeta,w)=\frac12\mathcal H(\zeta,w)+o(|w|^2)$. Hence the
degree-two Taylor part of $g(\arg\zeta,\cdot)$ is
$\frac12\mathcal H(\zeta,\cdot)$. The smooth local defining function
$|z|-1+g(\arg z,w)$ equals
$\operatorname{Re}z_0+\frac12\mathcal H(\zeta,w)+o(|z_0|+|w|^2)$,
because $|1+z_0|=1+\operatorname{Re}z_0+O(|z_0|^2)$ and
$g(\theta',w)-g(\theta,w)=O(|\theta'-\theta||w|^2)$. On the complex line
$\lambda\mapsto(\zeta,\lambda v)$ its degree-two part is
$\frac12(|\lambda|^2v^*\mathsf A(\zeta)v+\operatorname{Re}(\lambda^2v^{\mathsf T}\mathsf B(\zeta)v))$.
If $\mathsf A(\zeta)>0$, this is nonzero for every $v\ne0$, every line type is
$2$, and Lemma~\ref{lem:linear-multitype} gives the multitype
$(2,\ldots,2)$; the expansion above is then the normal
form~\eqref{eq:multitype-normal-form} with $H=\frac12\mathcal H(\zeta,\cdot)$.
If $v^*\mathsf A(\zeta)v=0$ for some $v\ne0$, then
$\operatorname{Re}(\lambda^2v^{\mathsf T}\mathsf B(\zeta)v)\ge0$ for all $\lambda$
by Lemma~\ref{lem:class-F}(a), so $v^{\mathsf T}\mathsf B(\zeta)v=0$ and the
line type of $v$ exceeds $2$. Finally, averaging~\eqref{eq:correctable}
over $w\mapsto e^{i\theta}w$ removes the pluriharmonic part and gives
$w^*\mathsf A(\zeta)w>0$ for $w\ne0$.
\end{proof}

\begin{proof}[Proof of Corollary~\ref{cor:winding}]
By Lemma~\ref{lem:type-two}, if $\partial D$ is smooth each point of
$\mathfrak c$ has type $2$ and model $\frac12\mathcal H(\zeta,\cdot)$, which
vanishes on a real line because $|\mathsf B|=\mathsf A$. Condition~\eqref{eq:correctable}
reads $\mathsf A|w|^2+\operatorname{Re}((\mathsf B+2\bar\zeta y)w^2)>0$ for $w\ne0$, with
$y$ continuous on $\overline\Delta$ and holomorphic on $\Delta$, that is,
\begin{equation}\label{eq:scalar-correctable}
 |2y(\zeta)+\zeta \mathsf B(\zeta)|<\mathsf A(\zeta)\qquad(\zeta\in\partial\Delta).
\end{equation}
Let $O$ be the outer function with $|O|=\mathsf A/2$ on
$\partial\Delta$~\cite[Chapter~II]{Garnett}. Since $\log \mathsf A$ is smooth,
$O$ is smooth up to the boundary and has no zeros on $\overline\Delta$,
so $\operatorname{wind}O=0$. The function
$\upsilon=-\zeta \mathsf B/(2O)$ is smooth and unimodular on $\partial\Delta$,
with $\operatorname{wind}\upsilon=1+\operatorname{wind}\mathsf B=:d$,
and~\eqref{eq:scalar-correctable} says $|f-\upsilon|<1$ on
$\partial\Delta$ for $f=y/O$.

If $d\ge0$, write $\upsilon=\zeta^de^{i\theta}$ with $\theta$ smooth and
real, let $\tilde\theta$ be its harmonic conjugate, which is smooth up
to the boundary, and put $f=c\zeta^de^{i(\theta+i\tilde\theta)}$, the
exponential being of $i$ times the holomorphic function with real
part $\mathcal P[\theta]$ and imaginary part $\tilde\theta$. On
$\partial\Delta$, $f\bar\upsilon=ce^{-\tilde\theta}\in(0,1]$ for
$0<c\le e^{\min\tilde\theta}$, so $|f-\upsilon|=1-ce^{-\tilde\theta}<1$,
and $y=Of$ is smooth up to the boundary. Conversely,
if~\eqref{eq:scalar-correctable} holds, then
$\operatorname{Re}(f\bar\upsilon)>0$ on $\partial\Delta$, so $f$ has no
zeros there and $\operatorname{wind}f=d$; by the argument principle
applied to $f(r\,\cdot\,)$, $r\uparrow1$, this is the number of zeros of
$f$ in $\Delta$, so $d\ge0$. The last assertion is
Theorem~\ref{thm:correction} with $T_{\mathbf m}=1$.
\end{proof}

\subsection{A quantitative gap for negative winding}
We now prove Theorem~\ref{thm:winding-gap}. The additional $C^4$
regularity is used only to bound the contact in a real null direction
by a fourth power. No invariance under a circle action is assumed.

\begin{lemma}[Null directions on a double cover]\label{lem:double-cover}
Suppose that $\mathsf A,\mathsf B$ are smooth on $\partial\Delta$,
$\mathsf A>0$, $|\mathsf B|=\mathsf A$, and
$\operatorname{wind}\mathsf B=-k$, $k\ge2$. There is a nowhere zero
$E\in\mathcal O(\Delta)\cap C^\infty(\overline\Delta)$ such that, for every
$a\in\Delta$, the holomorphic function
\[
 P_a(t)=(t^2-a)(1-\bar a t^2)t^{k-2}E(t)
\]
satisfies $\mathsf B(t^2)P_a(t)^2<0$ on $|t|=1$. For each compact
$K\subset\Delta$, $P_a$ and its derivative are uniformly bounded on
$\overline\Delta$, and $|P_a|$ is uniformly bounded below on a fixed
outer annulus, for $a\in K$.
\end{lemma}
\begin{proof}
The circle-valued function
$-\mathsf B(t^2)t^{2k}/\mathsf A(t^2)$ has winding zero, so it is
$e^{i\vartheta(t)}$ for a smooth real function $\vartheta$ on the
circle. Solve the Dirichlet problem with boundary data
$-\vartheta/2$, take a harmonic conjugate, and exponentiate the
holomorphic function whose imaginary part is that solution. The
resulting $E$ has $\arg E=-\vartheta/2$ on the circle. Consequently
$\mathsf B(t^2)t^{2k}E(t)^2$ is negative real. On $|t|=1$,
$(t^2-a)(1-\bar a t^2)=t^2|t^2-a|^2$, proving the assertion.
The uniform bounds follow from $\max_K|a|<1$ and compactness.
\end{proof}

\begin{proof}[Proof of Theorem~\ref{thm:winding-gap}]
The uniform graph expansion makes $D$ a domain of class (F) with
$m_1=2$. The case $\operatorname{wind}\mathsf B\ge-1$ is
Corollary~\ref{cor:winding}. Suppose that
$\operatorname{wind}\mathsf B=-k\le-2$.

\emph{Boundary depth.} By Lemma~\ref{lem:class-F}, near
$\mathfrak c$ the boundary is $|z|=1-g(\arg z,w)$, with $g\ge0$ of
class $C^4$. Put
\[
 Q_\theta(w)=2g(\theta,w)-g(\theta,w)^2,\qquad
 r_D(z,w)=|z|^2-1+Q_{\arg z}(w).
\]
Then $Q\ge0$, $D=\{r_D<0\}$ in this neighborhood, and
$-r_D\asymp\delta_D$ inside it. The quadratic Taylor polynomial of
$Q_\theta$ at $w=0$ is $\mathcal H(e^{i\theta},w)$: this follows by
substituting the boundary graph in (F2), using
$g(\theta,w)=O(|w|^2)$ and the uniform $o(|w|^2)$ remainder.
If $v$ is a unit null vector of this quadratic form, then, uniformly
in $\theta$ and either choice of $v$,
\begin{equation}\label{eq:null-quartic-upper}
 0\le Q_\theta(Xv+iYv)\le C(X^4+Y^2)\qquad (X,Y\in\R\text{ small}).
\end{equation}
Indeed, the quadratic part is a positive multiple of $Y^2$.
Nonnegativity on the line $Y=0$ makes its pure cubic coefficient
zero. Every remaining cubic term is bounded by a constant times
$X^4+Y^2$, by Young's inequality, and the fourth-order remainder has
the same bound. Compactness makes the constants uniform. This
argument uses Taylor expansion through order three with a bounded
fourth-order remainder, not a positive quartic coefficient.

\emph{An interior family.} Take $P_a$ from
Lemma~\ref{lem:double-cover} and write
\[
 \Psi_a(t,\tau)=(t^2,\tau P_a(t)),\qquad s(t)=1-|t|^4,\qquad
 X_h=\{\tau:|\operatorname{Re}\tau|<c h^{1/4},\
                   |\operatorname{Im}\tau|<c h^{1/2}\}.
\]
Fix $K\Subset\Delta$. On a fixed outer annulus, the change from
$P_a(t/|t|)$ to $P_a(t)$ is $O_K(s)$, and the former vector is null
for $\mathcal H((t/|t|)^2,\cdot)$. Thus
\eqref{eq:null-quartic-upper} gives, uniformly for $a\in K$,
\[
 Q_{\arg(t^2)}(\tau P_a(t))
 \le C_K\bigl((\operatorname{Re}\tau)^4+
                 (\operatorname{Im}\tau)^2+s^2|\tau|^2\bigr).
\]
For a sufficiently small fixed $c>0$, if $\tau\in X_{s(t)}$ then
\begin{equation}\label{eq:cover-depth}
 \tfrac12s(t)\le-r_D(\Psi_a(t,\tau))\le s(t).
\end{equation}
For $p$ sufficiently large, $\tau\in X_{1/p}$ therefore makes
$\Psi_a(\cdot,\tau)$ map $\{|t|^4<1-1/p\}$ into $D$.
On the outer annulus this follows from \eqref{eq:cover-depth}; on the
remaining closed disc it follows because $\tau P_a(t)\to0$ uniformly
and the corresponding segment of the axis is compactly contained in
$D$. These choices are uniform on $K$.

On the shell $S_p=\{1/p<s(t)<2/p\}$, the complex Jacobian of
$\Psi_a$ is $2tP_a(t)$, bounded above and below away from zero. The
map has multiplicity at most two there: $z=t^2$ has at most two
solutions and then $w=\tau P_a(t)$ determines $\tau$. The area formula,
\eqref{eq:cover-depth}, and $d\mu\asymp\delta_D^\alpha dV$ give
\begin{equation}\label{eq:cover-area}
 \int_D|f|^p\,d\mu\ge c_Kp^{-\alpha}
       \int_{X_{1/p}}\int_{S_p}|f(\Psi_a(t,\tau))|^p\,dA(t)dA(\tau).
\end{equation}
This is valid also for $-1<\alpha<0$, since the depth on the image is
comparable to $1/p$ in both directions.

\emph{Preserving the jet.} Normalize $f(a,0)=0$, $f_z(a,0)=1$, and
write $c_f=f_w(a,0)$. The function
\[
 G_\tau(t)=\frac{f(t^2,\tau P_a(t))}{\mob_a(t^2)}
\]
is holomorphic on $\{|t|^4<1-1/p\}$. At a root $b^2=a$ it has value
\begin{equation}\label{eq:cover-jet}
 G_\tau(b)=(1-|a|^2)(1+c_f\tau C_a),\qquad
 C_a=(1-|a|^2)b^{k-2}E(b).
\end{equation}
At $a=b=0$ interpret $b^0=1$ when $k=2$. The division is still
legitimate there: $P_0(t)=t^kE(t)$, so the numerator vanishes to order
at least two and its coefficient of $t^2$ gives precisely
\eqref{eq:cover-jet}. At nonzero $b$ the identity follows by cancelling
$t^2-a$; thus no inverse branch at $a=0$ is used.

On $S_p$, $|\mob_a(t^2)|^p\ge c_K$. The Poisson inequality for
$|G_\tau|^p$, evaluated at $b$ on each circle $|t|=r$ in $S_p$, has
Poisson kernel bounded by $C_K$. Since the shell has area comparable
to $p^{-1}$, it follows that
\[
 \int_{S_p}|f(\Psi_a(t,\tau))|^p\,dA(t)
 \ge c_Kp^{-1}(1-|a|^2)^p|1+c_f\tau C_a|^p.
\]
The rectangle $X_{1/p}$ is centrally symmetric, and convexity gives
$|1+y|^p+|1-y|^p\ge2$ for $p\ge1$. Integrating in $\tau$ in
\eqref{eq:cover-area} therefore proves
\[
 \|f\|_{L^p(\mu)}^p\ge c_K(1-|a|^2)^p p^{-7/4-\alpha}.
\]
The axis is a retract, so $\gamma_D((a,0);e_z)=(1-|a|^2)^{-1}$.
Taking the infimum over $f$ and using $E_{p,\mu}\ge1$ proves the
asserted upper bound, uniformly for $a\in K$.
\end{proof}

\begin{corollary}[Quartic contact gives equality]\label{cor:quartic-equality}
Under the negative-winding hypotheses of
Theorem~\ref{thm:winding-gap}, suppose additionally that the boundary
graph above satisfies
\begin{equation}\label{eq:null-quartic-lower}
 Q_\theta(Xv+iYv)\ge c(X^4+Y^2)
\end{equation}
uniformly for unit null vectors $v$ of $\mathcal H(e^{i\theta},\cdot)$
and small real $X,Y$. Then, for every $K\Subset\Delta$,
\[
 p\log\frac{B_{p,\mu}}{\gamma_D}((a,0);e_z)
       =(\tfrac74+\alpha)\log p+O_K(1),\qquad a\in K.
\]
\end{corollary}
\begin{proof}
Put $s=1-|z|^2$. By (F1), the full fiber over $z$ lies in the boundary
graph neighborhood when $s$ is small. Its area is at most $Cs^{3/4}$
by \eqref{eq:null-quartic-lower}. The maximum $d_*$ of $\delta_D$ on
this convex fiber satisfies $cs\le d_*\le Cs$: use the axis point for
the lower bound, and $D\subset\Delta\times\C$ for the upper bound.
For $\alpha\ge0$ this immediately bounds its weighted area by
$Cs^{\alpha+3/4}$. For $-1<\alpha<0$, concavity of $\delta_D$ on a
convex domain implies, by contraction toward a maximum point,
\[
 |\{w\text{ in the fiber}:\delta_D(z,w)<h\}|
       \le 2(h/d_*)\,|\{w:(z,w)\in D\}|\quad(0<h<d_*).
\]
Layer integration consequently bounds the weighted area by
$[1+2(-\alpha)/(1+\alpha)]d_*^\alpha$ times its area. Thus the
pushforward under $z$ has density at most $Cs^{\alpha+3/4}$ near the
circle. The part over a fixed smaller disc contributes an
exponentially small amount to the moment of $\mob_a(z)$, uniformly
for $a\in K$. Integrating on the collar gives
$\int_D|\mob_a(z)|^p\,d\mu\le C_Kp^{-7/4-\alpha}$.
This jet competitor, together with $E_{p,\mu}^p\le C_K$
(Lemma~\ref{lem:basic-comparison}), proves the opposite bound.
\end{proof}

\begin{corollary}[A gap under H\"older control of second derivatives]
\label{cor:holder-gap}
In Theorem~\ref{thm:winding-gap}, replace $C^4$ regularity near
$\mathfrak c$ by $C^{2,\eta}$, $0<\eta\le1$, retaining the smooth
quadratic data. If $d\le-2$, then, uniformly for $a\in K\Subset\Delta$,
\[
 p\log(B_{p,\mu}/\gamma_D)((a,0);e_z)
 \le\left(\tfrac32+\alpha+\frac1{2+\eta}\right)\log p+C_K.
\]
In particular the upper coefficient is strictly less than $2+\alpha$.
\end{corollary}
\begin{proof}
Taylor's formula with a H\"older remainder gives
$0\le Q_\theta(Xv+iYv)\le C(Y^2+|X|^{2+\eta})$ in null coordinates,
uniformly in $\theta$. Use the same $P_a,\Psi_a$ as in
Lemma~\ref{lem:double-cover}, replacing $X_h$ by
$\{|\operatorname{Re}\tau|<ch^{1/(2+\eta)},\
|\operatorname{Im}\tau|<ch^{1/2}\}$.
The radial error is still $O_K(s|\tau|)$; its contribution to the
quadratic term is $O_K(s^2|\tau|^2)$ and to the remainder is
$O_K(s^{2+\eta}|\tau|^{2+\eta})$. Both are absorbed in
\eqref{eq:cover-depth}. Thus \eqref{eq:cover-area} and
\eqref{eq:cover-jet} hold unchanged. The rectangle now has area
$4c^2p^{-1/2-1/(2+\eta)}$. The same Poisson inequality and symmetric
averaging give
$\|f\|_p^p\ge c_K(1-|a|^2)^p
p^{-3/2-\alpha-1/(2+\eta)}$, which proves the assertion.
\end{proof}

\subsection{Twisted domains}
Throughout, $\ell\ge1$, $0<\eps\le1/(2\ell^2)$,
$\rho_\ell(z,w)=|z|^2-1+\eps(|w|^2-\operatorname{Re}(\bar z^\ell w^2))+|w|^4$,
$\Omega_\ell$ is as in Theorem~\ref{thm:noncoercive-main}, and
$\mathcal Z=\overline\Delta\times\{|w|\le\sqrt2\}$. On $\mathfrak c$ the zero
line $\zeta^{\ell/2}\R$ of the model polynomial turns through the angle
$\ell\pi$ as $\zeta$ runs once around the circle.

\begin{figure}[tbp]
\centering
\includegraphics[width=\textwidth]{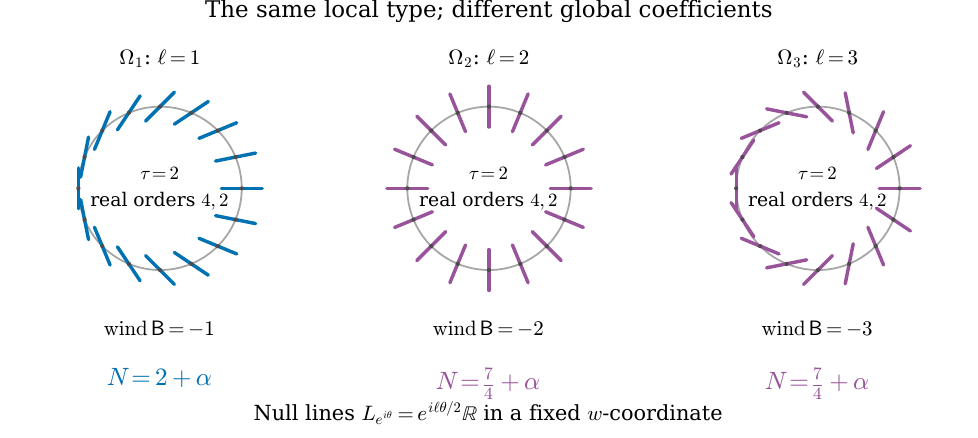}
\caption{The null-line fields of $\Omega_1$, $\Omega_2$, and $\Omega_3$.
Each circle parametrizes $\zeta\in\partial\Delta$; the segments are
unoriented lines in the fixed $w$-coordinate, translated to the
corresponding parameter points for display. For odd $\ell$ the unit
vector changes sign after a full turn, but the unoriented line is
single-valued. The values of $N$ are
those proved in Theorem~\ref{thm:noncoercive-main}, for axis jets and
$d\mu\asymp\delta^\alpha dV$; they hold for every allowed $\eps$.
The general negative-winding theorem gives an upper bound, with
equality under the additional quartic lower bound.}
\label{fig:null-winding}
\end{figure}

Figure~\ref{fig:null-winding} displays this null-line field for the first
three domains. Their local real contact orders agree, whereas the
coefficient changes between $\ell=1$ and $\ell=2$.

\begin{lemma}\label{lem:twisted-convex}
$\rho_\ell$ is convex on $\mathcal Z$.
\end{lemma}

\begin{proof}
For $X=(\xi,\upsilon)\in\C^2$ write
$\nabla^2f[X]=\frac{d^2}{dt^2}f(Z+tX)|_{t=0}$. Then $\nabla^2|z|^2[X]=2|\xi|^2$,
$\nabla^2|w|^2[X]=2|\upsilon|^2$,
$\nabla^2|w|^4[X]=8(\operatorname{Re}\bar w\upsilon)^2+4|w|^2|\upsilon|^2$, and
\[
 \nabla^2(\bar z^\ell w^2)[X]=\ell(\ell-1)\bar z^{\ell-2}\bar\xi^2w^2
 +4\ell\bar z^{\ell-1}\bar\xi w\upsilon+2\bar z^\ell\upsilon^2 .
\]
For $|z|\le1$ and $|w|^2\le2$, since $2\eps\ell(\ell-1)\le(\ell-1)/\ell<1$,
\[
 \nabla^2\rho_\ell[X]\ge(2-\eps\ell(\ell-1)|w|^2)|\xi|^2-4\ell\eps|w||\xi||\upsilon|
 +4|w|^2|\upsilon|^2
 \ge|\xi|^2-4\ell\eps|w||\xi||\upsilon|+4|w|^2|\upsilon|^2,
\]
and Young's inequality also gives
$\nabla^2\rho_\ell[X]\ge\frac12|\xi|^2+
(4-8\ell^2\eps^2)|w|^2|\upsilon|^2$.
Here $4-8\ell^2\eps^2\ge2$, so the real Hessian is positive when
$w\ne0$.
\end{proof}

\begin{proposition}[Geometry of $\Omega_\ell$]\label{prop:twisted-geometry}
$\Omega_\ell=\{Z\in\mathcal Z:\rho_\ell(Z)<0\}$ is a bounded convex
domain with real-algebraic strongly pseudoconvex boundary,
$\Omega_\ell\subset\Delta\times\{|w|<1\}$, and $-\rho_\ell\asymp\delta_{\Omega_\ell}$.
It is of class (F) with $k=1$, $m_1=2$,
$\mathcal H=\eps(|w|^2-\operatorname{Re}(\bar z^\ell w^2))$ and
$R=|w|^4$, so $\mathsf A=\eps$, $\mathsf B=-\eps\bar z^\ell$ and
$\operatorname{wind}(\mathsf B|_{\partial\Delta})=-\ell$. Every point
$(\zeta,0)\in\mathfrak c$ is strongly pseudoconvex and non-coercive, of
tangential type $2$, with model polynomial
$\eps(\operatorname{Im}(\bar\zeta^{\ell/2}w))^2$ for either square root.
If $g$ is a complex geodesic of $\Omega_\ell$ with $g(\lambda_0)=(a,0)$
and $g'(\lambda_0)\in\C^*e_z$, then $g=\varphi\circ \mathfrak a$ for a disc
automorphism $\mathfrak a$. In particular $N_\varphi\equiv2+\alpha$, and so
is the infimum of $N_\psi(\lambda_0)$ over complex geodesics $\psi$
with $\psi(\lambda_0)=(a,0)$ and $\psi'(\lambda_0)\in\C^*e_z$.
\end{proposition}

\begin{proof}
On $\mathcal Z$,
$|w|^2-\operatorname{Re}(\bar z^\ell w^2)\ge|w|^2(1-|z|^\ell)\ge0$, so
$\rho_\ell\ge|z|^2-1+|w|^4$. On $\{|z|=1\}\cap\mathcal Z$,
$\rho_\ell\ge|w|^4$, with equality only at $w=0$; on
$\{|w|=\sqrt2\}\cap\mathcal Z$, $\rho_\ell\ge3$. Hence
$\Omega_\ell=\{Z\in\mathcal Z:\rho_\ell(Z)<0\}$, which is convex by
Lemma~\ref{lem:twisted-convex}, contains $0$ and lies in
$\Delta\times\{|w|<1\}$. If $\rho_\ell(Z)=0$ with $Z\in\mathcal Z$,
convexity on $[0,Z]$ gives $d\rho_\ell(Z)[Z]\ge\rho_\ell(Z)-\rho_\ell(0)=1$,
so $\partial\Omega_\ell=\{Z\in\mathcal Z:\rho_\ell=0\}$ is a compact
real-algebraic hypersurface with $d\rho_\ell\ne0$, and
$\overline{\Omega_\ell}\cap\{|z|=1\}=\mathfrak c$. For $|z|\ge1$ and small
$|w|$,
$\rho_\ell\ge(|z|-1)\bigl(|z|+1-\eps|w|^2\sum_{j<\ell}|z|^j\bigr)+|w|^4\ge0$,
so $\Omega_\ell=\{\rho_\ell<0\}$ near every boundary point; thus
$\partial\Omega_\ell$ is smooth, $-\rho_\ell\asymp\delta_{\Omega_\ell}$,
and (F) holds. The real Hessian bound in
Lemma~\ref{lem:twisted-convex} gives strong pseudoconvexity wherever
$w\ne0$. At the remaining boundary points $(\zeta,0)$, the
tangential Levi form is $\eps|dw|^2>0$. Thus the whole boundary is
strongly pseudoconvex, in particular of finite type.

At $(\zeta,0)$ the complex tangent line is $\{0\}\times\C$, the Levi
form there is $\partial_w\partial_{\bar w}\rho_\ell=\eps>0$, and
$\rho_\ell(\zeta,\lambda)=2\eps(\operatorname{Im}(\bar\zeta^{\ell/2}\lambda))^2+|\lambda|^4$.
By Lemma~\ref{lem:type-two} the type is $2$ and the model polynomial is
$\frac12\mathcal H(\zeta,\cdot)=\eps(\operatorname{Im}(\bar\zeta^{\ell/2}w))^2$,
which vanishes on $\zeta^{\ell/2}\R$; so $(\zeta,0)$ is not coercive,
and~\eqref{eq:real-tangential-contact} fails.

Since $\Delta\times\{0\}\subset\Omega_\ell\subset\Delta\times\C$, the
projection $z$ is a left inverse of $\varphi$, so $\varphi$ is a complex
geodesic and the Kobayashi and Carath\'eodory metrics at $(a,0)$ in the
direction $e_z$ equal $1/(1-|a|^2)$. If $g$ is as stated, then
$z\circ g:\Delta\to\Delta$ maps $\lambda_0$ to $a$ with
$|(z\circ g)'(\lambda_0)|=(1-|a|^2)/(1-|\lambda_0|^2)$ by the isometry
property, so it is an automorphism $\mathfrak a$ by the equality case of the
Schwarz--Pick lemma, and $g\circ \mathfrak a^{-1}=(\lambda,\tilde g(\lambda))$ with
$\tilde g$ holomorphic and $|\tilde g(\lambda)|^4<1-|\lambda|^2$. Thus
$\tilde g$ has radial limits zero, and $\tilde g\equiv0$ (this is also
a case of Corollary~\ref{cor:uniqueness}). Finally
$T\equiv1$ on $\mathfrak c$, and $N_{\varphi\circ \mathfrak a}=N_\varphi\circ \mathfrak a$
by conformal invariance of harmonic measure.
\end{proof}

\begin{lemma}[The fiber exponent]\label{lem:twisted-fibers}
Let $\mu_0=c_0(-\rho_\ell)^\alpha dV$ be a probability measure on
$\Omega_\ell$, $\alpha>-1$. Then $z_*\mu_0=c_0s^{\alpha+3/4}J(s)\,dA$
with a function $J\asymp1$ of $s\in(0,1]$, and $z_*\mu\asymp s^{\alpha+3/4}dA$ for every
probability measure $d\mu\asymp\delta_{\Omega_\ell}^\alpha dV$.
\end{lemma}

\begin{proof}
The fiber over $z=re^{i\theta}$ is $\{w:\mathcal Q_z(w)<s\}$ with
$\mathcal Q_z(w)=\eps(|w|^2-\operatorname{Re}(\bar z^\ell w^2))+|w|^4$, and
$-\rho_\ell=s-\mathcal Q_z$ on it. With $w=e^{i\ell\theta/2}(x+iy)$ one has
$\mathcal Q_z=\eps(1-r^\ell)x^2+\eps(1+r^\ell)y^2+(x^2+y^2)^2$, and the
substitution $x=s^{1/4}X$, $y=s^{1/2}Y$ gives
\[
 \int_{\{\mathcal Q_z<s\}}(-\rho_\ell)^\alpha\,dV(w)=s^{\alpha+3/4}J(s),\qquad
 J(s)=\int_{\{\Phi_s<1\}}(1-\Phi_s)^\alpha\,dX\,dY,
\]
\[
 \Phi_s=\eps(1-r^\ell)s^{-1/2}X^2+\eps(1+r^\ell)Y^2+(X^2+s^{1/2}Y^2)^2 .
\]
Since $s/2\le1-r\le1-r^\ell\le\ell s$ and $s\le1$,
$\eps Y^2+X^4\le\Phi_s\le\eps\ell X^2+2\eps Y^2+2X^4+2Y^4$, so
$\{\Phi_s<1\}$ and $\{\Phi_s<1/2\}$ contain, and are contained in, fixed
bounded sets of positive area. This gives $J\asymp1$ for $\alpha\ge0$,
and the lower bound for $\alpha<0$. For $\alpha<0$, $\Phi_s$ is convex
and vanishes at $0$, so on each ray $\Phi_s(t\vartheta)\le t/t_\vartheta$
up to the exit radius $t_\vartheta$ of $\{\Phi_s<1\}$, and polar
coordinates give
$J\le2|\{\Phi_s<1\}|/((1+\alpha)(2+\alpha))$. The last assertion follows
from $\mu\asymp\mu_0$, since $-\rho_\ell\asymp\delta_{\Omega_\ell}$.
\end{proof}

\begin{proof}[Proof of Theorem~\ref{thm:noncoercive-main}]
The geometric assertions are Proposition~\ref{prop:twisted-geometry}. By
Lemma~\ref{lem:basic-comparison} we may replace $\mu$ by $\mu_0$, and by
homogeneity take $v=1$.

(i) For $\ell=1$, the constant map $\Upsilon\equiv\frac\eps2w^2$ gives
$\mathcal H(\zeta,w)+2\operatorname{Re}(\bar\zeta\frac\eps2w^2)=\eps|w|^2$,
so $\Omega_1$ is correctable along $\mathfrak c$, and
Theorem~\ref{thm:correction} applies; equivalently,
$\operatorname{wind}\mathsf B=-1$ in Corollary~\ref{cor:winding}. In fact
$u=z-\frac\eps2w^2$ satisfies
$1-|u|^2=-\rho_1+\eps|w|^2+(1-\eps^2/4)|w|^4$, so it maps $\Omega_1$
into $\Delta$ directly.

(ii) For $\ell\ge2$, Theorem~\ref{thm:winding-gap}(ii) gives the
upper bound, since $\operatorname{wind}\mathsf B=-\ell$. For the
opposite bound put $s=1-|z|^2$. By Lemma~\ref{lem:twisted-fibers},
$z_*\mu_0\le Cs^{\alpha+3/4}dA$. For $a\in K$, the function
$f=(1-|a|^2)\mob_a(z)$ is a normalized jet competitor and
$1-|\mob_a(z)|^2\ge c_Ks$. Hence
\[
 (\gamma I_{p,\mu_0})^p\le\int_\Delta|\mob_a(z)|^p\,d(z_*\mu_0)
 \le C_K\int_0^1e^{-c_Kps/2}s^{\alpha+3/4}\,ds
 \le C_Kp^{-7/4-\alpha}.
\]
Using $E_{p,\mu_0}^p\le C_K$ from
Lemma~\ref{lem:basic-comparison} proves the lower bound. Both constants
are uniform for $a\in K$; the argument applies to every integer
$\ell\ge2$.
\end{proof}

\begin{proposition}[Exact finite-$p$ formulas for the illustrations]
\label{prop:finite-p-illustrations}
Let $m\ge1$ be an integer, $\alpha>-1$, and put
\[
 E_m=\{(z,w)\in\C^2:|z|^2+|w|^{2m}<1\},\qquad
 d\mu_m=Z_m^{-1}(1-|z|^2-|w|^{2m})^\alpha dV,
\]
where $Z_m$ normalizes the measure. For every $p\ge1$, at $(0;e_z)$,
\begin{equation}\label{eq:illustration-ellipsoid}
 \left(\frac{B_{p,\mu_m}}{\gamma_{E_m}}\right)^p
 =\frac{\Gamma(p/2+N_m+1)}{\Gamma(N_m+1)\Gamma(p/2+1)},
 \qquad N_m=1+\alpha+1/m.
\end{equation}
On $\Omega_2$, take $0<\eps\le1/8$ and normalized Euclidean volume.
At $(0;e_z)$ one has
\begin{equation}\label{eq:illustration-twisted}
 \left(\frac{B_p}{\gamma_{\Omega_2}}\right)^p
 =\frac{\int_0^1 a_\eps(s)\,ds}
        {\int_0^1(1-s)^{p/2}a_\eps(s)\,ds},
\end{equation}
where the fiber area is given by
\[
 a_\eps(s)=\frac12\int_0^{2\pi}
 \frac{2s\,d\theta}{\sqrt{b_s(\theta)^2+4s}+b_s(\theta)},\qquad
 b_s(\theta)=\eps\bigl(1-(1-s)\cos(2\theta)\bigr).
\]
\end{proposition}
\begin{proof}
Both domains and their stated measures are invariant under
$R_t(z,w)=(e^{it}z,e^{it}w)$ and $S(z,w)=(z,-w)$. The operators
\[
 \Pi_j f=\frac1{2\pi}\int_0^{2\pi}e^{-ijt}f\circ R_t\,dt
 \quad(j=0,1),\qquad Af=\tfrac12(f+f\circ S)
\]
are $L^p$ contractions for $p\ge1$, by Jensen's inequality and invariance
of the measure. Taylor expansion near $0$, followed by the identity
theorem on the connected domain, gives
$\Pi_0f=f(0)$ and $A\Pi_1f=f_z(0)z$ on the whole domain.
Thus $1$ is a value extremal and $z$ minimizes the norm under
$f(0)=0$, $f_z(0)=1$; the measures have total mass one, so $E_p(0)=1$.
Both domains retract onto the axis disc, giving $\gamma(0;e_z)=1$
and $B_p(0;e_z)=\|z\|_p^{-1}$.

On $E_m$, integration in $w$ gives the base density proportional to
$(1-|z|^2)^{N_m-1}dA(z)$. Normalization and the beta integral give
\[
 \|z\|_{L^p(\mu_m)}^p
 =N_m\int_0^1(1-s)^{p/2}s^{N_m-1}\,ds
 =\frac{\Gamma(N_m+1)\Gamma(p/2+1)}{\Gamma(p/2+N_m+1)}.
\]
This proves \eqref{eq:illustration-ellipsoid} for all $\alpha>-1$.
On $\Omega_2$, rotate $w$ by $\arg z$ (arbitrarily if $z=0$), and put
$s=1-|z|^2$. The fiber is
$\{\eps s x^2+\eps(2-s)y^2+(x^2+y^2)^2<s\}$.
Its squared polar radius is the positive root of
$t^2+b_s(\theta)t=s$, giving the formula for $a_\eps(s)$.
Integration in $z$ proves \eqref{eq:illustration-twisted}; the common
factor $\pi$ cancels between the moment and the total volume.
\end{proof}

These exact identities use the specified measures. Replacing them by
comparable densities preserves the limiting coefficients, but need not
preserve the finite-$p$ values.

\begin{proposition}[Perturbations without the weighted circle action]
\label{prop:perturbed-twisting}
Let $\ell\ge2$, $0<\eps\le1/(2\ell^2)$, and $0<|\eta|\le10^{-4}$.
On $\mathcal Z=\overline\Delta\times\{|w|\le\sqrt2\}$ put
\[
 \rho_{\ell,\eta}=\rho_\ell+\eta\operatorname{Re}z\,|w|^6,
 \qquad D_{\ell,\eta}=\{Z\in\mathcal Z:\rho_{\ell,\eta}(Z)<0\}.
\]
This is a bounded convex domain with smooth real-algebraic, strongly
pseudoconvex boundary. It has the same quadratic data along
$\mathfrak c$ as $\Omega_\ell$. For every probability measure
$d\mu\asymp\delta_{D_{\ell,\eta}}^\alpha dV$, $\alpha>-1$, and
$K\Subset\Delta$,
\[
 p\log(B_{p,\mu}/\gamma_{D_{\ell,\eta}})((a,0);e_z)
      =(\tfrac74+\alpha)\log p+O_K(1).
\]
For $\eta\ne0$ the weighted circle action
$(z,w)\mapsto(e^{2i\theta}z,e^{i\ell\theta}w)$ does not preserve this
domain.
\end{proposition}
\begin{proof}
Write $r=|w|$ and $X=(\xi,\upsilon)$. The calculation in
Lemma~\ref{lem:twisted-convex}, followed by Young's inequality, gives
\[
 \nabla^2\rho_\ell[X]\ge\tfrac12|\xi|^2+
                 (4-8\ell^2\eps^2)r^2|\upsilon|^2
 \ge\tfrac12|\xi|^2+\tfrac72r^2|\upsilon|^2.
\]
On $\mathcal Z$ the absolute value of the Hessian of
$\operatorname{Re}z\,|w|^6$ is at most
$12r^5|\xi||\upsilon|+30r^4|\upsilon|^2
\le24|\xi|^2+84r^2|\upsilon|^2$. Thus $\rho_{\ell,\eta}$ is convex,
with positive real Hessian when $w\ne0$.
Moreover $\rho_{\ell,\eta}\ge |z|^2-1+(1-2|\eta|)|w|^4$.
It follows that the artificial cylinder boundary is disjoint from
$\{\rho_{\ell,\eta}<0\}$, and its zero set there consists only of
$\mathfrak c$. Convexity and $\rho_{\ell,\eta}(0)=-1$ give
$d\rho_{\ell,\eta}(Z)[Z]\ge1$ on the zero set. Near $\mathfrak c$,
$\partial_r\rho_{\ell,\eta}(re^{i\theta},w)>0$ and
$\rho_{\ell,\eta}(e^{i\theta},w)\ge0$; hence the full local sublevel
set lies inside $|z|<1$. Thus the boundary is a regular polynomial
level set everywhere. It is strongly pseudoconvex off $w=0$ by the
Hessian bound; at $(\zeta,0)$ its tangential Levi form is
$\eps|dw|^2$.

It remains to check the hypothesis of Corollary~\ref{cor:quartic-equality}.
Write $z=r\zeta$ and $w=\zeta^{\ell/2}(X+iY)$ locally. At a boundary
point, $s=1-r^2=Q_{\arg\zeta}(w)$ and
\[
 s=\eps(1-r^\ell)X^2+\eps(1+r^\ell)Y^2+|w|^4
                      +\eta r\operatorname{Re}\zeta\,|w|^6.
\]
Since $1-r^\ell\le C_\ell s$, for small $w$ this identity gives
$c(X^4+Y^2)\le s\le C(X^4+Y^2)$, after absorbing $C_\ell sX^2$.
The estimates are uniform in $\zeta$. The quadratic data are unchanged,
so their winding is $-\ell$, and the corollary applies. Finally the
unperturbed defining function is invariant under the stated action,
whereas the added term changes with $\operatorname{Re}(e^{2i\theta}z)$.
At a boundary point with $zw\ne0$ this is not constant in $\theta$,
which proves the last assertion.
\end{proof}

\begin{proposition}[The $C^2$ obstruction]\label{prop:C2-obstruction}
The negative-winding conclusion of Theorem~\ref{thm:winding-gap}
need not hold with only $C^2$ boundary. Fix $0<\eps\le1/8$, put
$\Phi(0)=0$, $\Phi(t)=t/\log(e^2/t)$ for $0<t\le2$, and define on
$\mathcal Z$
\[
 \widehat\rho(z,w)=\rho_2(z,w)+\Phi(|w|^2),\qquad
 \widehat D=\{Z\in\mathcal Z:\widehat\rho(Z)<0\}.
\]
Then $\widehat D$ is bounded and convex with $C^2$, strongly
pseudoconvex boundary, is of class (F), and has
$\mathsf A=\eps$, $\mathsf B(\zeta)=-\eps\bar\zeta^2$ on the circle.
For $d\mu\asymp\delta_{\widehat D}^\alpha dV$ a probability measure,
$\alpha>-1$,
\begin{equation}\label{eq:C2-loglog}
 p\log(B_{p,\mu}/\gamma_{\widehat D})((0,0);e_z)
   =(2+\alpha)\log p-\tfrac12\log\log p+O(1).
\end{equation}
\end{proposition}
\begin{proof}
With $L=\log(e^2/t)$, one has
$\Phi'=L^{-1}+L^{-2}>0$ and
$\Phi''=t^{-1}(L^{-2}+2L^{-3})>0$ for $0<t\le2$.
Thus $\Phi(|w|^2)$ is convex, extends as a $C^2$ function at zero with
zero Hessian there, and is $o(|w|^2)$. Convexity follows from
Lemma~\ref{lem:twisted-convex}. The same radial argument as in
Proposition~\ref{prop:twisted-geometry} makes the zero set regular,
with $-\widehat\rho\asymp\delta_{\widehat D}$; the positive added term
introduces no cylinder boundary. The real Hessian is positive off
$w=0$, and the Levi form at $(\zeta,0)$ is $\eps|dw|^2$. Along the
null direction $w=\zeta x$, the boundary graph satisfies
$Q_{\arg\zeta}(\zeta x)=(x^4+\Phi(x^2))/(1-\eps x^2)$. Its
quadratic Taylor polynomial on this line is zero, but the graph is
not $O(|x|^{2+\eta})$ for any $\eta>0$. Thus the boundary is not
$C^{2,\eta}$ for any $\eta>0$.

Use the normalized measure $d\mu_0=c(-\widehat\rho)^\alpha dV$;
comparison with it changes the left side of \eqref{eq:C2-loglog} by
$O(1)$. For $z=r\zeta$, rotate $w=\zeta(x+iy)$ and put $s=1-r^2$.
The fiber is defined by
\[
 \eps s x^2+\eps(2-s)y^2+(x^2+y^2)^2+\Phi(x^2+y^2)<s.
\]
For small $s$ its area is comparable to
$s\sqrt{L_s}$, where $L_s=\log(e^2/s)$. Indeed,
$\Phi(t)<s$ implies $t\le C sL_s$, since $\Phi$ is increasing and
$\Phi(CsL_s)>s$ for fixed $C>2$ and small $s$; the quadratic term
also gives $|y|\le C\sqrt s$. Conversely a rectangle
$|x|<c\sqrt{sL_s}$, $|y|<c\sqrt s$ lies in the half-depth region for
small fixed $c$: monotonicity gives
$\Phi(c^2sL_s+c^2s)\le Cc^2s$, and the other terms obey the same
bound. Its area is comparable to $s\sqrt{L_s}$.
For negative $\alpha$, apply the concave-depth contraction estimate
from Corollary~\ref{cor:quartic-equality} to
$d(w)=-\widehat\rho(z,w)$, whose maximum is $s$; for nonnegative
$\alpha$ use $d\le s$. The half-depth rectangle supplies the lower
bound in either case. The pushforward density under $z$ is therefore
comparable to $s^{1+\alpha}\sqrt{L_s}$.

The group $(z,w)\mapsto(e^{2i\theta}z,\pm e^{2i\theta}w)$ preserves
$\widehat D$ and $\mu_0$. Haar averaging in its character of weight
$2$ sends every competitor $f(0)=0$, $f_z(0)=1$ to $z$: the surviving
Taylor monomials satisfy $2j+2i=2$ and $i$ is even. Averaging in the
trivial character likewise gives $E_{p,\mu_0}(0)=1$. Consequently
$B_{p,\mu_0}(0;e_z)=\|z\|_p^{-1}$ and
\[
 \|z\|_p^p\asymp\int_0^1(1-s)^{p/2}s^{1+\alpha}\sqrt{\log(e^2/s)}\,ds
             \asymp p^{-2-\alpha}\sqrt{\log p}.
\]
For the last upper bound put $t=ps$, use
$(1-s)^{p/2}\le e^{-ps/2}$ and
$\log(e^2p/t)\le\log p+2+|\log t|$; the lower bound follows by
integrating over $1/p<s<2/p$. Since $\gamma_{\widehat D}(0;e_z)=1$,
this proves \eqref{eq:C2-loglog}.
\end{proof}

\begin{corollary}[Pointwise unitarily equivalent data do not determine the coefficient]
\label{cor:noncoercive-answer}
For $\ell\ge2$ and every $\alpha>-1$, the coefficient of $\Omega_\ell$
at $((0,0);(v,0))$ is $\frac74+\alpha$, strictly smaller than the
infimum $2+\alpha$ of $N_\psi(\lambda_0)$ over complex geodesics $\psi$
with $\psi(\lambda_0)=(0,0)$, $\psi'(\lambda_0)\in\C^*e_z$, although every boundary
point of the axis disc is strongly pseudoconvex. At each point
$(\zeta,0)$ of $\mathfrak c$ there is a unitary map $(z,w)\mapsto(z,e^{i\vartheta}w)$
fixing $(\zeta,0)$ which carries the tangential type, the Levi form and the
model polynomial of $\Omega_1$ at $(\zeta,0)$ to those of $\Omega_2$,
both defined with the same $\eps\in(0,1/8]$, but
the coefficients at the center are $2+\alpha$ and $\frac74+\alpha$. Thus
data at each point of the extremal circle that agree up to a pointwise
unitary equivalence (type, Levi form, model polynomial) do not
determine the coefficient.
\end{corollary}

\begin{proof}
Combine Theorem~\ref{thm:noncoercive-main} and
Proposition~\ref{prop:twisted-geometry}. At $(\zeta,0)$ the model
polynomial of $\Omega_\ell$ is $\eps(\operatorname{Im}(\bar\zeta^{\ell/2}w))^2$,
the Levi form is $\eps$ and the type is $2$. With
$e^{i\vartheta}=\bar\zeta^{1/2}$ one has
$\eps(\operatorname{Im}(\bar\zeta^{1/2}e^{i\vartheta}w))^2
=\eps(\operatorname{Im}(\bar\zeta w))^2$.
\end{proof}

\begin{remark}\label{rem:twisted}
For every $\ell$ the fiber exponent of the projection is $\frac34$: the
slice $\{z=\lambda\}$ has area $\asymp s^{1/4}\cdot s^{1/2}$, the factor
$s^{1/4}$ coming from the quartic term along the zero line of the
model. This exponent is invisible to the multitype ($T=1$). The
projection $z$ is the inverse which Lemma~\ref{lem:affine-inverse}
produces from $G=z$, since then $g_{(z,w)}(\lambda)=z-\lambda$, and
Lemma~\ref{lem:nontangential-fibers}(b) fails for it at every point of
$\mathfrak c$. For $\ell=1$ this deficiency disappears after the curved
correction $z-\frac\eps2w^2$, whose collar exponent is $1$; for
$\ell\ge2$ it is genuine, and by Lemma~\ref{lem:collar-competitor} every
$u\in\mathcal E((0,0),e_z)$ has $\beta(u)\le\frac34$. In the fixed axis coordinate and tangential trivialization, the
winding number distinguishing $\Omega_1$ from $\Omega_\ell$, $\ell\ge2$, is
$\operatorname{wind}\mathsf B=-\ell$ of the pure holomorphic tangential Hessian
along $\mathfrak c$, measured in the trivialization in which the left
inverse is $z$; the factor $\bar\zeta$ in~\eqref{eq:correctable}, which
compares the differential $dz$ of the left inverse with the outer
normal $(\zeta,0)$, shifts the admissible winding by one.

The coefficient is unchanged under a biholomorphism $\Xi:D\to D'$ such
that $\Xi$ and $\Xi^{-1}$ are $C^1$ up to the boundary: $f\mapsto f\circ\Xi$
is an isometry from $L^p(\Xi_*\mu)$ to $L^p(\mu)$ preserving jets, and
$\Xi_*\mu$ is comparable to $\delta_{D'}^\alpha dV$, so
Lemma~\ref{lem:basic-comparison} applies. Coercivity is not invariant in
this sense (Section~\ref{sec:regular-reverse}), and the polynomial
automorphism $(z,w)\mapsto(z-\frac\eps2w^2,w)$ of $\C^2$ carries the
curved fibers of the corrected inverse of $\Omega_1$ to affine ones.
Within class (F), correctability, which admits holomorphic rather than
only affine changes of the left inverse, is therefore the more natural
notion.
\end{remark}

\section{An optimized lower bound and a remaining question}\label{sec:questions}
For a normalized extremal $u:D\to\Delta$, $u(x)=0$,
$du_x(v)=\gamma_D(x;v)$, define $\beta_u(\zeta)$ as the supremum of
the $\beta\ge0$ for which, for every $M\ge2$, there are $c,s_0>0$
such that
\[
 u_*\mu\le c\,s^{\alpha+\beta}dA
       \quad\hbox{on }\mathcal A_M(\zeta)\cap\{s<s_0\}.
\]
The supremum of the empty set is defined to be zero. Constants may
depend on $\zeta,M,\beta$; attainment is not required.

\begin{proposition}[Optimizing the fiber lower bound]\label{prop:optimized-lower}
On a bounded convex domain with $C^2$ boundary, for a probability
measure $d\mu\asymp\delta_D^\alpha dV$, $\alpha>-1$,
$\beta_u$ is Borel measurable and
\[
 \liminf_{p\to\infty}\frac p{\log p}
       \log\frac{B_{p,\mu}(x;v)}{\gamma_D(x;v)}
 \ge1+\alpha+\sup_u\int_{\partial\Delta}\beta_u\,d\omega_0.
\]
\end{proposition}
\begin{proof}
For rational $r>0$ and integers $M,j\ge2$, let $E_{r,M,j}$ be the
closed set used in the proof of Lemma~\ref{lem:stolz-transfer}, with
$u_*\mu\le j s^{\alpha+r}dA$ on
$\mathcal A_M(\zeta)\cap\{s<1/j\}$. For $q\ge0$,
\[
 \{\beta_u>q\}=\bigcup_{r\in\mathbb Q,\ r>q}
                   \bigcap_{M\ge2}\bigcup_{j\ge2} E_{r,M,j}.
\]
Indeed, a bound at an exponent implies the same bound at each smaller
nonnegative exponent, and integer apertures suffice. This proves
measurability. On $C=\{\beta_u>0\}$ choose increasing, finite-valued
nonnegative simple functions $b_j<\beta_u$ with $b_j\uparrow\beta_u$.
Each $b_j(\zeta)$ is admissible, so the transfer lemma applied to
$C,b_j$ gives the bound with $\int_C b_j\,d\omega_0$. Monotone
convergence and then the supremum over $u$ prove the result.
Theorem~\ref{thm:general-main}(i) also implies
$\sup_u\int\beta_u\,d\omega_0\le n-1$.
\end{proof}

\emph{Conjecture.} On a smooth convex domain of finite type, the
coefficient at every jet exists and equals
$1+\alpha+\sup_u\int\beta_u\,d\omega_0$.

The proposition proves the lower bound without assuming that the
endpoint exponents are attained or finite-valued. Theorem~\ref{thm:regular-reverse}
proves equality for a.e.\ coercive geodesics. In general the missing
statement is a matching lower norm estimate for all competitors in terms of
the best pushforward exponents; constructing a single left inverse
does not supply it. The double-cover argument proves such an estimate
in the quartic negative-winding class, but does not cover arbitrary
non-coercive contact in higher dimensions.

The negative-winding question is settled by
Theorem~\ref{thm:winding-gap} for $C^4$ boundaries, and its $C^2$
version is false by Proposition~\ref{prop:C2-obstruction}. Corollary~\ref{cor:holder-gap} gives a strict gap under any positive
H\"older control of second derivatives. Non-power moduli of
continuity, and the coefficients for contact flatter than quartic
with negative winding, remain open here.

\medskip
\noindent\textbf{Acknowledgments.}
The authors used AI to help explore proof strategies and calculations.
The authors are solely responsible for verifying all arguments and
calculations and for the final content.

\defbibheading{nomark}{\section*{References}}
\renewcommand*{\bibfont}{\footnotesize}
\printbibliography[heading=nomark]
\end{document}